% File jytex.tex, for jyTeX version 2.6M (June 1992)
% Copyright (c) 1991, 1992 by Jonathan P. Yamron
% For full documentation, "get jydoc" from hep-ph@xxx.lanl.gov
%   Problems?  Contact brahm@theory3.caltech.edu.

\catcode`\@=11

%*****************************************************************************

\message{Loading jyTeX fonts...}

%************************************************************
%*
%*             Available fonts
%*
%************************************************************

%************** 5-point fonts *******************************

\font\vptrm=cmr5 \font\vptmit=cmmi5 \font\vptsy=cmsy5 \font\vptbf=cmbx5

\skewchar\vptmit='177 \skewchar\vptsy='60 \fontdimen16
\vptsy=\the\fontdimen17 \vptsy

\def\vpt{\ifmmode\err@badsizechange\else
     \@mathfontinit
     \textfont0=\vptrm  \scriptfont0=\vptrm  \scriptscriptfont0=\vptrm
     \textfont1=\vptmit \scriptfont1=\vptmit \scriptscriptfont1=\vptmit
     \textfont2=\vptsy  \scriptfont2=\vptsy  \scriptscriptfont2=\vptsy
     \textfont3=\xptex  \scriptfont3=\xptex  \scriptscriptfont3=\xptex
     \textfont\bffam=\vptbf
     \scriptfont\bffam=\vptbf
     \scriptscriptfont\bffam=\vptbf
     \@fontstyleinit
     \def\rm{\vptrm\fam=\z@}%
     \def\bf{\vptbf\fam=\bffam}%
     \def\oldstyle{\vptmit\fam=\@ne}%
     \rm\fi}

%************** 6-point fonts *******************************

\font\viptrm=cmr6 \font\viptmit=cmmi6 \font\viptsy=cmsy6
\font\viptbf=cmbx6

\skewchar\viptmit='177 \skewchar\viptsy='60 \fontdimen16
\viptsy=\the\fontdimen17 \viptsy

\def\vipt{\ifmmode\err@badsizechange\else
     \@mathfontinit
     \textfont0=\viptrm  \scriptfont0=\vptrm  \scriptscriptfont0=\vptrm
     \textfont1=\viptmit \scriptfont1=\vptmit \scriptscriptfont1=\vptmit
     \textfont2=\viptsy  \scriptfont2=\vptsy  \scriptscriptfont2=\vptsy
     \textfont3=\xptex   \scriptfont3=\xptex  \scriptscriptfont3=\xptex
     \textfont\bffam=\viptbf
     \scriptfont\bffam=\vptbf
     \scriptscriptfont\bffam=\vptbf
     \@fontstyleinit
     \def\rm{\viptrm\fam=\z@}%
     \def\bf{\viptbf\fam=\bffam}%
     \def\oldstyle{\viptmit\fam=\@ne}%
     \rm\fi}
%************** 7-point fonts *******************************

\font\viiptrm=cmr7 \font\viiptmit=cmmi7 \font\viiptsy=cmsy7
\font\viiptit=cmti7 \font\viiptbf=cmbx7

\skewchar\viiptmit='177 \skewchar\viiptsy='60 \fontdimen16
\viiptsy=\the\fontdimen17 \viiptsy

\def\viipt{\ifmmode\err@badsizechange\else
     \@mathfontinit
     \textfont0=\viiptrm  \scriptfont0=\vptrm  \scriptscriptfont0=\vptrm
     \textfont1=\viiptmit \scriptfont1=\vptmit \scriptscriptfont1=\vptmit
     \textfont2=\viiptsy  \scriptfont2=\vptsy  \scriptscriptfont2=\vptsy
     \textfont3=\xptex    \scriptfont3=\xptex  \scriptscriptfont3=\xptex
     \textfont\itfam=\viiptit
     \scriptfont\itfam=\viiptit
     \scriptscriptfont\itfam=\viiptit
     \textfont\bffam=\viiptbf
     \scriptfont\bffam=\vptbf
     \scriptscriptfont\bffam=\vptbf
     \@fontstyleinit
     \def\rm{\viiptrm\fam=\z@}%
     \def\it{\viiptit\fam=\itfam}%
     \def\bf{\viiptbf\fam=\bffam}%
     \def\oldstyle{\viiptmit\fam=\@ne}%
     \rm\fi}

%************** 8-point fonts *******************************

\font\viiiptrm=cmr8 \font\viiiptmit=cmmi8 \font\viiiptsy=cmsy8
\font\viiiptit=cmti8
%\font\viiiptsl=cmsl8
\font\viiiptbf=cmbx8
%\font\viiipttt=cmtt8
%\font\viiiptss=cmss8

\skewchar\viiiptmit='177 \skewchar\viiiptsy='60 \fontdimen16
\viiiptsy=\the\fontdimen17 \viiiptsy

\def\viiipt{\ifmmode\err@badsizechange\else
     \@mathfontinit
     \textfont0=\viiiptrm  \scriptfont0=\viptrm  \scriptscriptfont0=\vptrm
     \textfont1=\viiiptmit \scriptfont1=\viptmit \scriptscriptfont1=\vptmit
     \textfont2=\viiiptsy  \scriptfont2=\viptsy  \scriptscriptfont2=\vptsy
     \textfont3=\xptex     \scriptfont3=\xptex   \scriptscriptfont3=\xptex
     \textfont\itfam=\viiiptit
     \scriptfont\itfam=\viiptit
     \scriptscriptfont\itfam=\viiptit
     \textfont\bffam=\viiiptbf
     \scriptfont\bffam=\viptbf
     \scriptscriptfont\bffam=\vptbf
     \@fontstyleinit
     \def\rm{\viiiptrm\fam=\z@}%
     \def\it{\viiiptit\fam=\itfam}%
     \def\bf{\viiiptbf\fam=\bffam}%
     \def\oldstyle{\viiiptmit\fam=\@ne}%
     \rm\fi}

%************** Optional 9-point fonts **********************

\def\getixpt{%
     \font\ixptrm=cmr9
     \font\ixptmit=cmmi9
     \font\ixptsy=cmsy9
     \font\ixptit=cmti9
%     \font\ixptsl=cmsl9
     \font\ixptbf=cmbx9
%     \font\ixpttt=cmtt9
%     \font\ixptss=cmss9
     \skewchar\ixptmit='177 \skewchar\ixptsy='60
     \fontdimen16 \ixptsy=\the\fontdimen17 \ixptsy}

\def\ixpt{\ifmmode\err@badsizechange\else
     \@mathfontinit
     \textfont0=\ixptrm  \scriptfont0=\viiptrm  \scriptscriptfont0=\vptrm
     \textfont1=\ixptmit \scriptfont1=\viiptmit \scriptscriptfont1=\vptmit
     \textfont2=\ixptsy  \scriptfont2=\viiptsy  \scriptscriptfont2=\vptsy
     \textfont3=\xptex   \scriptfont3=\xptex    \scriptscriptfont3=\xptex
     \textfont\itfam=\ixptit
     \scriptfont\itfam=\viiptit
     \scriptscriptfont\itfam=\viiptit
     \textfont\bffam=\ixptbf
     \scriptfont\bffam=\viiptbf
     \scriptscriptfont\bffam=\vptbf
     \@fontstyleinit
     \def\rm{\ixptrm\fam=\z@}%
     \def\it{\ixptit\fam=\itfam}%
     \def\bf{\ixptbf\fam=\bffam}%
     \def\oldstyle{\ixptmit\fam=\@ne}%
     \rm\fi}

%************** 10-point fonts ******************************

\font\xptrm=cmr10 \font\xptmit=cmmi10 \font\xptsy=cmsy10
\font\xptex=cmex10 \font\xptit=cmti10 \font\xptsl=cmsl10
\font\xptbf=cmbx10 \font\xpttt=cmtt10 \font\xptss=cmss10
\font\xptsc=cmcsc10 \font\xptbfs=cmb10 \font\xptbmit=cmmib10

\skewchar\xptmit='177 \skewchar\xptbmit='177 \skewchar\xptsy='60
\fontdimen16 \xptsy=\the\fontdimen17 \xptsy

\def\xpt{\ifmmode\err@badsizechange\else
     \@mathfontinit
     \textfont0=\xptrm  \scriptfont0=\viiptrm  \scriptscriptfont0=\vptrm
     \textfont1=\xptmit \scriptfont1=\viiptmit \scriptscriptfont1=\vptmit
     \textfont2=\xptsy  \scriptfont2=\viiptsy  \scriptscriptfont2=\vptsy
     \textfont3=\xptex  \scriptfont3=\xptex    \scriptscriptfont3=\xptex
     \textfont\itfam=\xptit
     \scriptfont\itfam=\viiptit
     \scriptscriptfont\itfam=\viiptit
     \textfont\bffam=\xptbf
     \scriptfont\bffam=\viiptbf
     \scriptscriptfont\bffam=\vptbf
     \textfont\bfsfam=\xptbfs
     \scriptfont\bfsfam=\viiptbf
     \scriptscriptfont\bfsfam=\vptbf
     \textfont\bmitfam=\xptbmit
     \scriptfont\bmitfam=\viiptmit
     \scriptscriptfont\bmitfam=\vptmit
     \@fontstyleinit
     \def\rm{\xptrm\fam=\z@}%
     \def\it{\xptit\fam=\itfam}%
     \def\sl{\xptsl}%
     \def\bf{\xptbf\fam=\bffam}%
     \def\tt{\xpttt}%
     \def\ss{\xptss}%
     \def\sc{\xptsc}%
     \def\bfs{\xptbfs\fam=\bfsfam}%
     \def\bmit{\fam=\bmitfam}%
     \def\oldstyle{\xptmit\fam=\@ne}%
     \rm\fi}

%************** Optional 11-point fonts *********************

\def\getxipt{%
     \font\xiptrm=cmr10  scaled\magstephalf
     \font\xiptmit=cmmi10 scaled\magstephalf
     \font\xiptsy=cmsy10 scaled\magstephalf
     \font\xiptex=cmex10 scaled\magstephalf
     \font\xiptit=cmti10 scaled\magstephalf
     \font\xiptsl=cmsl10 scaled\magstephalf
     \font\xiptbf=cmbx10 scaled\magstephalf
     \font\xipttt=cmtt10 scaled\magstephalf
     \font\xiptss=cmss10 scaled\magstephalf
     \skewchar\xiptmit='177 \skewchar\xiptsy='60
     \fontdimen16 \xiptsy=\the\fontdimen17 \xiptsy}

\def\xipt{\ifmmode\err@badsizechange\else
     \@mathfontinit
     \textfont0=\xiptrm  \scriptfont0=\viiiptrm  \scriptscriptfont0=\viptrm
     \textfont1=\xiptmit \scriptfont1=\viiiptmit \scriptscriptfont1=\viptmit
     \textfont2=\xiptsy  \scriptfont2=\viiiptsy  \scriptscriptfont2=\viptsy
     \textfont3=\xiptex  \scriptfont3=\xptex     \scriptscriptfont3=\xptex
     \textfont\itfam=\xiptit
     \scriptfont\itfam=\viiiptit
     \scriptscriptfont\itfam=\viiptit
     \textfont\bffam=\xiptbf
     \scriptfont\bffam=\viiiptbf
     \scriptscriptfont\bffam=\viptbf
     \@fontstyleinit
     \def\rm{\xiptrm\fam=\z@}%
     \def\it{\xiptit\fam=\itfam}%
     \def\sl{\xiptsl}%
     \def\bf{\xiptbf\fam=\bffam}%
     \def\tt{\xipttt}%
     \def\ss{\xiptss}%
     \def\oldstyle{\xiptmit\fam=\@ne}%
     \rm\fi}

%************** 12-point fonts ******************************

\font\xiiptrm=cmr12 \font\xiiptmit=cmmi12 \font\xiiptsy=cmsy10
scaled\magstep1 \font\xiiptex=cmex10  scaled\magstep1 \font\xiiptit=cmti12
\font\xiiptsl=cmsl12 \font\xiiptbf=cmbx12
%\font\xiipttt=cmtt12
\font\xiiptss=cmss12 \font\xiiptsc=cmcsc10 scaled\magstep1
\font\xiiptbfs=cmb10  scaled\magstep1 \font\xiiptbmit=cmmib10
scaled\magstep1

\skewchar\xiiptmit='177 \skewchar\xiiptbmit='177 \skewchar\xiiptsy='60
\fontdimen16 \xiiptsy=\the\fontdimen17 \xiiptsy

\def\xiipt{\ifmmode\err@badsizechange\else
     \@mathfontinit
     \textfont0=\xiiptrm  \scriptfont0=\viiiptrm  \scriptscriptfont0=\viptrm
     \textfont1=\xiiptmit \scriptfont1=\viiiptmit \scriptscriptfont1=\viptmit
     \textfont2=\xiiptsy  \scriptfont2=\viiiptsy  \scriptscriptfont2=\viptsy
     \textfont3=\xiiptex  \scriptfont3=\xptex     \scriptscriptfont3=\xptex
     \textfont\itfam=\xiiptit
     \scriptfont\itfam=\viiiptit
     \scriptscriptfont\itfam=\viiptit
     \textfont\bffam=\xiiptbf
     \scriptfont\bffam=\viiiptbf
     \scriptscriptfont\bffam=\viptbf
     \textfont\bfsfam=\xiiptbfs
     \scriptfont\bfsfam=\viiiptbf
     \scriptscriptfont\bfsfam=\viptbf
     \textfont\bmitfam=\xiiptbmit
     \scriptfont\bmitfam=\viiiptmit
     \scriptscriptfont\bmitfam=\viptmit
     \@fontstyleinit
     \def\rm{\xiiptrm\fam=\z@}%
     \def\it{\xiiptit\fam=\itfam}%
     \def\sl{\xiiptsl}%
     \def\bf{\xiiptbf\fam=\bffam}%
     \def\tt{\xiipttt}%
     \def\ss{\xiiptss}%
     \def\sc{\xiiptsc}%
     \def\bfs{\xiiptbfs\fam=\bfsfam}%
     \def\bmit{\fam=\bmitfam}%
     \def\oldstyle{\xiiptmit\fam=\@ne}%
     \rm\fi}

%************** Optional 13-point fonts *********************

\def\getxiiipt{%
     \font\xiiiptrm=cmr12  scaled\magstephalf
     \font\xiiiptmit=cmmi12 scaled\magstephalf
     \font\xiiiptsy=cmsy9  scaled\magstep2
     \font\xiiiptit=cmti12 scaled\magstephalf
     \font\xiiiptsl=cmsl12 scaled\magstephalf
     \font\xiiiptbf=cmbx12 scaled\magstephalf
     \font\xiiipttt=cmtt12 scaled\magstephalf
     \font\xiiiptss=cmss12 scaled\magstephalf
     \skewchar\xiiiptmit='177 \skewchar\xiiiptsy='60
     \fontdimen16 \xiiiptsy=\the\fontdimen17 \xiiiptsy}

\def\xiiipt{\ifmmode\err@badsizechange\else
     \@mathfontinit
     \textfont0=\xiiiptrm  \scriptfont0=\xptrm  \scriptscriptfont0=\viiptrm
     \textfont1=\xiiiptmit \scriptfont1=\xptmit \scriptscriptfont1=\viiptmit
     \textfont2=\xiiiptsy  \scriptfont2=\xptsy  \scriptscriptfont2=\viiptsy
     \textfont3=\xivptex   \scriptfont3=\xptex  \scriptscriptfont3=\xptex
     \textfont\itfam=\xiiiptit
     \scriptfont\itfam=\xptit
     \scriptscriptfont\itfam=\viiptit
     \textfont\bffam=\xiiiptbf
     \scriptfont\bffam=\xptbf
     \scriptscriptfont\bffam=\viiptbf
     \@fontstyleinit
     \def\rm{\xiiiptrm\fam=\z@}%
     \def\it{\xiiiptit\fam=\itfam}%
     \def\sl{\xiiiptsl}%
     \def\bf{\xiiiptbf\fam=\bffam}%
     \def\tt{\xiiipttt}%
     \def\ss{\xiiiptss}%
     \def\oldstyle{\xiiiptmit\fam=\@ne}%
     \rm\fi}

%************** 14-point fonts ******************************

\font\xivptrm=cmr12   scaled\magstep1 \font\xivptmit=cmmi12
scaled\magstep1 \font\xivptsy=cmsy10  scaled\magstep2
\font\xivptex=cmex10  scaled\magstep2 \font\xivptit=cmti12 scaled\magstep1
\font\xivptsl=cmsl12  scaled\magstep1 \font\xivptbf=cmbx12
scaled\magstep1
%\font\xivpttt=cmtt12  scaled\magstep1
\font\xivptss=cmss12  scaled\magstep1 \font\xivptsc=cmcsc10
scaled\magstep2 \font\xivptbfs=cmb10  scaled\magstep2
\font\xivptbmit=cmmib10 scaled\magstep2

\skewchar\xivptmit='177 \skewchar\xivptbmit='177 \skewchar\xivptsy='60
\fontdimen16 \xivptsy=\the\fontdimen17 \xivptsy

\def\xivpt{\ifmmode\err@badsizechange\else
     \@mathfontinit
     \textfont0=\xivptrm  \scriptfont0=\xptrm  \scriptscriptfont0=\viiptrm
     \textfont1=\xivptmit \scriptfont1=\xptmit \scriptscriptfont1=\viiptmit
     \textfont2=\xivptsy  \scriptfont2=\xptsy  \scriptscriptfont2=\viiptsy
     \textfont3=\xivptex  \scriptfont3=\xptex  \scriptscriptfont3=\xptex
     \textfont\itfam=\xivptit
     \scriptfont\itfam=\xptit
     \scriptscriptfont\itfam=\viiptit
     \textfont\bffam=\xivptbf
     \scriptfont\bffam=\xptbf
     \scriptscriptfont\bffam=\viiptbf
     \textfont\bfsfam=\xivptbfs
     \scriptfont\bfsfam=\xptbfs
     \scriptscriptfont\bfsfam=\viiptbf
     \textfont\bmitfam=\xivptbmit
     \scriptfont\bmitfam=\xptbmit
     \scriptscriptfont\bmitfam=\viiptmit
     \@fontstyleinit
     \def\rm{\xivptrm\fam=\z@}%
     \def\it{\xivptit\fam=\itfam}%
     \def\sl{\xivptsl}%
     \def\bf{\xivptbf\fam=\bffam}%
     \def\tt{\xivpttt}%
     \def\ss{\xivptss}%
     \def\sc{\xivptsc}%
     \def\bfs{\xivptbfs\fam=\bfsfam}%
     \def\bmit{\fam=\bmitfam}%
     \def\oldstyle{\xivptmit\fam=\@ne}%
     \rm\fi}

%************** 17-point fonts ******************************

\font\xviiptrm=cmr17 \font\xviiptmit=cmmi12 scaled\magstep2
\font\xviiptsy=cmsy10 scaled\magstep3 \font\xviiptex=cmex10
scaled\magstep3 \font\xviiptit=cmti12 scaled\magstep2
\font\xviiptbf=cmbx12 scaled\magstep2 \font\xviiptbfs=cmb10
scaled\magstep3

\skewchar\xviiptmit='177 \skewchar\xviiptsy='60 \fontdimen16
\xviiptsy=\the\fontdimen17 \xviiptsy

\def\xviipt{\ifmmode\err@badsizechange\else
     \@mathfontinit
     \textfont0=\xviiptrm  \scriptfont0=\xiiptrm  \scriptscriptfont0=\viiiptrm
     \textfont1=\xviiptmit \scriptfont1=\xiiptmit \scriptscriptfont1=\viiiptmit
     \textfont2=\xviiptsy  \scriptfont2=\xiiptsy  \scriptscriptfont2=\viiiptsy
     \textfont3=\xviiptex  \scriptfont3=\xiiptex  \scriptscriptfont3=\xptex
     \textfont\itfam=\xviiptit
     \scriptfont\itfam=\xiiptit
     \scriptscriptfont\itfam=\viiiptit
     \textfont\bffam=\xviiptbf
     \scriptfont\bffam=\xiiptbf
     \scriptscriptfont\bffam=\viiiptbf
     \textfont\bfsfam=\xviiptbfs
     \scriptfont\bfsfam=\xiiptbfs
     \scriptscriptfont\bfsfam=\viiiptbf
     \@fontstyleinit
     \def\rm{\xviiptrm\fam=\z@}%
     \def\it{\xviiptit\fam=\itfam}%
     \def\bf{\xviiptbf\fam=\bffam}%
     \def\bfs{\xviiptbfs\fam=\bfsfam}%
     \def\oldstyle{\xviiptmit\fam=\@ne}%
     \rm\fi}

%************** 21-point fonts ******************************

\font\xxiptrm=cmr17  scaled\magstep1
%\font\xxiptmit=cmmi12 scaled\magstep3
%\font\xxiptsy=cmsy10 scaled\magstep4
%\font\xxiptex=cmex10 scaled\magstep4
%\font\xxiptbf=cmbx12 scaled\magstep3

%\skewchar\xxiptmit='177 \skewchar\xxiptsy='60
%\fontdimen16 \xxiptsy=\the\fontdimen17 \xxiptsy

\def\xxipt{\ifmmode\err@badsizechange\else
     \@mathfontinit
%     \textfont0=\xxiptrm  \scriptfont0=\xivptrm  \scriptscriptfont0=\xptrm
%     \textfont1=\xxiptmit \scriptfont1=\xivptmit \scriptscriptfont1=\xptmit
%     \textfont2=\xxiptsy  \scriptfont2=\xivptsy  \scriptscriptfont2=\xptsy
%     \textfont3=\xxiptex  \scriptfont3=\xivptex  \scriptscriptfont3=\xptex
%     \textfont\bffam=\xxiptbf
%     \scriptfont\bffam=\xivptbf
%     \scriptscriptfont\bffam=\xptbf
     \@fontstyleinit
     \def\rm{\xxiptrm\fam=\z@}%
     \rm\fi}

%************** 25-point fonts ******************************

\font\xxvptrm=cmr17  scaled\magstep2
%\font\xxvptmit=cmmi12 scaled\magstep4
%\font\xxvptsy=cmsy10 scaled\magstep5
%\font\xxvptex=cmex10 scaled\magstep5
%\font\xxvptbf=cmbx12 scaled\magstep4

%\skewchar\xxvptmit='177 \skewchar\xxvptsy='60
%\fontdimen16 \xxvptsy=\the\fontdimen17 \xxvptsy

\def\xxvpt{\ifmmode\err@badsizechange\else
     \@mathfontinit
%     \textfont0=\xxvptrm  \scriptfont0=\xviiptrm  \scriptscriptfont0=\xiiptrm
%     \textfont1=\xxvptmit \scriptfont1=\xviiptmit \scriptscriptfont1=\xiiptmit
%     \textfont2=\xxvptsy  \scriptfont2=\xviiptsy  \scriptscriptfont2=\xiiptsy
%     \textfont3=\xxvptex  \scriptfont3=\xviiptex  \scriptscriptfont3=\xiiptex
%     \textfont\bffam=\xxvptbf
%     \scriptfont\bffam=\xviiptbf
%     \scriptscriptfont\bffam=\xiiptbf
     \@fontstyleinit
     \def\rm{\xxvptrm\fam=\z@}%
     \rm\fi}

%************** Other fonts *********************************

%\font\dummy=dummy

%******************************************************************************

\message{Loading jyTeX macros...}

%************************************************************
%*
%*              Simple modifications to plain
%*
%************************************************************
\message{modifications to plain.tex,}

% The "\outer" qualifier is removed from the definitions of \newcount through
% \newif so that they may be used in definitions.  \newif is also changed to
% make \if commands globally defined.

\def\newcount{\alloc@0\count\countdef\insc@unt}
\def\newdimen{\alloc@1\dimen\dimendef\insc@unt}
\def\newskip{\alloc@2\skip\skipdef\insc@unt}
\def\newmuskip{\alloc@3\muskip\muskipdef\@cclvi}
\def\newbox{\alloc@4\box\chardef\insc@unt}
\def\newtoks{\alloc@5\toks\toksdef\@cclvi}
\def\newhelp#1#2{\newtoks#1\global#1\expandafter{\csname#2\endcsname}}
\def\newread{\alloc@6\read\chardef\sixt@@n}
\def\newwrite{\alloc@7\write\chardef\sixt@@n}
\def\newfam{\alloc@8\fam\chardef\sixt@@n}
\def\newinsert#1{\global\advance\insc@unt by\m@ne
     \ch@ck0\insc@unt\count
     \ch@ck1\insc@unt\dimen
     \ch@ck2\insc@unt\skip
     \ch@ck4\insc@unt\box
     \allocationnumber=\insc@unt
     \global\chardef#1=\allocationnumber
     \wlog{\string#1=\string\insert\the\allocationnumber}}
\def\newif#1{\count@\escapechar \escapechar\m@ne
     \expandafter\expandafter\expandafter
          \xdef\@if#1{true}{\let\noexpand#1=\noexpand\iftrue}%
     \expandafter\expandafter\expandafter
          \xdef\@if#1{false}{\let\noexpand#1=\noexpand\iffalse}%
     \global\@if#1{false}\escapechar=\count@}

%************** Some parameter changes **********************

\newlinechar=`\^^J
\overfullrule=0pt

%************** Font-related modifications ******************

% The plain fonts are mapped onto the corresponding jyTeX fonts

% Some control sequences are disabled.

\let\itfam=\undefined

\let\bffam=\undefined

\count18=3

% German sharp s is given a new name (\ss is already taken)

\chardef\sharps="19

% The mathcode assignments of characters in the math italic font are changed to
% allow for switching to boldface.

\mathchardef\alpha="710B \mathchardef\beta="710C
\mathchardef\gamma="710D \mathchardef\delta="710E
\mathchardef\epsilon="710F \mathchardef\zeta="7110
\mathchardef\eta="7111 \mathchardef\theta="7112 \mathchardef\iota="7113
\mathchardef\kappa="7114 \mathchardef\lambda="7115
\mathchardef\mu="7116 \mathchardef\nu="7117 \mathchardef\xi="7118
\mathchardef\pi="7119 \mathchardef\rho="711A \mathchardef\sigma="711B
\mathchardef\tau="711C \mathchardef\upsilon="711D
\mathchardef\phi="711E \mathchardef\chi="711F \mathchardef\psi="7120
\mathchardef\omega="7121 \mathchardef\varepsilon="7122
\mathchardef\vartheta="7123 \mathchardef\varpi="7124
\mathchardef\varrho="7125 \mathchardef\varsigma="7126
\mathchardef\varphi="7127 \mathchardef\imath="717B
\mathchardef\jmath="717C \mathchardef\ell="7160 \mathchardef\wp="717D
\mathchardef\partial="7140 \mathchardef\flat="715B
\mathchardef\natural="715C \mathchardef\sharp="715D

%************** Miscellaneous changes ***********************

% The dimension \p@ (1pt) is replaced with \rp@ (relative pt, defined below),
% whose size is determined by the base type size of the document.

\def\angle{{\vbox{\ialign{$\m@th\scriptstyle##$\crcr
     \not\mathrel{\mkern14mu}\crcr
     \noalign{\nointerlineskip}
     \mkern2.5mu\leaders\hrule height.34\rp@\hfill\mkern2.5mu\crcr}}}}
\def\vdots{\vbox{\baselineskip4\rp@ \lineskiplimit\z@
     \kern6\rp@\hbox{.}\hbox{.}\hbox{.}}}
\def\ddots{\mathinner{\mkern1mu\raise7\rp@\vbox{\kern7\rp@\hbox{.}}\mkern2mu
     \raise4\rp@\hbox{.}\mkern2mu\raise\rp@\hbox{.}\mkern1mu}}
\def\overrightarrow#1{\vbox{\ialign{##\crcr
     \rightarrowfill\crcr
     \noalign{\kern-\rp@\nointerlineskip}
     $\hfil\displaystyle{#1}\hfil$\crcr}}}
\def\overleftarrow#1{\vbox{\ialign{##\crcr
     \leftarrowfill\crcr
     \noalign{\kern-\rp@\nointerlineskip}
     $\hfil\displaystyle{#1}\hfil$\crcr}}}
\def\overbrace#1{\mathop{\vbox{\ialign{##\crcr
     \noalign{\kern3\rp@}
     \downbracefill\crcr
     \noalign{\kern3\rp@\nointerlineskip}
     $\hfil\displaystyle{#1}\hfil$\crcr}}}\limits}
\def\underbrace#1{\mathop{\vtop{\ialign{##\crcr
     $\hfil\displaystyle{#1}\hfil$\crcr
     \noalign{\kern3\rp@\nointerlineskip}
     \upbracefill\crcr
     \noalign{\kern3\rp@}}}}\limits}
\def\big#1{{\hbox{$\left#1\vbox to8.5\rp@ {}\right.\n@space$}}}
\def\Big#1{{\hbox{$\left#1\vbox to11.5\rp@ {}\right.\n@space$}}}
\def\bigg#1{{\hbox{$\left#1\vbox to14.5\rp@ {}\right.\n@space$}}}
\def\Bigg#1{{\hbox{$\left#1\vbox to17.5\rp@ {}\right.\n@space$}}}
\def\@vereq#1#2{\lower.5\rp@\vbox{\baselineskip\z@skip\lineskip-.5\rp@
     \ialign{$\m@th#1\hfil##\hfil$\crcr#2\crcr=\crcr}}}
\def\rlh@#1{\vcenter{\hbox{\ooalign{\raise2\rp@
     \hbox{$#1\rightharpoonup$}\crcr
     $#1\leftharpoondown$}}}}
\def\bordermatrix#1{\begingroup\m@th
     \setbox\z@\vbox{%
          \def\cr{\crcr\noalign{\kern2\rp@\global\let\cr\endline}}%
          \ialign{$##$\hfil\kern2\rp@\kern\p@renwd
               &\thinspace\hfil$##$\hfil&&\quad\hfil$##$\hfil\crcr
               \omit\strut\hfil\crcr
               \noalign{\kern-\baselineskip}%
               #1\crcr\omit\strut\cr}}%
     \setbox\tw@\vbox{\unvcopy\z@\global\setbox\@ne\lastbox}%
     \setbox\tw@\hbox{\unhbox\@ne\unskip\global\setbox\@ne\lastbox}%
     \setbox\tw@\hbox{$\kern\wd\@ne\kern-\p@renwd\left(\kern-\wd\@ne
          \global\setbox\@ne\vbox{\box\@ne\kern2\rp@}%
          \vcenter{\kern-\ht\@ne\unvbox\z@\kern-\baselineskip}%
          \,\right)$}%
     \null\;\vbox{\kern\ht\@ne\box\tw@}\endgroup}
\def\endinsert{\egroup
     \if@mid\dimen@\ht\z@
          \advance\dimen@\dp\z@
          \advance\dimen@12\rp@
          \advance\dimen@\pagetotal
          \ifdim\dimen@>\pagegoal\@midfalse\p@gefalse\fi
     \fi
     \if@mid\bigskip\box\z@
          \bigbreak
     \else\insert\topins{\penalty100 \splittopskip\z@skip
               \splitmaxdepth\maxdimen\floatingpenalty\z@
               \ifp@ge\dimen@\dp\z@
                    \vbox to\vsize{\unvbox\z@\kern-\dimen@}%
               \else\box\z@\nobreak\bigskip
               \fi}%
     \fi
     \endgroup}

% \normalbaselines is removed from \cases and \matrix.

\def\cases#1{\left\{\,\vcenter{\m@th
     \ialign{$##\hfil$&\quad##\hfil\crcr#1\crcr}}\right.}
\def\matrix#1{\null\,\vcenter{\m@th
     \ialign{\hfil$##$\hfil&&\quad\hfil$##$\hfil\crcr
          \mathstrut\crcr
          \noalign{\kern-\baselineskip}
          #1\crcr
          \mathstrut\crcr
          \noalign{\kern-\baselineskip}}}\,}

% \raggedbottom modified slightly

\newif\ifraggedbottom

\def\raggedbottom{\ifraggedbottom\else
     \advance\topskip by\z@ plus60pt \raggedbottomtrue\fi}%
\def\normalbottom{\ifraggedbottom
     \advance\topskip by\z@ plus-60pt \raggedbottomfalse\fi}

%************************************************************
%*
%*              Miscellaneous definitions
%*
%************************************************************
\message{hacks,}

%************** Hack registers ******************************

\toksdef\toks@i=1 \toksdef\toks@ii=2

%************** Basic macros ********************************

\def\TeX{T\kern-.1667em \lower.5ex \hbox{E}\kern-.125em X\null}
\def\jyTeX{{\leavevmode
     \raise.587ex \hbox{\it\j}\kern-.1em \lower.048ex \hbox{\it y}\kern-.12em
     \TeX}}

\let\then=\iftrue
\def\ifnoarg#1\then{\def\hack@{#1}\ifx\hack@\empty}
\def\ifundefined#1\then{%
     \expandafter\ifx\csname\expandafter\blank\string#1\endcsname\relax}
\def\useif#1\then{\csname#1\endcsname}
\def\usename#1{\csname#1\endcsname}
\def\useafter#1#2{\expandafter#1\csname#2\endcsname}

% Modify so that I can have \loop's within \loop's?
\long\def\loop#1\repeat{\def\@iterate{#1\expandafter\@iterate\fi}\@iterate
     \let\@iterate=\relax}
%\long\def\loop#1\repeat{\def\@loopbody{#1}\@iterate}
%\def\@iterate{\@loopbody\let\next=\@iterate\else\let\next=\relax\fi\next}

\let\TeXend=\end
\def\begin#1{\begingroup\def\@@blockname{#1}\usename{begin#1}}
\def\end#1{\usename{end#1}\def\hack@{#1}%
     \ifx\@@blockname\hack@
          \endgroup
     \else\err@badgroup\hack@\@@blockname
     \fi}
\def\@@blockname{}

\def\defaultoption[#1]#2{%
     \def\hack@{\ifx\hack@ii[\toks@={#2}\else\toks@={#2[#1]}\fi\the\toks@}%
     \futurelet\hack@ii\hack@}

\def\markup#1{\let\@@marksf=\empty
     \ifhmode\edef\@@marksf{\spacefactor=\the\spacefactor\relax}\/\fi
     ${}^{\hbox{\subscriptfonts#1}}$\@@marksf}

%************** Time registers ******************************

\newtoks\shortyear
\newtoks\militaryhour
\newtoks\standardhour
\newtoks\minute
\newtoks\amorpm

\def\settime{\count@=\time\divide\count@ by60
     \militaryhour=\expandafter{\number\count@}%
     {\multiply\count@ by-60 \advance\count@ by\time
          \xdef\hack@{\ifnum\count@<10 0\fi\number\count@}}%
     \minute=\expandafter{\hack@}%
     \ifnum\count@<12
          \amorpm={am}
     \else\amorpm={pm}
          \ifnum\count@>12 \advance\count@ by-12 \fi
     \fi
     \standardhour=\expandafter{\number\count@}%
     \def\hack@19##1##2{\shortyear={##1##2}}%
          \expandafter\hack@\the\year}

\def\monthword#1{%
     \ifcase#1
          $\bullet$\err@badcountervalue{monthword}%
          \or January\or February\or March\or April\or May\or June%
          \or July\or August\or September\or October\or November\or December%
     \else$\bullet$\err@badcountervalue{monthword}%
     \fi}

\def\monthabbr#1{%
     \ifcase#1
          $\bullet$\err@badcountervalue{monthabbr}%
          \or Jan\or Feb\or Mar\or Apr\or May\or Jun%
          \or Jul\or Aug\or Sep\or Oct\or Nov\or Dec%
     \else$\bullet$\err@badcountervalue{monthabbr}%
     \fi}

\def\militarytime{\the\militaryhour:\the\minute}
\def\standardtime{\the\standardhour:\the\minute}

%************** Number styles *******************************

\def\@setnumstyle#1#2{\expandafter\global\expandafter\expandafter
     \expandafter\let\expandafter\expandafter
     \csname @\expandafter\blank\string#1style\endcsname
     \csname#2\endcsname}
\def\numstyle#1{\usename{@\expandafter\blank\string#1style}#1}
\def\ifblank#1\then{\useafter\ifx{@\expandafter\blank\string#1}\blank}

\def\blank#1{}

\def\Roman#1{\expandafter\uppercase\expandafter{\romannumeral#1}}
\def\alphabetic#1{%
     \ifcase#1
          $\bullet$\err@badcountervalue{alphabetic}%
          \or a\or b\or c\or d\or e\or f\or g\or h\or i\or j\or k\or l\or m%
          \or n\or o\or p\or q\or r\or s\or t\or u\or v\or w\or x\or y\or z%
     \else$\bullet$\err@badcountervalue{alphabetic}%
     \fi}
\def\Alphabetic#1{\expandafter\uppercase\expandafter{\alphabetic{#1}}}
\def\symbols#1{%
     \ifcase#1
          $\bullet$\err@badcountervalue{symbols}%
          \or*\or\dag\or\ddag\or\S\or$\|$%
          \or**\or\dag\dag\or\ddag\ddag\or\S\S\or$\|\|$%
     \else$\bullet$\err@badcountervalue{symbols}%
     \fi}

%************** String macros *******************************

\catcode`\^^?=13 \def^^?{\relax}

\def\trimleading#1\to#2{\edef#2{#1}%
     \expandafter\@trimleading\expandafter#2#2^^?^^?}
\def\@trimleading#1#2#3^^?{\ifx#2^^?\def#1{}\else\def#1{#2#3}\fi}

\def\trimtrailing#1\to#2{\edef#2{#1}%
     \expandafter\@trimtrailing\expandafter#2#2^^? ^^?\relax}
\def\@trimtrailing#1#2 ^^?#3{\ifx#3\relax\toks@={}%
     \else\def#1{#2}\toks@={\trimtrailing#1\to#1}\fi
     \the\toks@}

\def\trim#1\to#2{\trimleading#1\to#2\trimtrailing#2\to#2}

\catcode`\^^?=15

%************** List macros *********************************

\long\def\additemL#1\to#2{\toks@={\^^\{#1}}\toks@ii=\expandafter{#2}%
     \xdef#2{\the\toks@\the\toks@ii}}

\long\def\additemR#1\to#2{\toks@={\^^\{#1}}\toks@ii=\expandafter{#2}%
     \xdef#2{\the\toks@ii\the\toks@}}

\def\getitemL#1\to#2{\expandafter\@getitemL#1\hack@#1#2}
\def\@getitemL\^^\#1#2\hack@#3#4{\def#4{#1}\def#3{#2}}

%************************************************************
%*
%*             Font-related macros
%*
%************************************************************
\message{font macros,}

%************** Font set-up *********************************

\newdimen\rp@
\newcount\@@sizeindex \@@sizeindex=0
\newcount\@@factori
\newcount\@@factorii
\newcount\@@factoriii
\newcount\@@factoriv

\countdef\maxfam=18
\newfam\itfam
\newfam\bffam
\newfam\bfsfam
\newfam\bmitfam

\def\@mathfontinit{\count@=4
     \loop\textfont\count@=\nullfont
          \scriptfont\count@=\nullfont
          \scriptscriptfont\count@=\nullfont
          \ifnum\count@<\maxfam\advance\count@ by\@ne
     \repeat}

\def\@fontstyleinit{%
     \def\it{\err@fontnotavailable\it}%
     \def\bf{\err@fontnotavailable\bf}%
     \def\bfs{\err@bfstobf}%
     \def\bmit{\err@fontnotavailable\bmit}%
     \def\sc{\err@fontnotavailable\sc}%
     \def\sl{\err@sltoit}%
     \def\ss{\err@fontnotavailable\ss}%
     \def\tt{\err@fontnotavailable\tt}}

\def\@parameterinit#1{\rm\rp@=.1em \@getscaling{#1}%
     \let\^^\=\@doscaling\scalingskipslist
     \setbox\strutbox=\hbox{\vrule
          height.708\baselineskip depth.292\baselineskip width\z@}}

\def\@getfactor#1#2#3#4{\@@factori=#1 \@@factorii=#2
     \@@factoriii=#3 \@@factoriv=#4}

\def\@getscaling#1{\count@=#1 \advance\count@ by-\@@sizeindex\@@sizeindex=#1
     \ifnum\count@<0
          \let\@mulordiv=\divide
          \let\@divormul=\multiply
          \multiply\count@ by\m@ne
     \else\let\@mulordiv=\multiply
          \let\@divormul=\divide
     \fi
     \edef\@@scratcha{\ifcase\count@                {1}{1}{1}{1}\or
          {1}{7}{23}{3}\or     {2}{5}{3}{1}\or      {9}{89}{13}{1}\or
          {6}{25}{6}{1}\or     {8}{71}{14}{1}\or    {6}{25}{36}{5}\or
          {1}{7}{53}{4}\or     {12}{125}{108}{5}\or {3}{14}{53}{5}\or
          {6}{41}{17}{1}\or    {13}{31}{13}{2}\or   {9}{107}{71}{2}\or
          {11}{139}{124}{3}\or {1}{6}{43}{2}\or     {10}{107}{42}{1}\or
          {1}{5}{43}{2}\or     {5}{69}{65}{1}\or    {11}{97}{91}{2}\fi}%
     \expandafter\@getfactor\@@scratcha}

\def\@doscaling#1{\@mulordiv#1by\@@factori\@divormul#1by\@@factorii
     \@mulordiv#1by\@@factoriii\@divormul#1by\@@factoriv}

%************* Size-changing commands ***********************

\newskip\headskip
\newskip\footskip

\def\typesize=#1pt{\count@=#1 \advance\count@ by-10
     \ifcase\count@
          \@setsizex\or\err@badtypesize\or
          \@setsizexii\or\err@badtypesize\or
          \@setsizexiv
     \else\err@badtypesize
     \fi}

\def\@setsizex{\getixpt
     \def\subsubscriptfonts{\vpt}%
          \def\subsubscriptsize{\vpt\@parameterinit{-8}}%
     \def\subscriptfonts{\viipt}\def\subscriptsize{\viipt\@parameterinit{-4}}%
     \def\footnotefonts{\viiipt}\def\footnotesize{\viiipt\@parameterinit{-2}}%
     \def\smallfonts{\ixpt}\def\smallsize{\ixpt\@parameterinit{-1}}%
     \def\normalfonts{\xpt}\def\normalsize{\xpt\@parameterinit{0}}%
     \def\bigfonts{\xiipt}\def\bigsize{\xiipt\@parameterinit{2}}%
     \def\Bigfonts{\xivpt}\def\Bigsize{\xivpt\@parameterinit{4}}%
     \def\biggfonts{\xviipt}\def\biggsize{\xviipt\@parameterinit{6}}%
     \def\Biggfonts{\xxipt}\def\Biggsize{\xxipt\@parameterinit{8}}%
     \def\tinyfonts{\vpt}\def\tinysize{\vpt\@parameterinit{-8}}%
     \def\HUGEFONTS{\xxvpt}\def\HUGESIZE{\xxvpt\@parameterinit{10}}%
     \normalsize\fixedskipslist}

\def\@setsizexii{\getxipt
     \def\subsubscriptfonts{\vipt}%
          \def\subsubscriptsize{\vipt\@parameterinit{-6}}%
     \def\subscriptfonts{\viiipt}%
          \def\subscriptsize{\viiipt\@parameterinit{-2}}%
     \def\footnotefonts{\xpt}\def\footnotesize{\xpt\@parameterinit{0}}%
     \def\smallfonts{\xipt}\def\smallsize{\xipt\@parameterinit{1}}%
     \def\normalfonts{\xiipt}\def\normalsize{\xiipt\@parameterinit{2}}%
     \def\bigfonts{\xivpt}\def\bigsize{\xivpt\@parameterinit{4}}%
     \def\Bigfonts{\xviipt}\def\Bigsize{\xviipt\@parameterinit{6}}%
     \def\biggfonts{\xxipt}\def\biggsize{\xxipt\@parameterinit{8}}%
     \def\Biggfonts{\xxvpt}\def\Biggsize{\xxvpt\@parameterinit{10}}%
     \def\tinyfonts{\vpt}\def\tinysize{\vpt\@parameterinit{-8}}%
     \def\HUGEFONTS{\xxvpt}\def\HUGESIZE{\xxvpt\@parameterinit{10}}%
     \normalsize\fixedskipslist}

\def\@setsizexiv{\getxiiipt
     \def\subsubscriptfonts{\viipt}%
          \def\subsubscriptsize{\viipt\@parameterinit{-4}}%
     \def\subscriptfonts{\xpt}\def\subscriptsize{\xpt\@parameterinit{0}}%
     \def\footnotefonts{\xiipt}\def\footnotesize{\xiipt\@parameterinit{2}}%
     \def\smallfonts{\xiiipt}\def\smallsize{\xiiipt\@parameterinit{3}}%
     \def\normalfonts{\xivpt}\def\normalsize{\xivpt\@parameterinit{4}}%
     \def\bigfonts{\xviipt}\def\bigsize{\xviipt\@parameterinit{6}}%
     \def\Bigfonts{\xxipt}\def\Bigsize{\xxipt\@parameterinit{8}}%
     \def\biggfonts{\xxvpt}\def\biggsize{\xxvpt\@parameterinit{10}}%
     \def\Biggfonts{\err@sizetoolarge\Biggfonts\HUGEFONTS}%
          \def\Biggsize{\err@sizetoolarge\Biggsize\HUGESIZE}%
     \def\tinyfonts{\vpt}\def\tinysize{\vpt\@parameterinit{-8}}%
     \def\HUGEFONTS{\xxvpt}\def\HUGESIZE{\xxvpt\@parameterinit{10}}%
     \normalsize\fixedskipslist}

\def\subsubscriptfonts{\vpt} \def\subsubscriptsize{\vpt\@parameterinit{-8}}
\def\subscriptfonts{\viipt}  \def\subscriptsize{\viipt\@parameterinit{-4}}
\def\footnotefonts{\viiipt}  \def\footnotesize{\viiipt\@parameterinit{-2}}
\def\smallfonts{\err@sizenotavailable\smallfonts}
                             \def\smallsize{\ixpt\@parameterinit{-1}}
\def\normalfonts{\xpt}       \def\normalsize{\xpt\@parameterinit{0}}
\def\bigfonts{\xiipt}        \def\bigsize{\xiipt\@parameterinit{2}}
\def\Bigfonts{\xivpt}        \def\Bigsize{\xivpt\@parameterinit{4}}
\def\biggfonts{\xviipt}      \def\biggsize{\xviipt\@parameterinit{6}}
\def\Biggfonts{\xxipt}       \def\Biggsize{\xxipt\@parameterinit{8}}
\def\tinyfonts{\vpt}         \def\tinysize{\vpt\@parameterinit{-8}}
\def\HUGEFONTS{\xxvpt}       \def\HUGESIZE{\xxvpt\@parameterinit{10}}

%************************************************************
%*
%*             Document layout
%*
%************************************************************
\message{document layout,}

%************** Page format *********************************

\newtoks\everyoutput \everyoutput={}
\newdimen\depthofpage
\newcount\pagenum \pagenum=0

\newdimen\oddtopmargin  \newdimen\eventopmargin
\newdimen\oddleftmargin \newdimen\evenleftmargin
\newtoks\oddhead        \newtoks\evenhead
\newtoks\oddfoot        \newtoks\evenfoot

\def\topmargin{\afterassignment\@seteventop\oddtopmargin}
\def\leftmargin{\afterassignment\@setevenleft\oddleftmargin}
\def\head{\afterassignment\@setevenhead\oddhead}
\def\foot{\afterassignment\@setevenfoot\oddfoot}

\def\@seteventop{\eventopmargin=\oddtopmargin}
\def\@setevenleft{\evenleftmargin=\oddleftmargin}
\def\@setevenhead{\evenhead=\oddhead}
\def\@setevenfoot{\evenfoot=\oddfoot}

\def\pagenumstyle#1{\@setnumstyle\pagenum{#1}}

\newif\ifdraft
\def\draft{\drafttrue\leftmargin=.5in \overfullrule=5pt }

\def\outputstyle#1{\global\expandafter\let\expandafter
          \@outputstyle\csname#1output\endcsname
     \usename{#1setup}}

\output={\@outputstyle}

\def\normaloutput{\the\everyoutput
     \global\advance\pagenum by\@ne
     \ifodd\pagenum
          \voffset=\oddtopmargin \hoffset=\oddleftmargin
     \else\voffset=\eventopmargin \hoffset=\evenleftmargin
     \fi
     \advance\voffset by-1in  \advance\hoffset by-1in
     \count0=\pagenum
     \expandafter\shipout\pagebox
     \ifnum\outputpenalty>-\@MM\else\dosupereject\fi}

\newdimen\fullhsize
\newbox\leftpage
\newcount\leftpagenum
\newcount\outputpagenum \outputpagenum=0
\let\leftorright=L

\def\twoupoutput{\the\everyoutput
     \global\advance\pagenum by\@ne
     \if L\leftorright
          \global\setbox\leftpage=\leftline{\pagebox}%
          \global\leftpagenum=\pagenum
          \global\let\leftorright=R%
     \else\global\advance\outputpagenum by\@ne
          \ifodd\outputpagenum
               \voffset=\oddtopmargin \hoffset=\oddleftmargin
          \else\voffset=\eventopmargin \hoffset=\evenleftmargin
          \fi
          \advance\voffset by-1in  \advance\hoffset by-1in
          \count0=\leftpagenum \count1=\pagenum
          \shipout\vbox{\hbox to\fullhsize
               {\box\leftpage\hfil\leftline{\pagebox}}}%
          \global\let\leftorright=L%
     \fi
     \ifnum\outputpenalty>-\@MM
     \else\dosupereject
          \if R\leftorright
               \globaldefs=\@ne\head={\hfil}\foot={\hfil}\globaldefs=\z@
               \null\newpage
          \fi
     \fi}

\def\pagebox{\vbox{\makeheadline\pagebody\makefootline}}

\def\makeheadline{%
     \vbox to\z@{\baselinestretch=\@m
          \vskip\topskip\vskip-.708\baselineskip\vskip-\headskip
          \line{\vbox to\ht\strutbox{}%
               \ifodd\pagenum\the\oddhead\else\the\evenhead\fi}%
          \vss}%
     \nointerlineskip}

\def\pagebody{\vbox to\vsize{%
     \boxmaxdepth\maxdepth
     \ifvoid\topins\else\unvbox\topins\fi
     \depthofpage=\dp255
     \unvbox255
     \ifraggedbottom\kern-\depthofpage\vfil\fi
     \ifvoid\footins
     \else\vskip\skip\footins
          \footnoterule
          \unvbox\footins
          \vskip-\footnoteskip
     \fi}}

\def\makefootline{\baselineskip=\footskip
     \line{\ifodd\pagenum\the\oddfoot\else\the\evenfoot\fi}}

%************** Sectioning commands *************************

\newskip\abovechapterskip
\newskip\belowchapterskip
\newskip\abovesectionskip
\newskip\belowsectionskip
\newskip\abovesubsectionskip
\newskip\belowsubsectionskip

\def\chapterstyle#1{\global\expandafter\let\expandafter\@chapterstyle
     \csname#1text\endcsname}
\def\sectionstyle#1{\global\expandafter\let\expandafter\@sectionstyle
     \csname#1text\endcsname}
\def\subsectionstyle#1{\global\expandafter\let\expandafter\@subsectionstyle
     \csname#1text\endcsname}

\def\chapter#1{%
     \ifdim\lastskip=17sp \else\chapterbreak\vskip\abovechapterskip\fi
     \@chapterstyle{\ifblank\chapternumstyle\then
          \else\newchapternum=\next\chapternumformat\ \fi#1}%
     \nobreak\vskip\belowchapterskip\vskip17sp }

\def\section#1{%
     \ifdim\lastskip=17sp \else\sectionbreak\vskip\abovesectionskip\fi
     \@sectionstyle{\ifblank\sectionnumstyle\then
          \else\newsectionnum=\next\sectionnumformat\ \fi#1}%
     \nobreak\vskip\belowsectionskip\vskip17sp }

\def\subsection#1{%
     \ifdim\lastskip=17sp \else\subsectionbreak\vskip\abovesubsectionskip\fi
     \@subsectionstyle{\ifblank\subsectionnumstyle\then
          \else\newsubsectionnum=\next\subsectionnumformat\ \fi#1}%
     \nobreak\vskip\belowsubsectionskip\vskip17sp }

%************** Text formatting commands ********************

\let\TeXunderline=\underline
\let\TeXoverline=\overline
\def\underline#1{\relax\ifmmode\TeXunderline{#1}\else
     $\TeXunderline{\hbox{#1}}$\fi}
\def\overline#1{\relax\ifmmode\TeXoverline{#1}\else
     $\TeXoverline{\hbox{#1}}$\fi}

\def\baselinestretch{\afterassignment\@baselinestretch\count@}
\def\@baselinestretch{\baselineskip=\normalbaselineskip
     \divide\baselineskip by\@m\baselineskip=\count@\baselineskip
     \setbox\strutbox=\hbox{\vrule
          height.708\baselineskip depth.292\baselineskip width\z@}%
     \bigskipamount=\the\baselineskip
          plus.25\baselineskip minus.25\baselineskip
     \medskipamount=.5\baselineskip
          plus.125\baselineskip minus.125\baselineskip
     \smallskipamount=.25\baselineskip
          plus.0625\baselineskip minus.0625\baselineskip}

\def\\{\ifhmode\ifnum\lastpenalty=-\@M\else\hfil\penalty-\@M\fi\fi
     \ignorespaces}
\def\newpage{\vfil\break}

\def\lefttext#1{\par{\@text\leftskip=\z@\rightskip=\centering
     \noindent#1\par}}
\def\righttext#1{\par{\@text\leftskip=\centering\rightskip=\z@
     \noindent#1\par}}
\def\centertext#1{\par{\@text\leftskip=\centering\rightskip=\centering
     \noindent#1\par}}
\def\@text{\parindent=\z@ \parfillskip=\z@ \everypar={}%
     \spaceskip=.3333em \xspaceskip=.5em
     \def\\{\ifhmode\ifnum\lastpenalty=-\@M\else\penalty-\@M\fi\fi
          \ignorespaces}}

\def\beginleft{\par\@text\leftskip=\z@ \rightskip=\centering}
     
\def\beginright{\par\@text\leftskip=\centering\rightskip=\z@ }
     
\def\begincenter{\par\@text\leftskip=\centering\rightskip=\centering}

\def\beginnarrow{\defaultoption[\parindent]\@beginnarrow}
\def\@beginnarrow[#1]{\par\advance\leftskip by#1\advance\rightskip by#1}

\begingroup
\catcode`\[=1 \catcode`\{=11 \gdef\beginignore[\endgroup\bgroup
     \catcode`\e=0 \catcode`\\=12 \catcode`\{=11 \catcode`\f=12 \let\or=\relax
     \let\nd{ignor=\fi \let\}=\egroup
     \iffalse}
\endgroup

\long\def\marginnote#1{\leavevmode
     \edef\@marginsf{\spacefactor=\the\spacefactor\relax}%
     \ifdraft\strut\vadjust{%
          \hbox to\z@{\hskip\hsize\hskip.1in
               \vbox to\z@{\vskip-\dp\strutbox
                    \marginnoteformat
                    \vskip-\ht\strutbox
                    \noindent\strut#1\par
                    \vss}%
               \hss}}%
     \fi
     \@marginsf}

%************** The \bye command ****************************

\newtoks\everybye \everybye={\par\vfil}
\outer\def\bye{\the\everybye
     \footnotecheck
     \prelabelcheck
     \streamcheck
     \supereject
     \TeXend}

%************************************************************
%*
%*             Footnotes
%*
%************************************************************
\message{footnotes,}

\newcount\footnotenum \footnotenum=0
\newskip\footnoteskip
\let\@footnotelist=\empty

\def\footnotenumstyle#1{\@setnumstyle\footnotenum{#1}%
     \useafter\ifx{@footnotenumstyle}\symbols
          \global\let\@footup=\empty
     \else\global\let\@footup=\markup
     \fi}

\def\footnote{\footnotecheck\defaultoption[]\@footnote}
\def\@footnote[#1]{\@footnotemark[#1]\@footnotetext}

\def\footnotemark{\defaultoption[]\@footnotemark}
\def\@footnotemark[#1]{\let\@footsf=\empty
     \ifhmode\edef\@footsf{\spacefactor=\the\spacefactor\relax}\/\fi
     \ifnoarg#1\then
          \global\advance\footnotenum by\@ne
          \@footup{\footnotenumformat}%
          \edef\@@foota{\footnotenum=\the\footnotenum\relax}%
          \expandafter\additemR\expandafter\@footup\expandafter
               {\@@foota\footnotenumformat}\to\@footnotelist
          \global\let\@footnotelist=\@footnotelist
     \else\markup{#1}%
          \additemR\markup{#1}\to\@footnotelist
          \global\let\@footnotelist=\@footnotelist
     \fi
     \@footsf}

\def\footnotetext{%
     \ifx\@footnotelist\empty\err@extrafootnotetext\else\@footnotetext\fi}
\def\@footnotetext{%
     \getitemL\@footnotelist\to\@@foota
     \global\let\@footnotelist=\@footnotelist
     \insert\footins\bgroup
     \footnoteformat
     \splittopskip=\ht\strutbox\splitmaxdepth=\dp\strutbox
     \interlinepenalty=\interfootnotelinepenalty\floatingpenalty=\@MM
     \noindent\llap{\@@foota}\strut
     \bgroup\aftergroup\@footnoteend
     \let\@@scratcha=}
\def\@footnoteend{\strut\par\vskip\footnoteskip\egroup}

\def\footnoterule{\normalfonts
     \kern-.3em \hrule width2in height.04em \kern .26em }

\def\footnotecheck{%
     \ifx\@footnotelist\empty
     \else\err@extrafootnotemark
          \global\let\@footnotelist=\empty
     \fi}

%************************************************************
%*
%*             Labelling macros
%*
%************************************************************
\message{labels,}

\let\@@labeldef=\xdef
\newif\if@labelfile
\newwrite\@labelfile
\let\@prelabellist=\empty

\def\label#1#2{\trim#1\to\@@labarg\edef\@@labtext{#2}%
     \edef\@@labname{lab@\@@labarg}%
     \useafter\ifundefined\@@labname\then\else\@yeslab\fi
     \useafter\@@labeldef\@@labname{#2}%
     \ifstreaming
          \expandafter\toks@\expandafter\expandafter\expandafter
               {\csname\@@labname\endcsname}%
          \immediate\write\streamout{\noexpand\label{\@@labarg}{\the\toks@}}%
     \fi}
\def\@yeslab{%
     \useafter\ifundefined{if\@@labname}\then
          \err@labelredef\@@labarg
     \else\useif{if\@@labname}\then
               \err@labelredef\@@labarg
          \else\global\usename{\@@labname true}%
               \useafter\ifundefined{pre\@@labname}\then
               \else\useafter\ifx{pre\@@labname}\@@labtext
                    \else\err@badlabelmatch\@@labarg
                    \fi
               \fi
               \if@labelfile
               \else\global\@labelfiletrue
                    \immediate\write\sixt@@n{--> Creating file \jobname.lab}%
                    \immediate\openout\@labelfile=\jobname.lab
               \fi
               \immediate\write\@labelfile
                    {\noexpand\prelabel{\@@labarg}{\@@labtext}}%
          \fi
     \fi}

\def\putlab#1{\trim#1\to\@@labarg\edef\@@labname{lab@\@@labarg}%
     \useafter\ifundefined\@@labname\then\@nolab\else\usename\@@labname\fi}
\def\@nolab{%
     \useafter\ifundefined{pre\@@labname}\then
          \undefinedlabelformat
          \err@needlabel\@@labarg
          \useafter\xdef\@@labname{\undefinedlabelformat}%
     \else\usename{pre\@@labname}%
          \useafter\xdef\@@labname{\usename{pre\@@labname}}%
     \fi
     \useafter\newif{if\@@labname}%
     \expandafter\additemR\@@labarg\to\@prelabellist}

\def\prelabel#1{\useafter\gdef{prelab@#1}}

\def\ifundefinedlabel#1\then{%
     \expandafter\ifx\csname lab@#1\endcsname\relax}
\def\useiflab#1\then{\csname iflab@#1\endcsname}

\def\prelabelcheck{{%
     \def\^^\##1{\useiflab{##1}\then\else\err@undefinedlabel{##1}\fi}%
     \@prelabellist}}

%************************************************************
%*
%*             Equation numbering
%*
%************************************************************
\message{equation numbering,}

\newcount\chapternum
\newcount\sectionnum
\newcount\subsectionnum
\newcount\equationnum
\newcount\subequationnum
\newcount\figurenum
\newcount\subfigurenum
\newcount\tablenum
\newcount\subtablenum

\newif\if@subeqncount
\newif\if@subfigcount
\newif\if@subtblcount

\def\newchapternum{\newsectionnum=\z@\@resetnum\chapternum}
\def\newsectionnum{\newsubsectionnum=\z@\@resetnum\sectionnum}
\def\newsubsectionnum{\newequationnum=\z@\newfigurenum=\z@\newtablenum=\z@
     \@resetnum\subsectionnum}
\def\newequationnum{\newsubequationnum=\z@\@resetnum\equationnum}
\def\newsubequationnum{\@resetnum\subequationnum}
\def\newfigurenum{\newsubfigurenum=\z@\@resetnum\figurenum}
\def\newsubfigurenum{\@resetnum\subfigurenum}
\def\newtablenum{\newsubtablenum=\z@\@resetnum\tablenum}
\def\newsubtablenum{\@resetnum\subtablenum}

\def\@resetnum#1{\global\advance#1by1 \edef\next{\the#1\relax}\global#1}

\newchapternum=0

\def\chapternumstyle#1{\@setnumstyle\chapternum{#1}}
\def\sectionnumstyle#1{\@setnumstyle\sectionnum{#1}}
\def\subsectionnumstyle#1{\@setnumstyle\subsectionnum{#1}}
\def\equationnumstyle#1{\@setnumstyle\equationnum{#1}}
\def\subequationnumstyle#1{\@setnumstyle\subequationnum{#1}%
     \ifblank\subequationnumstyle\then\global\@subeqncountfalse\fi
     \ignorespaces}
\def\figurenumstyle#1{\@setnumstyle\figurenum{#1}}
\def\subfigurenumstyle#1{\@setnumstyle\subfigurenum{#1}%
     \ifblank\subfigurenumstyle\then\global\@subfigcountfalse\fi
     \ignorespaces}
\def\tablenumstyle#1{\@setnumstyle\tablenum{#1}}
\def\subtablenumstyle#1{\@setnumstyle\subtablenum{#1}%
     \ifblank\subtablenumstyle\then\global\@subtblcountfalse\fi
     \ignorespaces}

\def\eqnlabel#1{%
     \if@subeqncount
          \newsubequationnum=\next
     \else\newequationnum=\next
          \ifblank\subequationnumstyle\then
          \else\global\@subeqncounttrue
               \newsubequationnum=\@ne
          \fi
     \fi
     \label{#1}{\puteqnformat}(\puteqn{#1})%
     \ifdraft\rlap{\hskip.1in{\tt#1}}\fi}

\let\puteqn=\putlab

\def\equation#1#2{\useafter\gdef{eqn@#1}{#2\eqno\eqnlabel{#1}}}
\def\Equation#1{\useafter\gdef{eqn@#1}}

\def\putequation#1{\useafter\ifundefined{eqn@#1}\then
     \err@undefinedeqn{#1}\else\usename{eqn@#1}\fi}

\def\eqnseriesstyle#1{\gdef\@eqnseriesstyle{#1}}
\def\begineqnseries{\subequationnumstyle{\@eqnseriesstyle}%
     \defaultoption[]\@begineqnseries}
\def\@begineqnseries[#1]{\edef\@@eqnname{#1}}
\def\endeqnseries{\subequationnumstyle{blank}%
     \expandafter\ifnoarg\@@eqnname\then
     \else\label\@@eqnname{\puteqnformat}%
     \fi
     \aftergroup\ignorespaces}

\def\figlabel#1{%
     \if@subfigcount
          \newsubfigurenum=\next
     \else\newfigurenum=\next
          \ifblank\subfigurenumstyle\then
          \else\global\@subfigcounttrue
               \newsubfigurenum=\@ne
          \fi
     \fi
     \label{#1}{\putfigformat}\putfig{#1}%
     {\def\marginnoteformat{\tt}\marginnote{#1}}}

\let\putfig=\putlab

\def\figseriesstyle#1{\gdef\@figseriesstyle{#1}}
\def\beginfigseries{\subfigurenumstyle{\@figseriesstyle}%
     \defaultoption[]\@beginfigseries}
\def\@beginfigseries[#1]{\edef\@@figname{#1}}
\def\endfigseries{\subfigurenumstyle{blank}%
     \expandafter\ifnoarg\@@figname\then
     \else\label\@@figname{\putfigformat}%
     \fi
     \aftergroup\ignorespaces}

\def\tbllabel#1{%
     \if@subtblcount
          \newsubtablenum=\next
     \else\newtablenum=\next
          \ifblank\subtablenumstyle\then
          \else\global\@subtblcounttrue
               \newsubtablenum=\@ne
          \fi
     \fi
     \label{#1}{\puttblformat}\puttbl{#1}%
     {\def\marginnoteformat{\tt}\marginnote{#1}}}

\let\puttbl=\putlab

\def\tblseriesstyle#1{\gdef\@tblseriesstyle{#1}}
\def\begintblseries{\subtablenumstyle{\@tblseriesstyle}%
     \defaultoption[]\@begintblseries}
\def\@begintblseries[#1]{\edef\@@tblname{#1}}
\def\endtblseries{\subtablenumstyle{blank}%
     \expandafter\ifnoarg\@@tblname\then
     \else\label\@@tblname{\puttblformat}%
     \fi
     \aftergroup\ignorespaces}

%************************************************************
%*
%*             Reference numbering
%*
%************************************************************
\message{reference numbering,}

\newcount\referencenum \referencenum=0
\newcount\@@prerefcount \@@prerefcount=0
\newcount\@@thisref
\newcount\@@lastref
\newcount\@@loopref
\newcount\@@refseq
\newdimen\refnumindent
\let\@undefreflist=\empty

\def\referencenumstyle#1{\@setnumstyle\referencenum{#1}}

\def\referencestyle#1{\usename{@ref#1}}

\def\@refsequential{%
     \gdef\@refpredef##1{\global\advance\referencenum by\@ne
          \let\^^\=0\label{##1}{\^^\{\the\referencenum}}%
          \useafter\gdef{ref@\the\referencenum}{{##1}{\undefinedlabelformat}}}%
     \gdef\@reference##1##2{%
          \ifundefinedlabel##1\then
          \else\def\^^\####1{\global\@@thisref=####1\relax}\putlab{##1}%
               \useafter\gdef{ref@\the\@@thisref}{{##1}{##2}}%
          \fi}%
     \gdef\endputreferences{%
          \loop\ifnum\@@loopref<\referencenum
                    \advance\@@loopref by\@ne
                    \expandafter\expandafter\expandafter\@printreference
                         \csname ref@\the\@@loopref\endcsname
          \repeat
          \par}}

\def\@refpreordered{%
     \gdef\@refpredef##1{\global\advance\referencenum by\@ne
          \additemR##1\to\@undefreflist}%
     \gdef\@reference##1##2{%
          \ifundefinedlabel##1\then
          \else\global\advance\@@loopref by\@ne
               {\let\^^\=0\label{##1}{\^^\{\the\@@loopref}}}%
               \@printreference{##1}{##2}%
          \fi}
     \gdef\endputreferences{%
          \def\^^\####1{\useiflab{####1}\then
               \else\reference{####1}{\undefinedlabelformat}\fi}%
          \@undefreflist
          \par}}

\def\beginprereferences{\par
     \def\reference##1##2{\global\advance\referencenum by1\@ne
          \let\^^\=0\label{##1}{\^^\{\the\referencenum}}%
          \useafter\gdef{ref@\the\referencenum}{{##1}{##2}}}}
\def\endprereferences{\global\@@prerefcount=\the\referencenum\par}

\def\beginputreferences{\par
     \refnumindent=\z@\@@loopref=\z@
     \loop\ifnum\@@loopref<\referencenum
               \advance\@@loopref by\@ne
               \setbox\z@=\hbox{\referencenum=\@@loopref
                    \referencenumformat\enskip}%
               \ifdim\wd\z@>\refnumindent\refnumindent=\wd\z@\fi
     \repeat
     \putreferenceformat
     \@@loopref=\z@
     \loop\ifnum\@@loopref<\@@prerefcount
               \advance\@@loopref by\@ne
               \expandafter\expandafter\expandafter\@printreference
                    \csname ref@\the\@@loopref\endcsname
     \repeat
     \let\reference=\@reference}

\def\@printreference#1#2{\ifx#2\undefinedlabelformat\err@undefinedref{#1}\fi
     \noindent\ifdraft\rlap{\hskip\hsize\hskip.1in \tt#1}\fi
     \llap{\referencenum=\@@loopref\referencenumformat\enskip}#2\par}

\def\reference#1#2{{\par\refnumindent=\z@\putreferenceformat\noindent#2\par}}

\def\putref#1{\trim#1\to\@@refarg
     \expandafter\ifnoarg\@@refarg\then
          \toks@={\relax}%
     \else\@@lastref=-\@m\def\@@refsep{}\def\@more{\@nextref}%
          \toks@={\@nextref#1,,}%
     \fi\the\toks@}
\def\@nextref#1,{\trim#1\to\@@refarg
     \expandafter\ifnoarg\@@refarg\then
          \let\@more=\relax
     \else\ifundefinedlabel\@@refarg\then
               \expandafter\@refpredef\expandafter{\@@refarg}%
          \fi
          \def\^^\##1{\global\@@thisref=##1\relax}%
          \global\@@thisref=\m@ne
          \setbox\z@=\hbox{\putlab\@@refarg}%
     \fi
     \advance\@@lastref by\@ne
     \ifnum\@@lastref=\@@thisref\advance\@@refseq by\@ne\else\@@refseq=\@ne\fi
     \ifnum\@@lastref<\z@
     \else\ifnum\@@refseq<\thr@@
               \@@refsep\def\@@refsep{,}%
               \ifnum\@@lastref>\z@
                    \advance\@@lastref by\m@ne
                    {\referencenum=\@@lastref\putrefformat}%
               \else\undefinedlabelformat
               \fi
          \else\def\@@refsep{--}%
          \fi
     \fi
     \@@lastref=\@@thisref
     \@more}

%************************************************************
%*
%*             Job streaming
%*
%************************************************************
\message{streaming,}

\newif\ifstreaming

\def\streamto{\defaultoption[\jobname]\@streamto}
\def\@streamto[#1]{\global\streamingtrue
     \immediate\write\sixt@@n{--> Streaming to #1.str}%
     \newwrite\streamout\immediate\openout\streamout=#1.str }

\def\streamfrom{\defaultoption[\jobname]\@streamfrom}
\def\@streamfrom[#1]{\newread\streamin\openin\streamin=#1.str
     \ifeof\streamin
          \expandafter\err@nostream\expandafter{#1.str}%
     \else\immediate\write\sixt@@n{--> Streaming from #1.str}%
          \let\@@labeldef=\gdef
          \ifstreaming
               \edef\@elc{\endlinechar=\the\endlinechar}%
               \endlinechar=\m@ne
               \loop\read\streamin to\@@scratcha
                    \ifeof\streamin
                         \streamingfalse
                    \else\toks@=\expandafter{\@@scratcha}%
                         \immediate\write\streamout{\the\toks@}%
                    \fi
                    \ifstreaming
               \repeat
               \@elc
               \input #1.str
               \streamingtrue
          \else\input #1.str
          \fi
          \let\@@labeldef=\xdef
     \fi}

\def\streamcheck{\ifstreaming
     \immediate\write\streamout{\pagenum=\the\pagenum}%
     \immediate\write\streamout{\footnotenum=\the\footnotenum}%
     \immediate\write\streamout{\referencenum=\the\referencenum}%
     \immediate\write\streamout{\chapternum=\the\chapternum}%
     \immediate\write\streamout{\sectionnum=\the\sectionnum}%
     \immediate\write\streamout{\subsectionnum=\the\subsectionnum}%
     \immediate\write\streamout{\equationnum=\the\equationnum}%
     \immediate\write\streamout{\subequationnum=\the\subequationnum}%
     \immediate\write\streamout{\figurenum=\the\figurenum}%
     \immediate\write\streamout{\subfigurenum=\the\subfigurenum}%
     \immediate\write\streamout{\tablenum=\the\tablenum}%
     \immediate\write\streamout{\subtablenum=\the\subtablenum}%
     \immediate\closeout\streamout
     \fi}

%************************************************************
%*
%*             Error messages
%*
%************************************************************

\def\err@badtypesize{%
     \errhelp={The limited availability of certain fonts requires^^J%
          that the base type size be 10pt, 12pt, or 14pt.^^J}%
     \errmessage{--> Illegal base type size}}

\def\err@badsizechange{\immediate\write\sixt@@n
     {--> Size change not allowed in math mode, ignored}}

\def\err@sizetoolarge#1{\immediate\write\sixt@@n
     {--> \noexpand#1 too big, substituting HUGE}}

\def\err@sizenotavailable#1{\immediate\write\sixt@@n
     {--> Size not available, \noexpand#1 ignored}}

\def\err@fontnotavailable#1{\immediate\write\sixt@@n
     {--> Font not available, \noexpand#1 ignored}}

\def\err@sltoit{\immediate\write\sixt@@n
     {--> Style \noexpand\sl not available, substituting \noexpand\it}%
     \it}

\def\err@bfstobf{\immediate\write\sixt@@n
     {--> Style \noexpand\bfs not available, substituting \noexpand\bf}%
     \bf}

\def\err@badgroup#1#2{%
     \errhelp={The block you have just tried to close was not the one^^J%
          most recently opened.^^J}%
     \errmessage{--> \noexpand\end{#1} doesn't match \noexpand\begin{#2}}}

\def\err@badcountervalue#1{\immediate\write\sixt@@n
     {--> Counter (#1) out of bounds}}

\def\err@extrafootnotemark{\immediate\write\sixt@@n
     {--> \noexpand\footnotemark command
          has no corresponding \noexpand\footnotetext}}

\def\err@extrafootnotetext{%
     \errhelp{You have given a \noexpand\footnotetext command without first
          specifying^^Ja \noexpand\footnotemark.^^J}%
     \errmessage{--> \noexpand\footnotetext command has no corresponding
          \noexpand\footnotemark}}

\def\err@labelredef#1{\immediate\write\sixt@@n
     {--> Label "#1" redefined}}

\def\err@badlabelmatch#1{\immediate\write\sixt@@n
     {--> Definition of label "#1" doesn't match value in \jobname.lab}}

\def\err@needlabel#1{\immediate\write\sixt@@n
     {--> Label "#1" cited before its definition}}

\def\err@undefinedlabel#1{\immediate\write\sixt@@n
     {--> Label "#1" cited but never defined}}

\def\err@undefinedeqn#1{\immediate\write\sixt@@n
     {--> Equation "#1" not defined}}

\def\err@undefinedref#1{\immediate\write\sixt@@n
     {--> Reference "#1" not defined}}

\def\err@nostream#1{%
     \errhelp={You have tried to input a stream file that doesn't exist.^^J}%
     \errmessage{--> Stream file #1 not found}}

%************************************************************
%*
%*             Initialization
%*
%************************************************************
\message{jyTeX initialization}

\everyjob{\immediate\write16{--> jyTeX version \fmtversion}%
     \edef\@@jobname{\jobname}%
%     \openin0=\inputpath jysupp
%     \ifeof0
%     \else\closein0
%          \immediate\write16{--> Additional macros loaded from jysupp.tex}%
%          \jyinput jysupp
%     \fi
%     \openin0=\inputpath jylocal
%     \ifeof0
%     \else\closein0
%          \immediate\write16{--> Additional macros loaded from jylocal.tex}%
%          \jyinput jylocal
%     \fi
     \edef\jobname{\@@jobname}%
     \settime
     \openin0=\jobname.lab
     \ifeof0
     \else\closein0
          \immediate\write16{--> Getting labels from file \jobname.lab}%
          \input\jobname.lab
     \fi}

%************** Spacing *************************************

\def\fixedskipslist{%
     \^^\{\topskip}%
     \^^\{\splittopskip}%
     \^^\{\maxdepth}%
     \^^\{\skip\topins}%
     \^^\{\skip\footins}%
     \^^\{\headskip}%
     \^^\{\footskip}}

\def\scalingskipslist{%
     \^^\{\p@renwd}%
     \^^\{\delimitershortfall}%
     \^^\{\nulldelimiterspace}%
     \^^\{\scriptspace}%
     \^^\{\jot}%
     \^^\{\normalbaselineskip}%
     \^^\{\normallineskip}%
     \^^\{\normallineskiplimit}%
     \^^\{\baselineskip}%
     \^^\{\lineskip}%
     \^^\{\lineskiplimit}%
     \^^\{\bigskipamount}%
     \^^\{\medskipamount}%
     \^^\{\smallskipamount}%
     \^^\{\parskip}%
     \^^\{\parindent}%
     \^^\{\abovedisplayskip}%
     \^^\{\belowdisplayskip}%
     \^^\{\abovedisplayshortskip}%
     \^^\{\belowdisplayshortskip}%
     \^^\{\abovechapterskip}%
     \^^\{\belowchapterskip}%
     \^^\{\abovesectionskip}%
     \^^\{\belowsectionskip}%
     \^^\{\abovesubsectionskip}%
     \^^\{\belowsubsectionskip}}

%************** Document layout *****************************

\def\twoupsetup{%                                % setup for twoup style
     \topmargin=.75in
     \leftmargin=.5in
     \vsize=6.9in
     \hsize=4.75in
     \fullhsize=10in
     \let\draft=\relax}

\outputstyle{normal}                             % page style

\def\marginnoteformat{\subscriptsize             % paragraphing of margin notes
     \hsize=1in \baselinestretch=1000 \everypar={}%
     \tolerance=5000 \hbadness=5000 \parskip=0pt \parindent=0pt
     \leftskip=0pt \rightskip=0pt \raggedright}

\head={\ifdraft\normalfonts\it\hfil DRAFT\hfil   % format of headline
     \llap{\number\day\ \monthword\month\ \militarytime}\else\hfil\fi}
\foot={\hfil\normalfonts\numstyle\pagenum\hfil}  % format of footline

\normalbaselineskip=12pt                         % usual \baselineskip
\normallineskip=0pt                              % usual \lineskip
\normallineskiplimit=0pt                         % usual \lineskiplimit
\normalbaselines                                 % set \baselineskip

\topskip=.85\baselineskip \splittopskip=\topskip \headskip=2\baselineskip
\footskip=\headskip

\pagenumstyle{arabic}                            % counter style

\parskip=0pt                                     % no skip between paragraphs
\parindent=20pt                                  % usual \parindent

\baselinestretch=1000                            % set \big-, \med-, \smallskip

%************** Sectioning **********************************

\chapterstyle{left}                              % position of heading
\chapternumstyle{blank}                          % counter style
\def\chapterbreak{\newpage}                      % break before heading
\abovechapterskip=0pt                            % space before heading
\belowchapterskip=1.5\baselineskip               % space after heading
     plus.38\baselineskip minus.38\baselineskip
\def\chapternumformat{\numstyle\chapternum.}     % format of heading counter

\sectionstyle{left}                              % position of heading
\sectionnumstyle{blank}                          % counter style
\def\sectionbreak{\vskip0pt plus4\baselineskip\penalty-100
     \vskip0pt plus-4\baselineskip}              % break before heading
\abovesectionskip=1.5\baselineskip               % space before heading
     plus.38\baselineskip minus.38\baselineskip
\belowsectionskip=\the\baselineskip              % space after heading
     plus.25\baselineskip minus.25\baselineskip
\def\sectionnumformat{%                          % format of heading counter
     \ifblank\chapternumstyle\then\else\numstyle\chapternum.\fi
     \numstyle\sectionnum.}

\subsectionstyle{left}                           % position of heading
\subsectionnumstyle{blank}                       % counter style
\def\subsectionbreak{\vskip0pt plus4\baselineskip\penalty-100
     \vskip0pt plus-4\baselineskip}              % break before heading
\abovesubsectionskip=\the\baselineskip           % space before heading
     plus.25\baselineskip minus.25\baselineskip
\belowsubsectionskip=.75\baselineskip            % space after heading
     plus.19\baselineskip minus.19\baselineskip
\def\subsectionnumformat{%                       % format of heading counter
     \ifblank\chapternumstyle\then\else\numstyle\chapternum.\fi
     \ifblank\sectionnumstyle\then\else\numstyle\sectionnum.\fi
     \numstyle\subsectionnum.}

%************** Footnotes ***********************************

\footnotenumstyle{symbols}                       % counter style
\footnoteskip=0pt                                % jyTeX spacing parameter
\def\footnotenumformat{\numstyle\footnotenum}    % \footnotemark format
\def\footnoteformat{\footnotesize                % paragraphing of text
     \everypar={}\parskip=0pt \parfillskip=0pt plus1fil
     \leftskip=1em \rightskip=0pt
     \spaceskip=0pt \xspaceskip=0pt
     \def\\{\ifhmode\ifnum\lastpenalty=-10000
          \else\hfil\penalty-10000 \fi\fi\ignorespaces}}

%************** Labels **************************************

\def\undefinedlabelformat{$\bullet$}             % mark for undefined label

%************** Equation numbering **************************

\equationnumstyle{arabic}                        % counter style
\subequationnumstyle{blank}                      % counter style
\figurenumstyle{arabic}                          % counter style
\subfigurenumstyle{blank}                        % counter style
\tablenumstyle{arabic}                           % counter style
\subtablenumstyle{blank}                         % counter style

\eqnseriesstyle{alphabetic}                      % sub-counter style for series
\figseriesstyle{alphabetic}                      % sub-counter style for series
\tblseriesstyle{alphabetic}                      % sub-counter style for series

\def\puteqnformat{\hbox{%                        % equation number format
     \ifblank\chapternumstyle\then\else\numstyle\chapternum.\fi
     \ifblank\sectionnumstyle\then\else\numstyle\sectionnum.\fi
     \ifblank\subsectionnumstyle\then\else\numstyle\subsectionnum.\fi
     \numstyle\equationnum
     \numstyle\subequationnum}}
\def\putfigformat{\hbox{%                        % figure number format
     \ifblank\chapternumstyle\then\else\numstyle\chapternum.\fi
     \ifblank\sectionnumstyle\then\else\numstyle\sectionnum.\fi
     \ifblank\subsectionnumstyle\then\else\numstyle\subsectionnum.\fi
     \numstyle\figurenum
     \numstyle\subfigurenum}}
\def\puttblformat{\hbox{%                        % table number format
     \ifblank\chapternumstyle\then\else\numstyle\chapternum.\fi
     \ifblank\sectionnumstyle\then\else\numstyle\sectionnum.\fi
     \ifblank\subsectionnumstyle\then\else\numstyle\subsectionnum.\fi
     \numstyle\tablenum
     \numstyle\subtablenum}}

%************** Reference numbering *************************

\referencestyle{sequential}                      % referencing method
\referencenumstyle{arabic}                       % counter style
\def\putrefformat{\numstyle\referencenum}        % format of reference citation
\def\referencenumformat{\numstyle\referencenum.} % format of number in list
\def\putreferenceformat{%                        % paragraphing of list
     \everypar={\hangindent=1em \hangafter=1 }%
     \def\\{\hfil\break\null\hskip-1em \ignorespaces}%
     \leftskip=\refnumindent\parindent=0pt \interlinepenalty=1000 }

%************** Font initialization *************************

\normalsize

%*****************************************************************************

\def\fmtversion{2.6M (June 1992)}

\catcode`\@=12
% ------------------ End of jytex.tex -----------------

%\input jytex.tex   % available from hep-th
\typesize=10pt \magnification=1200 \baselineskip17truept
%\baselineskip25truept
\footnotenumstyle{arabic} \hsize=6truein\vsize=8.5truein
%\draft
%\leftmargin=1.25in
%\oddleftmargin=.5in
%\evenleftmargin=1.5in
\sectionnumstyle{blank}
\chapternumstyle{blank}
\chapternum=1
\sectionnum=1
\pagenum=0
%\referencestyle{preordered}
% title style follows

\def\begintitle{\pagenumstyle{blank}\parindent=0pt
\begin{narrow}[0.4in]}
\def\endtitle{\end{narrow}\newpage\pagenumstyle{arabic}}

% exercise style follows

\def\beginexercise{\vskip 20truept\parindent=0pt\begin{narrow}[10
truept]}
\def\endexercise{\vskip 10truept\end{narrow}}

% **************    my jyTeX abbreviations   *****************

\def\eql#1{\eqno\eqnlabel{#1}}
\def\ref{\reference}
\def\peq{\puteqn}
\def\pref{\putref}

\def\mgn{\marginnote}
\def\bex{\begin{exercise}}
\def\eex{\end{exercise}}

% *********************** My definitions ************************

\font\open=msbm10 %scaled\magstep1 % For VAX. Borde p195.

 %scaled\magstep1 % For VAX. Borde p195.
%\font\open=msym10 %scaled\magstep1 % For Arbortxt on PC
%\font\opens=msym8 %scaled\magstep1 % For Arbortxt on PC
  % For Arbortxt on PC, and VAX. Borde p199
 
\def\StretchRtArr#1{{\count255=0\loop\relbar\joinrel\advance\count255 by1
\ifnum\count255<#1\repeat\rightarrow}}
\def\StretchLtArr#1{\,{\leftarrow\!\!\count255=0\loop\relbar
\joinrel\advance\count255 by1\ifnum\count255<#1\repeat}}

\def\StretchLRtArr#1{\,{\leftarrow\!\!\count255=0\loop\relbar\joinrel\advance
\count255 by1\ifnum\count255<#1\repeat\rightarrow\,\,}}

\def\mbox#1{{\leavevmode\hbox{#1}}}

\def\hspace#1{{\phantom{\mbox#1}}}
\def\oZ{\mbox{\open\char90}}

\def\al{\alpha}

\def\b0{{\bf0}}
 %in jyTeX
 %in jyTeX
 %in jyTeX
 %in jyTeX
 %in jyTeX
 %in jyTeX
 %in jyTeX
 %in jyTeX
 %in jyTeX
 %in jyTeX
 %in jyTeX
 %in jyTeX
 %in jyTeX
 %in jyTeX
% in jyTeX
% in jyTeX
% in jyTeX
% in jyTeX
% in jyTeX
% in jyTeX
\def\be{\beta}

\def\de{\delta}

\def\la{\lambda}

\def\th{\theta}
\def\Th{\Theta}

\def\De{\Delta}

\def\caG{{\cal G}}
\def\caF{{\cal F}}

\def\caH{{\cal H}}

\def\sc{{\rm sc }}

     % Newline

\def\frac#1/#2{\leavevmode\kern.1em
\raise.5ex\hbox{\the\scriptfont0 #1}\kern-.1em/\kern-.15em
\lower.25ex\hbox{\the\scriptfont0 #2}}
\def\sfrac#1/#2{\leavevmode\kern.1em
\raise.5ex\hbox{\the\scriptscriptfont0 #1}\kern-.1em/\kern-.15em
\lower.25ex\hbox{\the\scriptscriptfont0 #2}}

\def\gtorder{\mathrel{\raise.3ex\hbox{$>$}\mkern-14mu
             \lower0.6ex\hbox{$\sim$}}}
\def\ltorder{\mathrel{\raise.3ex\hbox{$<$}\mkern-14mu
             \lower0.6ex\hbox{$\sim$}}}

\def\semidirprod{\rlap{\ss C}\raise1pt\hbox{$\mkern.75mu\times$}}
\def\for{\lower6pt\hbox{$\Big|$}}
\def\fish{\kern-.25em{\phantom{abcde}\over \phantom{abcde}}\kern-.25em}

 %triple
%dot
 %double
%dot
 %double dot
%for small #1

\def\boxit#1{\vbox{\hrule\hbox{\vrule\kern3pt
        \vbox{\kern3pt#1\kern3pt}\kern3pt\vrule}\hrule}}
\def\dalemb#1#2{{\vbox{\hrule height .#2pt
        \hbox{\vrule width.#2pt height#1pt \kern#1pt \vrule
                width.#2pt} \hrule height.#2pt}}}

\def\ol{\overline}
        %double stroke
\def\frac#1#2{{{#1}\over{#2}}}
 %lower covariant deriv.
 %upper covariant deriv.
 %lower covariant deriv semicolon.
    %lower ordinary  deriv.
    %lower ordinary  deriv comma.

\def\noin{\noindent}

      %Connection
    %Connection'

\def\etc{{\it etc. }}

\def\eg{{\it e.g.}}
\def\ie{{\it i.e. }}
\def\cf{{\it cf }}

\def\bra#1{\langle#1\mid}
\def\ket#1{\mid#1\rangle}

 %gives average <#1>
 %gives thermal average <<#1>>
\def\br#1#2{\langle{#1}\mid{#2}\rangle}   %gives bracket <#1|#2>
   %gives comma bracket <#1,#2>
 %gives round bracket (#1,#2)
\def\rbr#1#2{({#1}\mid{#2})} %gives round bracket (#1,|#2)
 %gives big bracket <#1|#2>
  %gives
%matrix element <#1|#2|#3>

%gives reduced matrix element
%<#1||#2||#3>

\def\wt{\widetilde}

\def\3j#1#2#3#4#5#6{\left\lgroup\matrix{#1&#2&#3\cr#4&#5&#6\cr}
\right\rgroup}

\def\m?{\mgn{?}}
% KK's defs

\def\beq{\begin{eqnarray}}
\def\eeq{\end{eqnarray}}

%  *******************  Journal refs **********************

\def\aop#1#2#3{{\it Ann. Phys.} {\bf {#1}} ({#2}) #3}
\def\cjp#1#2#3{{\it Can. J. Phys.} {\bf {#1}} ({#2}) #3}
\def\cmp#1#2#3{{\it Comm. Math. Phys.} {\bf {#1}} ({#2}) #3}
\def\cqg#1#2#3{{\it Class. Quant. Grav.} {\bf {#1}} ({#2}) #3}

\def\ijmp#1#2#3{{\it Int. J. Mod. Phys.} {\bf {#1}} ({#2}) #3}

\def\jmp#1#2#3{{\it J. Math. Phys.} {\bf {#1}} ({#2}) #3}
\def\jpa#1#2#3{{\it J. Phys.} {\bf A{#1}} ({#2}) #3}
\def\jpamt#1#2#3{{\it J. Phys.A:Math.Theor.} {\bf{#1}} ({#2}) #3}
\def\jpc#1#2#3{{\it J. Phys.} {\bf C{#1}} ({#2}) #3}
\def\lnm#1#2#3{{\it Lect. Notes Math.} {\bf {#1}} ({#2}) #3}

\def\np#1#2#3{{\it Nucl. Phys.} {\bf B{#1}} ({#2}) #3}
\def\npa#1#2#3{{\it Nucl. Phys.} {\bf A{#1}} ({#2}) #3}
\def\pl#1#2#3{{\it Phys. Lett.} {\bf {#1}} ({#2}) #3}

\def\prp#1#2#3{{\it Phys. Rep.} {\bf {#1}} ({#2}) #3}
\def\pr#1#2#3{{\it Phys. Rev.} {\bf {#1}} ({#2}) #3}
\def\prA#1#2#3{{\it Phys. Rev.} {\bf A{#1}} ({#2}) #3}

\def\prD#1#2#3{{\it Phys. Rev.} {\bf D{#1}} ({#2}) #3}
\def\prE#1#2#3{{\it Phys. Rev.} {\bf E{#1}} ({#2}) #3}
\def\prl#1#2#3{{\it Phys. Rev. Lett.} {\bf #1} ({#2}) #3}

\def\rmp#1#2#3{{\it Rev. Mod. Phys.} {\bf {#1}} ({#2}) #3}

\def\zfp#1#2#3{{\it Z. f. Phys.} {\bf {#1}} ({#2}) #3}

\def\cras#1#2#3{{\it Comptes Rend. Acad. Sci. (Paris)} {\bf{#1}} (#2) #3}
\def\prs#1#2#3{{\it Proc. Roy. Soc.} {\bf A{#1}} ({#2}) #3}
\def\pcps#1#2#3{{\it Proc. Camb. Phil. Soc.} {\bf{#1}} ({#2}) #3}
\def\mpcps#1#2#3{{\it Math. Proc. Camb. Phil. Soc.} {\bf{#1}} ({#2}) #3}

\def\amsh#1#2#3{{\it Abh. Math. Sem. Ham.} {\bf {#1}} ({#2}) #3}
\def\am#1#2#3{{\it Acta Mathematica} {\bf {#1}} ({#2}) #3}
\def\aim#1#2#3{{\it Adv. in Math.} {\bf {#1}} ({#2}) #3}
\def\ajm#1#2#3{{\it Am. J. Math.} {\bf {#1}} ({#2}) #3}
\def\amm#1#2#3{{\it Am. Math. Mon.} {\bf {#1}} ({#2}) #3}

\def\aom#1#2#3{{\it Ann. of Math.} {\bf {#1}} ({#2}) #3}
\def\cjm#1#2#3{{\it Can. J. Math.} {\bf {#1}} ({#2}) #3}
\def\bams#1#2#3{{\it Bull.Am.Math.Soc.} {\bf {#1}} ({#2}) #3}

\def\cmh#1#2#3{{\it Comm. Math. Helv.} {\bf {#1}} ({#2}) #3}

\def\dmj#1#2#3{{\it Duke Math. J.} {\bf {#1}} ({#2}) #3}
\def\invm#1#2#3{{\it Invent. Math.} {\bf {#1}} ({#2}) #3}

\def\jdg#1#2#3{{\it J. Diff. Geom.} {\bf {#1}} ({#2}) #3}

\def\joa#1#2#3{{\it J. of Algebra} {\bf {#1}} ({#2}) #3}
\def\jram#1#2#3{{\it J. f. Reine u. Angew. Math.} {\bf {#1}} ({#2}) #3}
\def\jims#1#2#3{{\it J. Indian. Math. Soc.} {\bf {#1}} ({#2}) #3}
\def\jlms#1#2#3{{\it J. Lond. Math. Soc.} {\bf {#1}} ({#2}) #3}
\def\jmpa#1#2#3{{\it J. Math. Pures. Appl.} {\bf {#1}} ({#2}) #3}
\def\ma#1#2#3{{\it Math. Ann.} {\bf {#1}} ({#2}) #3}

\def\mz#1#2#3{{\it Math. Zeit.} {\bf {#1}} ({#2}) #3}
\def\ojm#1#2#3{{\it Osaka J.Math.} {\bf {#1}} ({#2}) #3}

\def\pems#1#2#3{{\it Proc. Edin. Math. Soc.} {\bf {#1}} ({#2}) #3}

\def\plb#1#2#3{{\it Phys. Letts.} {\bf {B#1}} ({#2}) #3}
\def\pla#1#2#3{{\it Phys. Letts.} {\bf {A#1}} ({#2}) #3}
\def\plms#1#2#3{{\it Proc. Lond. Math. Soc.} {\bf {#1}} ({#2}) #3}
\def\pgma#1#2#3{{\it Proc. Glasgow Math. Ass.} {\bf {#1}} ({#2}) #3}
\def\qjm#1#2#3{{\it Quart. J. Math.} {\bf {#1}} ({#2}) #3}
\def\qjpam#1#2#3{{\it Quart. J. Pure and Appl. Math.} {\bf {#1}} ({#2}) #3}

\def\rmjm#1#2#3{{\it Rocky Mountain J. Math.} {\bf {#1}} ({#2}) #3}

\def\tams#1#2#3{{\it Trans.Am.Math.Soc.} {\bf {#1}} ({#2}) #3}

% *******************   Main text *********************
\begin{title}
\vglue 0.5truein
%\righttext {MUTP/96/23}
%\righttext{hep-th/96}
\vskip15truept
%\leftline{\today}
%\vskip 30truept
\centertext {\Bigfonts \bf Poweroids revisited --} \vskip2truept
\vskip10truept\centertext{\Bigfonts \bf an old  symbolic approach}
 \vskip 20truept
\centertext{J.S.Dowker\footnote{dowker@man.ac.uk; dowkeruk@yahoo.co.uk}} \vskip
7truept \centertext{\it Theory Group,} \centertext{\it School of Physics and Astronomy,}
\centertext{\it The University of Manchester,} \centertext{\it Manchester, England} \vskip
7truept \centertext{}

\vskip 7truept

\vskip40truept
\begin{narrow}
Jeffery's 1861 computations using finite difference calculus are resurrected and extended
from forward differences to general delta operators and used to neatly prove theorems in
the Rota--Mullins theory of polynomials of binomial type (Steffensen's poweroids) allowing,
for example, compact treatments of umbral composition, the binomial property and the
connection constants. It is shown that it forms a legitimate alternative to the usual umbral
device and also anticipates a number of results obtained more recently.
\end{narrow}
\vskip 5truept
\vskip 60truept
%\righttext{Typeset in \jyTeX}
\vfil
\end{title}
\pagenum=0
\newpage

\section{\bf 1. Introduction}

As is well recorded, finite difference calculus can be developed symbolically, and was done
so, with varying intensity, from its beginnings in the 17th century. The concepts of operator
and operand became systematically more separated in the early 19th, particularly in
England with the work of Herschel, Murphy and Boole. The difference operator most often
used was the forwards one, although the advantage of moving the various expansions to a
central point, and thereby effectively introducing central differences, goes back to the time
of Newton.

Together with the concept of difference operator goes that of the {\it factorial} which plays
the role that a power does in differential calculus. In this paper I wish to draw attention to a
forgotten, purely symbolic, approach involving these particular notions. This is not in
response to any practical, or even theoretical, necessity, but simply to present what is, I
think, a different slant on these ideas which might amuse and could even have some
technical virtue. Even if not, I can at least advertise some ignored, and rather elegant,
symbolic work \footnote{ I distinguish between symbolic methods and symbolic notation.}
that appeared at the same time as umbral methods in the mid 19th century. I also make
some other historical points.

I will use a form for the difference operator  that encompasses the usual ones and is due, in
this connection, to Steffensen who first, [\pref{Steffensen2}], introduced a cut down version
in order to defend his definitions of the {\it factorial} function in [\pref{Steffensen}]. He
then extended it in his notion of `poweroids', or generalised powers (or factorials),
[\pref{Steffensen3}]. This theory becomes partly subsumed  into the theory of finite
operators exploited by Rota and coworkers in their polynomial underpinning of the umbral
technique which has been extended and applied in many, sometimes advanced ways. By
contrast, the methods used here are simple minded and explicit. For the most part my
notation will be that employed in the traditional works on finite differences despite the
attractions of the more upmarket Dirac bra(c)ket, used \eg\ by \ Roman,
[\pref{Roman,RandR}].

My aim is to see to what extent the early symbolic techniques can be accommodated to
general difference operators and how they compare with the finite operator and modern
umbral approaches. These have been used to analyse aspects of discretised quantum theory
and one can find here useful summaries. See, for example, the review by Levi {\it et al},
[\pref{L}].

\section{\bf2. The general difference operator. Poweroids}

\noin In this review section I start out with the special and then generalise.

Steffensen, [\pref{Steffensen2}], defines the operator
  $$
   \th={E^{\al}\big(E^\be-1\big)\over\be}
   \eql{th}
  $$
$E$ being, as usual, the unit translation operator, $E=e^D$. $D$ is the derivative,
represented by $d/dx$. $\th$ is the most general divided difference of the first order.

Writing, temporarily, $\th(\al,\be)$ for $\th$, the four standard differences are given by the
special values,
  $$\eqalign{
  &\De=\th(0,1)\,,\qquad\quad\ldots\quad{\rm forwards}\cr
  &\nabla=\th(0,-1)\,,\qquad\,\ldots\quad{\rm backwards}\cr
  &\de\,=\,\th(-1/2,1)\,,\quad\,\ldots\quad{\rm central}\cr
  &D=\th(0,0)\,,\qquad\quad\ldots\quad{\rm confluent}\,.
  }
  \eql{four}
  $$
I have adopted Aitken's terminology of {\it confluent} for the limiting case of the ordinary
differential.

The basic fact about $\th$ is its action on a particular generalised `factorial'  which is defined
by, (I here deviate from Steffensen's notation),
  $$
   x^{\{n\}}\equiv x\big(x-n\al-\be\big)\big(x-n\al-2\be\big)
   \ldots\big(x-n\al-(n-1)\be\big)\,,
   \eql{fact}
  $$
with $x^{\{0\}}=1$ and $x^{\{1\}}=x$.

Direct calculation shows the expected, and fundamental, behaviour,
  $$
    \th x^{\{n\}}=nx^{\{n-1\}}\,.
    \eql{act}
  $$
The four cases listed in (\peq{four}), give $x^{\{n\}}=$ $x^{(n)}$, $x^{(-n)}$,
$x^{[n]}$ and $x^n$ in turn in Steffensen's notation (regarding which consult Aitken,
[\pref{Aitken}]).

In [\pref{Steffensen3}] Steffensen generalised this whole structure and replaced the specific
operator, (\peq{th}), by a more general one,
  $$
  \th=\phi(D)\,,
  \eql{psi}
  $$
where, by definition, the function $\phi(D)\to D$  as $D\to 0$ and was taken as a formal
power series, $\sum_{\nu=1}^\infty k_\nu\,D^\nu$.

It is therefore clear that $\th^m 0^n$ is zero if $m>n$ and so $x^{\{n\}}$, now {\it
defined} (uniquely if  one adds $x^{\{0\}}=1$ and $0^{\{n\}}=\de_n^0$) by (\peq{act})
and (\peq{psi}), are polynomials of degree $n$ which always have a factor of $x$. They are
termed `poweroids' by Steffensen, `basic polynomials' by Sheffer and by Rota and
`associated polynomials' by Roman and Rota. They are also often referred to as `Sheffer
sequences' (of a certain type). Further terminology has $\th$ as a `delta operator', after
Hildebrand. I will use any, or all, of these terms.

One drawback of the notation $x^{\{n\}}$ is that it does not indicate the associated
operator, which is sometimes useful. Hence the alternative, $b^\th_n(x)$.
\section{\bf3. General difference theorems. Duality.}

As I have mentioned, there exists a considerable body of early work at the purely symbolic
level concerned with the standard differences, (\peq{four}) (and mostly then with the
forwards quantities). In the present paper I give a treatment, using the operator $\th$, of
selected aspects of the corresponding calculus, and base my approach on two papers by
H.M.Jeffery \footnote{ Henry Martyn Jeffery (1826-1891), of west country origins,
graduated sixth wrangler in 1849 and, after a spell as lecturer at the Civil Engineering
College in Putney, was appointed in 1852 as second master at Pate's Grammar School in
Cheltenham, an old foundation linked to his Cambridge college, St. Catherine's. He
remained there until retirement in 1882, becoming headmaster in 1868 and producing a
large number of papers, mostly on algebraic geometry, all of which have disappeared from
view, apart from the odd reference. The analytical works, however, are of more than a little
interest. His work was well regarded at the time and he was elected FRS in 1880. Along with
Blissard, Horner, Kirkman and others, he might be placed in that peculiarly English set of
`gentleman mathematicians', many clergymen, whose daily bread was earned elsewise.
The pages of The Quarterly Journal of Mathematics and The Messenger of Mathematics are
replete with their contributions, the Cambridge contingent particularly.},
[\pref{Jeffery,Jeffery2}] in 1861 and 1862.

As in an earlier work, [\pref{Dowcen}], It is convenient to begin with transcribing some
basic theorems assembled by Jeffery, [\pref{Jeffery}],\footnote{ The main objectives of
Jeffery's first paper are combinatorial, but he uses and developes results of general
operational applicability.} which are concerned with the transformation of symbolic
expressions. I firstly just state these generalised theorems,

  $$\eqalign{
  &{\rm A}:\quad f(\th)\,0^n=n{f(\th)\over D}\,0^{n-1}=\ldots
  =n!{f(\th)\over D^n}0^0\cr
  &{\rm B}:\quad f(\th)\,0^n={df(\th)\over dD}\,0^{n-1}=\ldots
  ={d^n f(\th)\over dD^n}0^0\cr
  &{\rm C}:\quad f(D)\,0^{\{n\}}=n{f(D)\over \th}\,0^{\{n-1\}}=\ldots
  =n!{ f(D)\over \th^n}0^{\{0\}}\cr
  &{\rm D}:\quad f(D)\,0^{\{n\}}={df(D)\over d\th}0^{\{n-1\}}=\ldots
  ={d^n f(D)\over d\th^n}0^{\{0\}}\,,\cr
  }
  \eql{JABCD}
  $$
which Jeffery, [\pref{Jeffery}], gives for the particular case, $\th=\De$.

A and C express the simple fact that the symbols $\th$ and $D$ equally commute when
acting upon $x$ at $x=0$ as they do when $x$ is current.

Theorem B can be more `universally' expressed as $f(E)\,0^n=Ef'(E)\,0^{n-1}$ and is
proved in \eg\ Boole, [\pref{Boole2}] p.28. Many of these relations go back at least as far
as Herschel, [\pref{Herschel}].

To derive the more significant D, I start from the very basic action on factorials,
(\peq{act}), which implies,
  $$
  F(\th)\,0^{\{r\}}=F(D)\,0^r\,.
  \eql{factact}
  $$

I can then define a new function $f$ by the functional relation,
  $$
  F(\th)=F\big(\phi(D)\big)\equiv f(D)=f\big(\phi^{-1}(\th)\big)\,,
  $$
and just replace $F(\th)$ by $f(D)$ on the left--hand side of (\peq{factact}) while, on the
right--hand side, replace $\th$ by $D$ in $f$, {\it regarded as a function of} $\th$. So, if I
define the operator $\eta$, by
   $$
    \eta=\phi^{-1}(D)\,,
    \eql{xidef}
   $$
I obtain the elegant result,
  $$
  f(D)\,0^{\{r\}}=f(\eta)\,0^r\,,
  \eql{factact2}
  $$
which is the dual of (\peq{factact}). This equation is due to Jeffery, [\pref{Jeffery}], derived
slightly differently.\footnote{ While Jeffery uses forward differences, his analysis is valid
generally. For him, $\eta=\zeta=\log(1+D)$.}

The proof of Theorem D now follows more or less directly. Applying Theorem B,  one finds,
  $$\eqalign{
  f(D)\,0^{\{r\}}&=f(\eta)\,0^r=f'(\eta)\,{d\over dD}\phi^{-1}(D)\,0^{r-1}\cr
  &=f'(\eta)\bigg[A_0+A_1D+A_2D^2+\ldots\bigg]{0}^{r-1}\cr
 &=f'(\eta)\bigg[A_00^{r-1}+(r-1)A_10^{r-2}+\ldots\bigg]\cr
  &=f'(D)\bigg[A_00^{\{r-1\}}+(r-1)A_10^{\{r-2\}}+\ldots\bigg]\cr
  &=f'(D){d\over d\th}\phi^{-1}(\th)\,0^{\{r-1\}}
  ={df(D)\over d\th}\,0^{\{r-1\}}\,,
  }
  \eql{thmD}
  $$
using (\peq{factact2}) again and $\phi^{-1}(\th)=D$. This is Theorem D.

It is possible to pass immediately from the first line to the last just on the basis of
(\peq{factact2}) which implies the duality replacements $\eta\leftrightarrow D$,
$D\leftrightarrow\th$ and $0^r\leftrightarrow 0^{\{r\}}$. The expansions, which aren't
needed explicitly anyway, serve merely to reinforce the validity of this formal
transformation by displaying its meaning.

\section{\bf4. Expansions.  Jeffery's equation}

The fundamental Maclaurin--like expansion of a function in poweroids is often required,
\footnote{ To be as general as possible, one should go into the question of remainders.
However for formal considerations, as here, this can be put aside or attention can be
restricted to polynomials} \ie
  $$
   f(x)=\sum_{\nu=0}^\infty x^{\{\nu\}}\,{\th^\nu f(0)\over \nu!}\,,
   \eql{fexp}
  $$
\eg\ Steffensen, [\pref{Steffensen3}] Equ.(7), and can be elevated to a Theorem.

 An important example is the power,
  $$
    x^n=\sum_{\nu=0}^n x^{\{\nu\}}\,{\th^\nu0^n\over \nu!}\,,\quad n\in \oZ\,,
    \eql{pow}
  $$
which introduces the (generalised) `differences of nothing', $\th^\nu 0^n$. These could be
regarded as fundamental data \footnote{ This attitude seems to date back to Brinkley,
1807, [\pref{Brinkley}].} and have been tabulated for the standard cases, (\peq{four})
from historical times.

Conversely, the (Maclaurin) expansion of the poweroid is,
  $$\eqalign{
  x^{\{n\}}&=\sum_{\nu=1}^n x^{\nu}\, {D^{\nu}0^{\{n\}}\over\nu!}=
  \sum_{\nu=1}^n x^{\nu}\, {\eta^{\nu}0^{n}\over\nu!}\cr
  &=e^{x\,D}\,0^{\{n\}}=e^{x\,\eta}\,0^n\,,\cr
  }
  \eql{cdf}
  $$
using (\peq{factact2}) to rewrite the differentials of nothing and also formally summing the
{\it power} series ($1.0^n=0$ for $n\ne0$). Comparing Jeffery, [\pref{Jeffery}], \S5, I
refer to (\peq{cdf}) as Jeffery's equation and $e^{x\eta}$ as the Jeffery operator. Again,
the standard coefficients have been tabulated, \eg\ [\pref{Jeffery}], [\pref{Steffensen}],
[\pref{Thiele}], [\pref{Shovelton}].

The coefficients in (\peq{pow}) and (\peq{cdf}) could be called `generalised factorial
coefficients of the second and first kind', respectively. Note that these numbers vanish if
$\nu>n$. They are examples of the `connection constants' of [\pref{RandM}] and
[\pref{RKO}] discussed later, in section 14.

\section{\bf5. Rodrigues--type recursion relation}

As an illustration of the use of the transformations (\peq{JABCD}) I derive Steffensen's
poweroid recursion, [\pref{Steffensen3}] eqn.(17),
  $$
   \th'\,x^{\{n\}-1}=x^{\{n-1\}}
   \eql{precur}
  $$
which here follows quickly from (\peq{cdf}) upon applying Theorem D with
$f(D)=e^{x\,D}$. Spelling out the details,
  $$\eqalign{
  x^{\{n\}}&=e^{x\,D}\,0^{\{n\}}=x{1\over \th'}\,e^{x\,D}\,0^{\{n-1\}}=
  x{1\over \th'}\,(0+x)^{\{n-1\}}\cr
  &=x{1\over \th'}\,x^{\{n-1\}}\,,\quad {\rm QED}\,.
  }
  $$
On the first line $\th'$ acts on $0$ while on the last its action has been transferred to $x$.

This recursion is the Rodrigues--type formula of Rota and Mullin, [\pref{RandM}] Theorem
4.4.
\section{\bf 6. Inverse poweroids}

Since $\eta$ is a function of $D$ vanishing at $D=0$, it can act as a delta operator, the
corresponding `$\eta$' operator being just $\th$. Hence it is possible to interchange
$\th\leftrightarrow\eta$ in the previous analysis if, at the same time, the poweroid,
$x^{\{n\}}$, is replaced by that, $x^{\{n\}^{-1}}$ say, associated with $\eta$ and given
by,
  $$
  x^{\{n\}^{-1}}=e^{x\th}\,0^n\,.
  \eql{powerd}
  $$
Corresponding to (\peq{factact2}), one has,
  $$
  f(D)\,0^{\{r\}^{-1}}=f(\th)\,0^r\,,
  \eql{factact3}
  $$
so the expansion (\peq{cdf}) turns into,
  $$\eqalign{
  x^{\{n\}^{-1}}&=\sum_{\nu=1}^n x^{\nu}\, {D^{\nu}0^{\{n\}^{-1}}\over\nu!}=
  \sum_{\nu=1}^n x^{\nu}\, {\th^{\nu}0^{n}\over\nu!}\,,\cr
  }
  \eql{cdf2}
  $$
taking us back to (\peq{powerd}).

An alternative notation is often advantageous, \cf\ Riordan, [\pref{Riordan}], and I set,
  $$\eqalign{
 &D^{\nu}0^{\{n\}}=\eta^\nu0^n\equiv\nu!\, g(n,\nu)\cr
  &D^{\nu}0^{\{ n\}^{-1}}=\th^\nu0^n\equiv\nu!\, \ol g(\nu,n)\,,
  }
  \eql{not}
  $$
so that,
  $$\eqalign{
   x^{\{n\}}&=\sum_{\nu=1}^n  g(n,\nu)\,x^{\nu} \cr
  x^{\{n\}^{-1}}&=\sum_{\nu=1}^n x^{\nu} \ol g(\nu,n)\, \cr
  }
  \eql{cdf3}
  $$
and
  $$
   x^{\nu}=\sum_{n=1}^n \ol g(\nu,n) x^{\{n\}}=\sum_{n=1}^n x^{\{n\}^{-1}} g(n,\nu)
  $$
  
For the forwards case, see (\peq{four}), $g=s$ and $\wt{\ol g}=S$, the Stirling numbers,
and for the central system, $g=t$ and $\wt{\ol g}=T$, in Riordan's notation,
[\pref{Riordan}]. Of course, for the confluent case, $g(n,\nu)=\ol g(\nu,n)=\de_n^\nu$,
the Kronecker delta.

This notation comes into play in section 15 in connection with representative notation.

\section{\bf7. Completeness}
As a further illustration of this symbolic formalism, I derive an intrinsic relation between the
first and second kind factorial numbers, \ie between the $\eta^m\,0^n$ and $\th^r\,0^s$.

I start from the definition of the second kind, (\peq{pow}),
  $$\eqalign{
  x^{n}&=\sum_{\nu=1}^n x^{\{\nu\}}\, {\th^{\nu}0^{n}\over\nu!}
  =\sum_{\nu=1}^n e^{x\eta}\,0^\nu\, {\th^{\nu}0^{n}\over\nu!}\cr
  &=e^{x\eta}\,e^{0.\th}\,0^n\cr
  &=e^{x\th}\,e^{0.\eta}\,0^n\,,\cr
  }
  \eql{pow2}
  $$
and equate powers of $x$ to get the compact statements of `orthogonality',
  $$\eqalign{
   &\eta^m\,e^{0.\,\th}\,0^n=n!\,\de_m^n\cr
   &\th^m\,e^{0.\,\eta}\,0^n=n!\,\de_m^n\,.\cr
   }
   \eql{inv}
  $$
The $\de$ on the right--hand side is a Kronecker delta.

Equation (\peq{inv}) is really a statement about inverses, and tantamount to completeness
(of polynomial bases). \cf\ Riordan, [\pref{Riordan}] p.213 for central factorial (Stirling)
numbers.

Another way of expressing (\peq{pow2}) (or (\peq{inv})), is,
  $$
  e^{xy}=e^{x\eta}\,e^{0.\th}\,e^{0 .y}=e^{x\th}\,e^{0.\eta}\,e^{0. y}\,.
  $$

Clearly one has the expansions equivalent to the basic (\peq{fexp}),
  $$
  f(x)=e^{x\eta}\,e^{0.\th}\,f(0)=e^{x\th}\,e^{0.\eta}\,f(0)\,.
  \eql{fexps}
  $$
  
A further useful relation is,
  $$
    f(D)\,0^r=f(\eta)\,e^{0.\th}\,0^r=f(\th)\,e^{0.\eta}\,0^r\,,
    \eql{cor}
  $$
a corollary of (\peq{inv}). This implies another (equivalent) expression,
  $$
   f(D)\,e^{0.t}=f(\eta)\,e^{0.\th}\,e^{0.t}=f(\th)\,e^{0.\eta}\,e^{0.t}\,.
    \eql{cor2}
  $$

\section{\bf8. Operator expansions}

From the function (possibly polynomial) expansion (\peq{fexp}), or (\peq{fexps}), follows a
generalised Taylor expansion. Expressed in operators, this is (\cf\ [\pref{Steffensen3}]
Equ.(30)),
  $$\eqalign{
   E^x&=\sum_{\nu=0}^\infty {x^{\{\nu\}}\,\th^\nu \over\nu!}=
   \sum_{\nu=0}^\infty {x^{\{\nu\}^{-1}}\,\eta^\nu \over\nu!}\cr
   &=e^{x\eta}\,e^{0.\th}=e^{x\th}\,e^{0.\eta}\,.
   }
   \eql{gtay}
  $$
The standard cases are conventionally discussed by Steffensen [\pref{Steffensen}] \S18
202,203.\mgn{Others}

It is sometimes convenient to remove the factor of $x$ that occurs in all poweroids, and ask
for expansions in terms of $x^{\{\nu+1\}-1}$. The operator form is obtained by
differentiating (\peq{gtay}) with respect to $D$ which gives
   $$\eqalign{
   E^x&={1\over x}\,e^{x\eta}\,0\,e^{0.\th}\th'\cr
   &=\sum_{\nu=0}^\infty {x^{\{\nu+1\}-1}\,\th^\nu \th'\over\nu!}\,,
   }
   \eql{gtay2}
  $$
[\pref{Steffensen3}], equ.(32).

Corresponding to (\peq{inv}) there are important operator expansions connecting multiple
derivatives and multiple differences. I derive these for general differences following the
forward difference treatment of Boole, [\pref{Boole}] p.24, which is an application of
Maclaurin's theorem in its secondary form,\footnote{It is very basic that Maclaurin's
theorem includes the Taylor expansion which is obtained by setting $t=D$ and $f(D)=E^h$
to give $ E^h=E^h\,e^{0.D}=e^{(0+h).D}=e^{hD}$.}
  $$
  f(t)=f(D)\,e^{0.t}\,,
  \eql{mt2}
  $$
where $D$ acts on $0$. This is a trivial consequence of the basic derivative,
$D^m\,0^n=m!\,\de_m^n$ (a Kronecker delta), an expression of duality.

$\th$ is a function of $D$ and so, setting $t=D$ and $f=\th^n$,
  $$
  \th^n=\th^n\,e^{0.D}\,,
  \eql{mcd}
  $$
where the $\th$ on the right--hand side acts on $0$. This is the required expression in
symbolic form. Expansion of the exponential gives a sum of powers of $D$,
  $$
  \th^n =\sum_{m=n}^\infty {\th^n0^m\over m!}\,D^m \,.
  \eql{mcdc}
  $$

Inversely, consider $f=D^m$ as a function of $\th$, and so set $t=\th$. According to
(\peq{mt2}), $\th$ in $D$ then has to be replaced by $D$. This gives the $\eta$ operator,
(\peq{xidef}), and the required symbolic expression is then,
  $$
   D^m=\eta^m\,e^{0.\th}\,,
   \eql{mcx}
  $$
yielding the generalised Newton series (\cf\ [\pref{Steffensen3}] Equ.(40)),
  $$
  D^m =\sum_{n=m}^\infty {\eta^m0^n\over n!}\,\,\th^m \,,
  $$
which could also be derived from (\peq{mcdc}) using (\peq{inv}) or from (\peq{gtay}) by
differentiating with respect to $h$. These expansions are classic for the standard delta
operators (\eg\ [\pref{Steffensen}] \S18 214).

\section{\bf9.  Umbral composition}

I remark that equation (\peq{cdf}) for the poweroid exhibits an umbral--like quality in that
the symbol, $0^n$, can be treated (legitimately) as a power. Jeffery's operator,
$e^{x\eta}$, linking the poweroid and $0^n$ acts as an `umbral operator' in this
approach. In [\pref{RandM}], this term signifies a (linear) operator on a polynomial
sequence that yields another such sequence and in [\pref{RandM}] an important property of
such operators is their composition. I  here consider that of $e^{x\eta_1}$ and
$e^{y\eta_2}$. The basic formula needed is again the Maclaurin expansion, (\peq{mt2}),
  $$
   f(\eta_2)=f(D)\,e^{0.\eta_2}
  $$
applied to,
  $$
  f(*)\equiv e^{x\eta_1(*)}
  $$
giving the composition,
  $$
  e^{x\eta_1(\eta_2)} =e^{x\eta_1}\,e^{0_1.\eta_2}\,,
  \eql{comp}
  $$
which shows, analogous to  [\pref{RandM}], that successive actions correspond to functional
composition, and form a group, say $\caG$.

It follows that the poweroid for the functionally composite operator,\break
$\phi_{12}(D)\equiv(\phi_2\circ\phi_1)(D)=\phi_2\big(\phi_1(D)\big)$, equals,
  $$\eqalign{
   x^{\{n\}_{12}}&=e^{x\eta_1}\,e^{0_1.\eta_2}\,0_2^n\cr
   &=e^{x\eta_1}\,0_1^{\{n\}_2}\cr
   &=x^{\{n\}_2}\big|_{x^i\to x^{\{i\}_1}}\,,\quad i=1,\ldots,n\,,
   }
   \eql{umbc}
  $$
which is referred to as `umbral composition'  in [\pref{RandM}] Theorem 6, which I have
therefore just proved. No explicit umbral notions are needed.

If $\phi_2$ is the inverse of $\phi_1$, the composite delta operator is simply $D$ which has
the ordinary power as its poweroid. Expressed in symbols
  $$
  x^{\{n\}^{-1}}\big|_{x^i\to x^{\{i\}}}=n^2\,,\quad \forall n\,,
  $$
so that the set of all poweroids forms a group under composition, with identity the sequence
of poweroids, $x^n$, $(n=0,1,\ldots)$. It is clear from (\peq{umbc}) that this group is
isomorphic to $\caG$.
\section{\bf10. Vector space interpretation. Interpolation.}

This similitude is no accident. It will be recognised that our development so far is nothing
more than an alternative symbolisation of the vector space approach to umbral calculus
advanced by Roman, [\pref{Roman}], following Rota. In fact Roman's `first umbral result',
[\pref{Roman}] Theorem 2.1.10, is just Theorem B of (\peq{JABCD}), which is probably
due to Herschel.

All relevant expressions take, can take, or include the form $ f(D)\,0^n$ and I have
manipulated $0^n$ as a power, performing summations, for example. Making an extension,
I replace $ f(D)\,0^n$ by $ f(D)\,h(0)$ where both $f(D)$ and $h(0)$ are formal power
series ($h$ a polynomial for ease) with, by ancient definition, $
f(D)\,h(0)=\big(f(D)h\big)(0)\equiv \big(f(D)\,h(x)\big)\big|_0$. The space \footnote{ I do
not use the correct terminology.}, $\caF$, of the $f(D)$ is generated by the set of divided
powers, $D^m/m!$ $(m=0,1,\ldots)$, while the space of the $h(x)$, $\caH$, has the basis
$x^n$ $(n=0,1,\ldots)$. The equation $D^m/m!\,0^n=\de^m_n$ (used before) allows a
duality to be set up between the two spaces with the bracket $\br{f(D)}{h(x)}\equiv
f(D)h(0)$, implying evaluation at zero. The previous equations can then be rewritten in
bra(c)ket form, \cf\ [\pref{Roman}], but I will not do so large scale. For example, the
expansion theorem, (\peq{fexp}), takes the form, (\cf\ [\pref{Roman}] Theorem 2.4.2),
  $$
   \ket{f(x)}={1\over\nu!}\ket{b^\th_\nu(x)}\br{\th^\nu(D_x)}{f(x)}\,.
  $$
or, expunging $x$ and $D$ as understood,\mgn{$\nu! $ absorbed?}
  $$
   \ket{f}={1\over\nu!}\ket{b^\th_\nu}\br{\th^\nu}{f}\,.
  $$
or even the abstract expression of Taylor expansion,
  $$
   \ket{b^\th_\nu}\bra{\th^\nu}=\nu!\,{\bf1}\,.
   \eql{abt}
  $$
Elements of $\caF$ are bras and elements of $\caH$ are kets and ${\bf 1}$ is the unit
operator in both spaces.

Of course, it is not necessary to use the $(D^m/m!,x^n)$ basis set. This is just a particular
case of the, equally dual, bases, $(\th^m/m!,b^\th_n)$, all having delta brackets,
$\br{\th^m}{b^\th_n}=m!\,\de^m_n$ (from (\peq{factact})).

A more symmetrical form of equation (\peq{abt}) arises in an approach to interpolation
espoused by Aitken, [\pref{Aitken2,Aitken3}]. (See also Curry, [\pref{Curry}], and
[\pref{Steffensen3}] \S\S 4,14), Aitken derives a generalised Gregory--Newton formula in
terms of a set of delta operators, $\th_i$, which, when all these are the same, is equivalent
to (\peq{fexp}), or (\peq{gtay}). His formalism involves an inverse operator to $\th$,
denoted by $\Th$, such that\footnote { Note that $\th$ and $\Th$ have the nonzero
commutator, $[\th,\Th]=L_0$.}
  $$
   \Th\th f(x)=f(x)-f(0)
  $$
or
  $$
  {\bf1} =L_0+\Th\th=\th\Th
   \eql{thinv}
  $$
where $L_0$ is the operator signifying evaluation at zero.\footnote{ This is analogous to the
definition of a Green function for an operator that has a zero mode, which is removed from
the delta function. $L_0$ is then the projection onto the zero mode.}

The iteration of (\peq{thinv}) produces,
  $$\eqalign{
  \bf 1&=L_0+\Th\th\cr
  &=L_0+\Th L_0\th+\Th^2\th^2\cr
  &=L_0+\Th L_0\th+\Th^2L_0\th^2+\Th^3\th^3\cr
  &=L_0+\Th L_0\th+\Th^2L_0\th^2+\Th^3L_0\th^3+\ldots\,.
  }
  $$

Since $L_0\th^\nu$ yields a constant, $\Th^\nu L_0\,\th^\nu=(\Th^\nu
x^0)\,L_0\,\th^\nu$, and Aitken's version of the interpolatory (\peq{fexp}) is then,
($1\equiv x^0$), \footnote{ For interpolation purposes, the $\th^\nu f(0)$ are taken as
numerical input data.  The series terminates for polynomials. Otherwise, as an
approximation, one can stop the sum at a finite point with an explicit expression for the
remainder, [\pref{Steffensen3}] equ. (161).}
  $$
    f(x)-f(0)=\sum_{\nu=1}^\infty \Th^\nu 1\,.\,   \th^\nu f(0)=
    \sum_{\nu=1}^\infty \Th^{\nu-1} x\,.\,   \th^\nu f(0)\,,
    \eql{aint}
  $$
$\Th^\nu 1$ being a polynomial of degree $\nu$ vanishing at $x=0$ and satisfying (\cf\
[\pref{Steffensen3}], equ.(29)).
  $$
  \th\,(\Th^\nu1)= \th\,\Th^\nu1=\Th^{\nu-1}1\,,
  $$
which shows that (\peq{aint}) is (\peq{fexp}) with $\Th^\nu1$ the divided poweroid,
$x^{\{\nu\}}/\nu!$. This can be quickly obtained by a formal backwards iteration of the
basic property, (\peq{act}), \ie $x^{\{\nu\}}=\nu\,\Th\,x^{\{\nu-1\}}$, leading to a
multiple summation.

The abstract completeness equation, (\peq{abt}), can therefore be compressed to a more
symmetrical looking,
  $$\eqalign{
  {\bf 1}&=\sum_{\nu=0}\ket{\Th^{\nu}}\bra {\th^{\nu}}\cr
   &=\ket1\bra1+\sum_{\nu=1}\ket{\Th^{\nu}}\bra {\th^{\nu}}\cr
   }
   \eql{ggn}
  $$
where $\ket1\bra1,\,=L_0,$ projects onto the (constant) zero mode, $\th1=0$.

Equation (\peq{ggn}) is my formal expression of a generalised Gregory--Newton
interpolation.
\section{\bf11. Specific expansions}
As Jeffery remarks, [\pref{Jeffery}] \S10, knowledge of the factorial numbers is useful in
many expansions.

From the explicit definition of $\eta$, (\peq{xidef}), the expansion of the powers of the
inverse $\phi$ function, whatever this is, is contained in,
  $$
   \big(\phi^{-1}(t)\big)^n=\eta^n\,e^{0.\,t}\,,
   \eql{invsh}
  $$
showing rapidly that the coefficients are just the factorial numbers of the first kind.

A more complicated example is,
  $$\eqalign{
  f\bigg(\phi^{-1}\big(\phi^{-1}(x)\big)\bigg)&=
  f\big(\phi^{-1}(\eta)\big)\,e^{0.x}\cr
  &=f\big(\phi^{-1}(D)\big)\,e^{0.\eta}\,e^{0.x}\cr
  &=f(\eta)\,e^{0.\eta}\,e^{0.x}\cr
  &=f(\eta)\,\bigg(0^{\{1\}}{x\over1!}+0^{\{2\}}\,{x^2\over2!}+\ldots\bigg)\,,
  }
  \eql{exp2}
  $$
where (\peq{cor}) and (\peq{cdf}) have been used to give the second and fourth lines. This
process can be continued.

As an example, consider
   $$\eqalign{
  \phi^{-2}(x)\equiv\phi^{-1}\big(\phi^{-1}(x)\big)
  =\eta\,0^{\{1\}}\cdot{x\over1!}+\eta\,0^{\{2\}}\cdot{x^2\over2!}+\ldots\,.
  }
  \eql{exp3}
  $$
and require the coefficient of $x^r/r!$, \ie
  $$
  \eta\,0^{\{r\}}=\sum_{s=0}^r{\eta\,0^s\over s!}\cdot\eta^s\,0^r\,.
  $$
which can be evaluated straightforwardly given the generalised differentials of nothing. (This
is really going back to the third line of (\peq{exp2}). A simple iteration of (\peq{invsh})
also yields the same result.)

Setting $x=\th$ gives the explicit expression for $\eta$ as a formal power series in $\th$,
  $$
  \eta=\eta\,0^{\{1\}}\cdot{\th\over1!}+\eta\,0^{\{2\}}\cdot{\th^2\over2!}+\ldots
  \eql{etapow}
  $$
and interchanging $\th$ and $\eta$ gives the inverse series,
  $$
  \th=\th\,0^{\{1\}^{-1}}\cdot{\eta\over1!}+\th\,0^{\{2\}^{-1}}\cdot{\eta^2\over2!}+\ldots
  $$

\section{\bf12. Eigenfunctions as generating fumctions. Another approach}

In calculus, the defining characteristic of the exponential is that it is reproduced upon
differentiation. With this in mind, one can ask for the eigenfunctions of the difference
operator, $\th$, \ie for functions $e(x,t)$ satisfying the difference equation,
  $$
   \th\,e(x,t)=t\,e(x,t)\,,
  $$
where $\th$ acts on $x$.

Proceeding in a standard way, as a trial solution assume $e(x,t)=\rho^x(t)$ (\eg\ Boole,
[\pref{Boole}]). Then, the form of $\th$, (\peq{psi}), gives the formal solution
  $$
  \rho(t)=e^{\phi^{-1}(t)}
  $$
 and so,
   $$\eqalign{
    e(x,t)&=e^{x\,\phi^{-1}(t)}\cr
    &=e^{x\eta}\,e^{0. t}\,.
    }
    \eql{eigf}
    $$
This last follows from the basic relation (say as an extension of (\peq{invsh}) or directly),
  $$
  f\big(\phi^{-1}(t)\big)=f(\eta)\,e^{0. t}\,.
  \eql{fxi}
  $$

In this approach, the poweroid $x^{\{n\}}$ is {\it defined} by,
  $$
  x^{\{n\}}\equiv e^{x\eta}\,0^n\,,
 \eql{fact3}
  $$
 because it {\it follows} from the eigenfunction relation,
   $$
   \th\,e^{x\eta}\,e^{0. t}=t\,e^{x\eta}\,e^{0. t}\,,
   $$
by expansion in $t$ that,
  $$
  \th\,x^{\{n\}}=n\,x^{\{n-1\}}\,,
  $$
the basic feature of the poweroid. It is seen that the eigenfunction is a poweroid generating
function.

The same analytical result, derived differently is given as Corollary 4 in [\pref{RKO}] p.693,
and Corollary 2 in [\pref{RandM}] p.189. See also Roman, [\pref{Roman}].

The other two properties result immediately from (\peq{fact3}), {\it viz.},
  $$
   0^{\{n\}}=0^n=\de_n^0\,\quad{\rm and}\quad x^{\{1\}}=(1+x)0^1=x\,.
  $$
  
Interchanging $\th$ and $\eta$ gives the eigenfunction of $\eta$  as $e^{x\th}\,e^{0.
t}=e^{x\phi(t)}$ which is the generating function of the inverse poweroids,
$x^{\{n\}^{-1}}$.
  
\noin  {\it Example. Central case.}

As a standard, but non--trivial example, I {\it derive} the central factorial, $x^{[n]}$, from
(\peq{fact3}). Then, specifically,
  $$
  \eta=2\sinh^{-1}{D\over2}\,,
  $$
and
  $$
   e^{x\eta}=\bigg({D\over2}+\big(1+{D^2\over4}\big)^{1/2}\bigg)^{2x}\,.
  $$
  
According to (\peq{fact3}), the expansion of $e^{x\eta}$ in powers of $D$ is required, in
particular the $n$th power.\footnote{ The forwards case is just the binomial expansion and
is given in Jeffery, [\pref{Jeffery}].} The direct analytical calculation is given by Hansen,
[{\pref{Hansen}] \S18, and I copy it out here as the reference is slightly obscure.

For convenience set $D/2=u$ and $\la=e^{x\eta}\big|_{D=2u}$. Differentiation with
respect to $u$ yields the differential equation,
  $$
   (1+u^2){d^2\la\over du^2}+u{d\la\over du}-4x^2\la=0\,.
   \eql{de}
  $$
Assuming the power series solution,
  $$
   \la=\sum_{m=0}^\infty {M_m\over m!}\, u^m\,,
  $$
substitution into (\peq{de}) easily yields the simple recursion,
  $$
   M_{m}=\big(4x^2-(m-2)^2\big)\,M_{m-2}
  $$
with the initial values, $M_1=2x$ and $M_2=x^2$.

Picking out the coefficient of $D^n/n!$, equation (\peq{fact3}) gives the established
factorial, $x^{[n]}$, \eg\ Steffensen, [\pref{Steffensen}], Thiele, [\pref{Thiele}].

An alternative method obtains the central factorial from the forwards one (assumed given)
by a translation to the midpoint (after removal of the universal factor of $x$) -- a
consequence of the relation between difference operators, $\de=E^{-1/2}\,\De$. This is
formalised in the more general Rodrigues formula of Rota and Mullin, [\pref{RandM}],
[\pref{RKO}]. The Gould polynomial poweroid, (\peq{fact}), can also be computed in this
way. Steffensen, [\pref{Steffensen3}], gives an equivalent derivation.

\section{\bf13. The binomial property}

From the definition (\peq{fact3}) quickly follows the important, sometimes considered
defining, binomial property of the poweroids (basic polynomials).

In terms of the eigenfunction, (\peq{eigf}),  this is contained in the exponent property,
(obtained straightaway from (\peq{eigf})),
  $$
   e(x,t)\,e(y,t)=e(x+y,t)\,,
   \eql{expnt}
  $$
after application of the easy symbolic theorem, (Horner, [\pref{Horner}], \S6),
  $$
  F(\th)e^{0. t}\times G(\th)e^{0. t}=F(\th)\times G(\th) e^{0 .t}\,.
  \eql{horn}
  $$
$\eta$ is a function of $\th$, $\eta=\phi^{-2}(\th)$. (The power series has been given,
(\peq{etapow})). Or one can simply replace $\th$ by $\eta$ according to the previously
explained general rule.

The polynomials appear on expansion in $t$. The binomial property is explicitly, and neatly,
expressed in the coefficient theorem obtained from (\peq{horn}),
  $$
   F(\eta)0^n \times G(\eta) 0^n=\big(F(\eta)\times G(\eta')\big) \big(0+0')^n\,.
  $$
(Horner, [\pref{Horner}], \S7) setting $F(\eta)=e^{x\eta}$ and $G(\eta)=e^{y\eta}$.
Here $\eta'$ and $0'$ are symbols equivalent to $\eta$ and $0$ in the manner fully
explained by Herschel, [\pref{Herschel,Herschel3}], long before umbral methods  or
Aronhold's notation in  invariant theory. In later language, $0$ and $0'$ are exchangeable,
or similar, umbrae. These early works also include the obvious multinomial extension via
$(0+0'+0''+\ldots)^n$.

The binomial property statement (for poweroids) is due to Steffensen. Its converse was
proved by Rota and Mullin, [\pref{RandM}]. I deal with the converse in the present
formalism which is cosmetically rapid.

The converse amounts to firstly being given the binomial relation, (\peq{expnt}), for a {\it
general} polynomial generating function, $e(x,t)$. Then an operator, $\th$, is {\it defined}
so that $e(x,t)$ is an eigenfunction, $\th e(x,t)=te(x,t)$, and finally one just needs to show
that $\th$ commutes with $D$ or equivalently with $E$. This  follows smartly as follows.
  $$\eqalign{
  \th E^y e(x,t)&=\th e(x+y,t)=\th e(x,t). e(y,t)=te(x,t)e(y,t)\cr
  &=te(x+y,t)=E^y te(x,t)=E^y\th\, e(x,t)\,.
  }
  $$
This commutation relation can be extended by linearity to arbitrary polynomials (or power
series) and so the required operator statement, $\th E^y=E^y\th$, follows.

\section{\bf14. The connection constants}

The main concern in [\pref{RandM}] was the computation of the  connection between two
different sets of basic polynomials.\footnote{ One might term these `transformation
coefficients' as in vector space theory, used, say, in quantum mechanics.} Between those
for $\th_1$ and $\th_2$, this is defined by the linear relation,
  $$\eqalign{
  e^{x\eta_1}0^n&=e^{x\eta_2}\sum_m 0^m\,\rbr {m;2} {n;1}\cr
  &\equiv e^{x\eta_2}\, 0^{\{n\}_3}
  }
  $$
which defines a new delta operator, $\th_3$. Looking back at the composition statement,
(\peq{umbc}), one deduces the relation $\th_1=\th_3\circ\th_2$ so that
$\th_3=\th_1\circ\th_2^{-1}$ whose poweroid can hence be calculated and  its coefficients
then read off to give the connection constants, $\rbr {m;2} {n;1}$. This reproduces the
conclusion of [\pref{RandM}] quite neatly.
  
\section{\bf15. Representative Notation}

Although my position is that representative notation (classical umbral calculus) is
unnecessary, it is illuminating to make use of it, the main advantage, for me, being its
notational convenience and suggestive power.

The essential points are well known and are as follows. Elements of the space of operators
(\ie formal power series in $D$) are written as $e^{\al\,D}$ in terms of the umbra $\al$
with the coefficients, $\al_n$, of the power series being represented by $\al^n$. The delta
operator $\th$ is then $e^{\al\,D}-1$ which is the forwards operator with umbral step
$\al$. Likewise, the inverse \footnote {Roman refers to this as the {\it conjugate}
operator.}, $\eta$, is, umbrally, $e^{\ol\al\,D}-1$ which defines the inverse or `conjugate'
umbra, $\ol\al$. Naturally, the coefficients of $\eta$ are obtained from those of $\th$ by
inversion.

The power polynomial basis, $x^n$, is special in that it allows (some) power series to be
exponentially summed, as exemplified by the expansion of the factorial, (\peq{cdf}). The
extension of this desirable feature to, say, (\peq{pow}) can be achieved by the use of
umbrae. Riordan, [\pref{Riordan}], uses this tempting device and in this section, for
completeness, I link it to the Jeffery operator formalism (in general form). Ray,
[\pref{Ray2,Ray3}], also employs this symbolism and gives useful summaries.

The representative, $g$, is heuristically introduced by the umbral equality,
   $$
    g^n(x)=g_n(x)\equiv\sum_{\nu=0}^n g(n,\nu)\,x^\nu
    =x^{\{n\}}=e^{x\eta}\,0^n\,,
    \eql{gee}
   $$
so that the definition, (\peq{pow}), can now be summed,
  $$
    x^n=e^{g\,\th}\,0^n\,,\quad g=g(x)\,.
    \eql{def2}
  $$
I also define its reciprocal, $\ol g$, by, \cf\ [\pref{Riordan}],
  $$\eqalign{
    {\ol g}^n(x)=\ol g_n(x)\equiv\sum_{\nu=0}^n \ol g(\nu,n)\,x^\nu
    =e^{x\th}\,0^n=x^{\{n\}^{-1}}
    }
    \eql{Gee}
   $$
employing (\peq{not}).

I can then umbrally sum (\peq{gee}) and (\peq{Gee}) to exponentials (for example),
  $$
e^{y\,g(x)}=e^{x\eta}\,e^{0 .y}
 \eql{cons}
 $$
 $$
  e^{y\,\ol g(x)}=e^{x\th}\,e^{0. y}\,,
  \eql{cons2}
  $$
where $x$ and $y$ can be complex numbers.

From the definitions of the operators $\eta$ and $\th$, the {\it generating functions},
(\peq{cons}) and (\peq{cons2}), take the explicit forms,
  $$
    e^{y\,g(x)}= e^{x\phi^{-1}(y)}\,,
     \eql{Bl4}
  $$
and
  $$
     e^{y\,\ol g(x)}=e^{x\phi(y)}\,,
     \eql{Bl3}
  $$
the two expressions being related by reciprocity as I now show.

I first derive the umbral statement of reciprocity. Umbrally, from  Eq.(\peq{cons2}) by
setting $x\to g(x)$ one arrives at the formally neat expression of
orthogonality/completeness,
   $$
    \ol g\circ g=1
    \eql{recip}
   $$
(This also easily results from (\peq{Gee}) followed by (\peq{gee}), which gives $({\ol
g}^n\circ g)(x)=x^n$.)

The reciprocal relation follows on iteration of (\peq{recip}),
  $$
    g\circ \ol g\circ g=g\,,
    \eql{recip2}
   $$
or, invoking completeness,
  $$
    g\circ \ol g=1\,.
    \eql{recip3}
   $$
One might therefore formally set $\ol g=g^{-1}$, \etc
 
Now, making the replacement $x\to \ol g(x)$ and employing (\peq{recip3}) turns
(\peq{Bl4}) into (\peq{Bl3}) as promised.

Eqs. (\peq{recip}) and (\peq{recip3}) are the classic umbral equivalents of the explicit
(\peq{inv}).

\section{\bf16. Explicit expansions}

In particular cases, Eq.(\peq{Bl4}) is identical with  classical expansions. Adjusting notation,
(\peq{Bl4}), under $y\to\phi(x)$ and $x\to r$, becomes, neatly,
  $$
     e^{rx}=e^{g\phi(x)}\,,\quad g=g(r)\,,
     \eql{Bl5}
  $$
which is the umbral representation of the expansion of the exponential in powers of $\phi$,
{\it viz.},
  $$
  e^{rx}=\sum_{n=0}^\infty {r^{\{n\}}\over n!}\,\phi^n(x)\,.
  \eql{expexp}
  $$
In the central case, when $\phi(x)=2\sinh(x/2)$, on choosing $r$ to be an integer, $n$, and
$x=i\phi$, some known, and very old, trigonometric expansions  for $\cos n\phi$ and $\sin
n\phi$ have thus been obtained with very little effort.\footnote{ These expansions, which
are valid for any $r$, are discussed by Sheppard, [\pref{Sheppard}], in the context of
difference calculus.} In the forwards and backwards cases, there only results an identity,
  $$
   e^{rx}=\sum_{n=0}^\infty{r^{(n)}\over n!}\,(e^{x}-1)^n=e^{rx}\,.
  $$

As a check, or as an illustration of the circularity of the relations, setting $x=D$ in
(\peq{expexp}) reproduces the operator expansion (\peq{gtay}), the generalised
Gregory--Newton equation. For the forwards case see Roman, [\pref{Roman}], p.58.
\begin{ignore}
which produces more complicated expansions for $\cosh rx$ and $\sinh rx$. , \eg,
  $$
  \cosh rx=\sum_{n=0}^\infty {2^n\,r^{(n)}\over n!}
  \big(\cosh (nx/2)\,P_n+\sinh( nx/2)\,(1-P_n)\big)\sinh ^n(x/2)\,,
  $$
where $P_b=(1+(-1)^n)/2$ projects onto even $n$.
\end{ignore}

Related to (\peq{Bl5}), I give a further illustration of the classic umbral notation.  Either as
an extension of (\peq{invsh}) or directly,
  $$
  f\big(\phi^{-1}(x)\big)=f(\eta)\,e^{0.x}=\sum_{s=0}^\infty\,
  {f(\eta)\,0^s\over s!}\,x^s\,.
  \eql{fxi2}
  $$
The umbral, primary form of Maclaurin's theorem is,
  $$
  f(\eta)=\sum_{s=0}^\infty{f^{(s)}(0)\over s!}\,\eta^s=e^{f\,\eta}\,,
  $$
where the power $f^n$ is to be replaced by the derivative, $f^{(n)}(0)$, after the
expansion of the exponential. Then the elegant result,
  $$\eqalign{
  f\big(\phi^{-1}(x)\big)&=\sum_{n=0}^\infty
  {e^{f\,\eta}\,0^n\over n!}\,x^n\cr
  &=\sum_{n=0}^\infty{f^{\{n\}}\over n!}\,x^n\,,\cr
  }
  \eql{el}
  $$
follows from (\peq{fxi}) on using (\peq{cdf}). One could write this result as,
  $$
 f(\eta)\,e^{0. x} =e^{fx}\,,
 \eql{snap}
  $$
where now the power, $f^r$, first represents the poweroid, $f^{\{r\}}$, in which the
powers of $f$ are still to be replaced by  the derivatives, $f^{(n)}(0)$, as above.

Equation (\peq{expexp}) is an example of (\peq{el}) in view of the symbolic relation
  $$
  \big(e^{rx}\big)^{\{n\}}=r^{\{n\}}
  $$

Roman, [\pref{Roman}], gives many examples of expansions for various delta operators
and associated polynomials.
 
\section{\bf17. Conclusion}
It has not been my intention to rederive all the results obtained by the modern umbral
calculus symbolisation. This would be a waste of effort. However, I think I have
demonstrated the pertinence of the older formalism, with some of the derivations going
through more smoothly than their modern versions. I might cite the connection constant
result of \S14, the binomial property of \S13 and the derivation of the Rodrigues--type
equation in \S5.

As a continuation, the combinatorial questions analysed by Jeffery, [\pref{Jeffery}], could
be extended to the other factorials, and their roots.
%\newpage
\vglue 15truept

\noin{\bf References.} \vskip5truept

\begin{putreferences}
   \ref{BaandB}{Balian,R. and Bloch,C. \aop{60}{1970}{401}.}
   \ref{Sheppard}{Sheppard,W.F. \plms{31}{1899}{449}.}
   \ref{Ray2}{Ray,N. \aim{61}{1986}{49}.}
   \ref{Ray3}{Ray,N. \tams{309}{1988}{191}.}
   \ref{RKO}{Rota,G-.C.,Kahaner,D. and Odlyzko,A. {\it J.Math.Anal.Appl.} {\bf42} (1973)
    684.}
   \ref{L}{Levi,D.,Tempesta,P. and Winternitz,P. \jmp{45}{2004}{4077}.}
   \ref{Hansen}{Hansen,P.A. {\it Abh.K\"on.S\"achs.Gesellsch.} {\bf 11} (1865) 505.}
   \ref{Curry}{Curry,H.B. {\it Port.Math.} {\bf10} (1951) 1.}
   \ref{Roman}{Roman,S. {\it The Umbral Calculus} (Academic Press, New York, 1984).}
   \ref{RandR}{Roman,S. and Rota,G-.C. \aim{27}{1978}{95}.}
   \ref{RandM}{Rota,G-.C. and Mullin,R {\it On the foundations of Combinatorial Theory, III}
   in {\it Graph Theory and its Applications} (Academic Press, New York, 1970).}
   \ref{Dowcen}{Dowker,J.S. {\it Central Differences, Euler numbers and
   symbolic methods} ArXiv: 1305.0500.}
   \ref{Aitken}{Aitken,A.C. {\it J. Inst. Actuaries} {\bf 64} (1933) 449.}
   \ref{Aitken2}{Aitken,A.C. \pems {1 } {1929} {199}.}
   \ref{Aitken3}{Aitken,A.C. {\it J. Inst. Actuaries} {\bf 61} (1930) 107.}
   \ref{Steffensen2}{Steffensen,J.F. {\it J. Inst. Actuaries} {\bf64 } (1933) 165.}
   \ref{Steffensen3}{Steffensen,J.F. \am{73}{1941}{333}.}
   \ref{Zhang}{Zhang,W.P. {\it The Fibonacci Quarterly} {\bf 36} (1998) 154.}
   \ref{Michel}{Michel, J.G.L. {\it J. Inst. Actuaries} {\bf 72} (1946) 470.}
   \ref{LiandZ}{Liu,G.D. and Zhang,W.P. {\it Acta Math. Sinica} {\bf24} (2008) 343.}
   \ref{JandS}{Joliffe, A.E. and Shovelton, S.T. \plms{13}{1913}{29}.}
   \ref{CaandR}{Carlitz,L. and Riordan,J. \cjm {15}{1963}{94}.}
   \ref{RandT2}{Rota, G-.C. and Taylor, B.D. {\it SIAM J.Math.Anal.} {\bf 25} (1994) 694.}
   \ref{Bickley}{Bickley, W.G. {\it J. Math. and Phys.}(MIT) {\bf 27} (1948) 183.}
   \ref{Schwatt}{Schwatt, I.J. {\it An introduction to the operations with series}
    (Univ.Pennsylvania Press, Philadelphia, 1924).}
    \ref{Loney}{Loney, S.L. {\it Plane Trigonometry} (CUP, Cambridge, 1893).}
   \ref{Teixeira}{Teixeira, F.G. \jram{116}{1896}{14}.}
   \ref{Riordan}{Riordan,J. {\it Combinatorial Identities} (Wiley, New York, 1968).}
   \ref{Charal}{Charalambides, Ch.A. {\it The Fibonacci Quarterly} {\bf 19.5} (1981) 451.}
    \ref{Stern}{Stern,W. \jram {79}{1875}{67}.}
    \ref{Thiele}{Thiele,T.N. {\it Interpolationsrechnung} (Teubner, Leipzig, 1909).}
    \ref{Dowren}{Dowker,J.S. \jpamt {46}{2013}{2254}.}
    \ref{KnandB}{Knuth,D.E. and Buckholtz,T.J. {\it Math.Comp.} {\bf 21} (1967) 663.}
    \ref{BSSV}{Butzer,P.L., Schmidt,M., Stark,E.L. and Vogt,I. {\it Numer.Funct.Anal.Optim.}
    {\bf 10} (1989) 419.}
    \ref{Ely}{Ely,G.S. \ajm{5}{1882}{337}.}
    \ref{BrandH}{Brent,R.P. and Harvey,D., {\it Fast Computation of Bernoulli, Tangent and
    Secant numbers}, ArXiv:1108.0286.}
   \ref{BaandBe}{Bayad and Beck  ArXiv.}
   \ref{Blissard}{Blissard, \qjm{}{}{}.}
   \ref{Gould2}{Gould,H.W. \amm{79}{1972}{} .}
   \ref{Worpitzky}{Worpitsky,J. \jram{}{}{}.}
   \ref{Saalschutz}{Saalsch\"utz,L. {\it Vorlesungen \"uber die Bernoullischen Zahlen.}
   (Springer--Verlag, Berlin, 1893).}
   \ref{Workman}{Workman, W.P. {\it Memoranda Mathematica}, (OUP, Oxford, 1912).}
   \ref{Nielsen}{Nielsen,N. {\it Handbuch der Theorie der Gammafunktion} (Teubner,
   Leipzig, 1906).}
   \ref{Steffensen}{Steffensen,J.F. {\it Interpolation}, (Williams and Wilkins,
    Baltimore, 1927).}
    \ref{Joffe}{Joffe,S.A. \qjm{47}{1916}{103}.}
    \ref{Shovelton}{Shovelton,S.T. \qjm{46}{1915}{220}.}
   \ref{Horner}{Horner,J. \qjm{4}{1861}{111}.}
   \ref{Jeffery}{Jeffery, H.M. \qjm{4}{1861}{364}.}
   \ref{Jeffery2}{Jeffery, H.M. \qjm{5}{1862}{91}.}
   \ref{Jeffery3}{Jeffery, H.M. \qjm{6}{1864}{179}.}
   \ref{Grunert}{Grunert, \jram{25}{1843}{240}.}
   \ref{AgandD}{Agoh,T. and Dilcher,K. {\it J. Number Theory} {\bf 124} (2007) 105.}
   \ref{Franssens}{Franssens,G.R. {\it Analysis Mathematica} {\bf 33} (2007) 17.}
   \ref{Andrews}{Andrews,G.E. {\it Ramaujan J.} {\bf 7} (2003) 385.}
    \ref{Glaisher2}{Glaisher,J.W.L. \qjm {40}{1909}{275}.}
    \ref{RandF}{Rubinstein, B.Y. and Fel,L.G., {\it Ramanujan J.} {\bf11}(2006)331.}
    \ref{Rubinstein}{Rubinstein, B.Y. {\it Ramanujan J.} {\bf15}(2008)177.}
   \ref{Dickson3}{Dickson,L.E. {\it History of the Theory of Numbers} vol.2.
   (Carnegie Institute, Washington, 1920).}
   \ref{Brioschi}{Brioschi,F. {\it Annali di sc. mat. e fis.} {\bf8} (1857) 5.}
   \ref{Fedosov}{Fedosov,B.V. {\it Sov.Math.Dokl.} {\bf 5} (1963) 1992}
   \ref{Chapman}{Chapman}
    \ref{Jordan}{Jordan,C. {\it Calculus of Finite Differences}, (Budapest, 1939).}
   \ref{PandB}{van der Pol, B. and Bremmer,H. {\it Operational Calculus} (CUP,
   Cambridge, 1964).}
   \ref{TandD}{Tauber,S. and Dean,D. {\it J.SIAM} {\bf 8} (1960) 174.}
   \ref{Traub}{Traub,J.F. {\it Math. of Comp.} {\bf 19} (1965) 177.}
   \ref{Chapman}{Chapman,S. \plms{9}{1911}{369}. }
   \ref{Moore}{Moore,D.H. \amm{69}{1962}{132}.}
   \ref{HaandR}{Hardy,G.H. and Riesz,M. {\it General theory of Dirichlet's series}
   (CUP, Cambridge, 1915).}
   \ref{Gupta}{Gupta,H. {\it Tables of partitions} (Royal.Soc.Math.Tables 4)
   (Cambridge, 1958).}
   \ref{Bromwich}{Bromwich, T.J.I'A. {\it Infinite Series},
  (Macmillan, London, 1947).}
  \ref{Hobson2}{Hobson,E.W. {\it Theory of functions of a real variable}, vol.2.
  (C.U.P., Cambridge,1907).}
   \ref{Bromwich2}{Bromwich, T.J.I'A. {\it Infinite Series}, 1st Edn.
  (Macmillan, London,1907).}
   \ref{Knopp3}{Knopp,K. {\it Theory of Infinite Series} (Blackie, London, 1928).}
   \ref{Vandiver}{Vandiver,H.S. \tams{51}{1942}{502}.}
   \ref{Adamchik}{Adamchik,V.S.{\it Appl.Math. and Comp.} {\bf 187} (2007) 3.}
   \ref{Cvijovic}{Cvijovi\'{c},D. {\it Appl.Math. and Comp.} {\bf 215} (2009) 3002.}
   \ref{Cvijovic1}{Cvijovi\'{c},D. {\it Comp.and Math. with Appl.} {\bf 62} (2011) 1879.}
   \ref{Carlitz2}{Carlitz,L. {\it Math.Magazine} {\bf 32} (1959) 247.}
   \ref{Israilov}{Israilov.M.I. {\it Sibirsk Mat. Zh.} {\bf 22} (1981) 121.}
   \ref{SandZ}{Sills,A.V. and Zeilberger,D. {\it Formulae for the number of partitions of
    $n$ into at most $m$ parts (using the quasi-polynomial ansatz)} ArXiv:1108.4391.}
   \ref{BandV}{Brion,M. and Vergne,M. {\it J.Am.Math.Soc.} {\bf 10} (1997) 371.}
   \ref{Knopf}{Knopf,P.M. {\it Math. Magazine} {\bf 76} (2003) 364.}
   \ref{Gould}{Gould,H.W. {\it Am. Math. Monthly} {\bf 86} (1978) 450.}
   \ref{Hoffman}{Hoffman,M.E. \amm{102}{1995}{23}.}
   \ref{Hoffman2}{Hoffman,M.E. {\it Electronic J.Comb.} {\bf 6} (1999) \#R21.}
   \ref{Boyadzhiev}{Boyadzhiev, K.N. {\it Fibonacci Quarterly} {\bf45} (2008) 291.}
   \ref{Sills}{Sills}
   \ref{Munagi}{Munagi,A.O. {\it Electronic J. Comb. Number Th.} {\bf 7} (2007) \#A25.}
   \ref{Dowsyl}{Dowker,J.S. {\it Relations between Ehrhart polynomials, the heat kernel
   and \break Sylvester waves} ArXiv:1108.1760}
   \ref{Alfonsin}{Alfonsin,J.L.R. {\it The Diophantine Frobenius Problem} (O.U.P.,
   Oxford, 2005).}
   \ref{GeandS}{Gel'fand,I.M. and Shilov,G.E. {\it Generalized Functions} Vol.1
   (Academic Press, N.Y., 1964).}
    \ref{Boole}{Boole, G. {\it Calculus of Finite Differences}, (MacMillan, Cambridge,
1860).}
   \ref{Boole2}{Boole, G. {\it Calculus of Finite Differences}, 2nd Edn. (MacMillan, Cambridge,
1872).}
     \ref{Cayley5}{Cayley, A.  {\it Phil. Trans. Roy. Soc. Lond.} {\bf 148} (1858) 47.}
     \ref{Milne-Thomson}{Milne-Thomson, L.M. {\it The Calculus of Finite Differences},
     (MacMillan, London, 1933).}
    \ref{Herschel}{Herschel, J.F.W.  {\it Phil. Trans. Roy. Soc. Lond.} {\bf 106} (1816) 25.}
    \ref{Herschel2}{Herschel, J.F.W.  {\it Phil. Trans. Roy. Soc. Lond.} {\bf 108} (1818) 144.}
     \ref{Herschel4}{Herschel, J.F.W.  {\it Phil. Trans. Roy. Soc. Lond.} {\bf 140} (1850) 399.}
      \ref{Herschel5}{Herschel, J.F.W.  {\it Phil. Trans. Roy. Soc. Lond.} {\bf 150} (1860) 319.}
    \ref{Herschel3}{Herschel,J.F.W. {\it Examples of the Applications of the Calculus of Finite
    Differences} (Deighton, Cambridge, 1820).}
    \ref{Brinkley}{Brinkley {\it Phil. Trans. Roy. Soc. Lond.} {\bf } (1807) 114 .}
    \ref{Littlewood2}{Littlewood,D.E. {\it The Theory of Group Characters}
    (Clarendon Press, Oxford, 1950).}
    \ref{Wright}{Wright,E.M. \amm{68}{1961}{144}.}
\ref{Carlitz}{Carlitz,L. \dmj{27}{1960}{401}.}
     \ref{Netto}{Netto,E. {\it Lehrbuch der Combinatorik} 2nd Edn. (Teubner, Leipzig, 1927).}
    \ref{FdeB}{Fa\`{a} de Bruno, F. {\it Th\'eorie des Formes Binaires} (Brero, Turin,1876).}
    \ref{Ehrhart}{Ehrhart,E. \jram{227}{25}{1967}.}
    \ref{Bell}{Bell,E.T.\ajm{65}{1943}{382}.}
    \ref{BandR}{Beck,M. and Robins,S. {\it Computing the Continuous Discretely,}
    (Springer, New York, 2007).}
    \ref{BandR2}{Beck, M. and Robins,S. {\it Discrete and Comp. Geom.} {\bf 27}(2002) 443.}
    \ref{Harmer}{Harmer,M. {\it J.Australian Math.Soc.} {\bf 84}(2008)217.}
    \ref{RandF}{Rubinstein, B.Y. and Fel,L.G., {\it Ramanujan J.} {\bf11}(2006)331.}
    \ref{BGK}{Beck,M., Gessel, I.M. and Komatsu,T. {\it Electronic Journal of Combinatorics}
    {\bf8}(2001) 1.}
    \ref{Sylvester}{Sylvester,J.J. \qjpam{1}{1858}{81}.}
    \ref{Sylvester2}{Sylvester,J.J. \qjpam{1}{1858}{142}.}
    \ref{Sylvester3}{Sylvester,J.J. \ajm{5}{1882}{79}.}
    \ref{Sylvester4}{Sylvester,J.J. \plms{28}{1896}{33}.}
    \ref{Dowgta}{J.S.Dowker, {\it Group theory aspects of spectral problems on spherical
    factors}, ArXiv.Math.DG: 0907.1309.}
    \ref{BDR}{Beck,M., Diaz and Robins,S. {\it J.Numb.Theory} {\bf 96} (2002) 1.}
    \ref{PandS}{P\'{o}lya, G. and Szeg\H{o},G. {\it Aufgaben und Lehrs\"atze aus der Analysis}
    (Springer--Verlag, Berlin, 1925).}
    \ref{EOS}{Elizalde,E., Odintsov, S.D. and Saharian, A.A. \prD{79}{2009}{065023}.}
    \ref{Cavalcanti}{Cavalcanti,R.M. \prD{69}{2004}{065015}.}
    \ref{MWK}{Milton, K.A., Wagner,J. and Kirsten,K. \prD{80}{2009}{125028}.}
    \ref{EBM2}{Ellingsen,S.A., Brevik,I. and Milton,K.A. \prE{81}{2010}{065031}.}
    \ref{EBM}{Ellingsen,S.A., Brevik,I. and Milton,K.A. \prE{80}{2009}{021125}.}
    \ref{BEM}{Brevik,I., Ellingsen,S.A. and Milton,K.A. \prE{79}{2009}{041120}.}
    \ref{FKW}{Fulling,S.A, Kaplan L. and Wilson,J.H. \prA{76}{2007}{012118}.}
    \ref{Lukosz}{Lukosz,W, {\it Physica} {\bf 56} (1971) 109; \zfp{258}{1973}{99}
    ;\zfp{262}{1973}{327}.}
    \ref{Gromes}{Gromes, D. \mz{94}{1966}{110}.}
    \ref{FandK1}{Kirsten,K. and Fulling,S.A. \prD{79}{2009}{065019} .}
    \ref{FandK2}{Fucci,G. and Kirsten,K, JHEP (2011), 1103:016.}
    \ref{dowgjms}{Dowker,J.S. {\it Determinants and conformal anomalies
    of GJMS operators on spheres}, ArXiv: 1007.3865.}
    \ref{Dowcascone}{dowker,J.S. \prD{36}{1987}{3095}.}
    \ref{Dowcos}{dowker,J.S. \prD{36}{1987}{3742}.}
    \ref{Dowtherm}{Dowker,J.S. \prD{18}{1978}{1856}.}
    \ref{Dowgeo}{Dowker,J.S. \cqg{11}{1994}{L55}.}
    \ref{ApandD2}{Dowker,J.S. and Apps,J.S. \cqg{12}{1995}{1363}.}
   \ref{HandW}{Hertzberg,M.P. and Wilczek,F. {\it Some calculable contributions to
   Entanglement Entropy}, ArXiv:1007.0993.}
   \ref{KandB}{Kamela,M. and Burgess,C.P. \cjp{77}{1999}{85}.}
   \ref{Dowhyp}{Dowker,J.S. \jpa{43}{2010}{445402}; ArXiv:1007.3865.}
   \ref{LNST}{Lohmayer,R., Neuberger,H, Schwimmer,A. and Theisen,S.
   \plb{685}{2010}{222}.}
   \ref{Allen2}{Allen,B. PhD Thesis, University of Cambridge, 1984.}
   \ref{MyandS}{Myers,R.C. and Sinha,A. {\it Seeing a c-theorem with
   holography}, ArXiv:1006.1263}
   \ref{MyandS2}{Myers,R.C. and Sinha,A. {\it Holographic c-theorems in
   arbitrary dimensions},\break ArXiv: 1011.5819.}
   \ref{RyandT}{Ryu,S. and Takayanagi,T. JHEP {\bf 0608}(2006)045.}
   \ref{CaandH}{Casini,H. and Huerta,M. {\it Entanglement entropy
   for the n--sphere},\break arXiv:1007.1813.}
   \ref{CaandH3}{Casini,H. and Huerta,M. \jpa {42}{2009}{504007}.}
   \ref{Solodukhin}{Solodukhin,S.N. \plb{665}{2008}{305}.}
   \ref{Solodukhin2}{Solodukhin,S.N. \plb{693}{2010}{605}.}
   \ref{CaandW}{Callan,C.G. and Wilczek,F. \plb{333}{1994}{55}.}
   \ref{FandS1}{Fursaev,D.V. and Solodukhin,S.N. \plb{365}{1996}{51}.}
   \ref{FandS2}{Fursaev,D.V. and Solodukhin,S.N. \prD{52}{1995}{2133}.}
   \ref{Fursaev}{Fursaev,D.V. \plb{334}{1994}{53}.}
   \ref{Donnelly2}{Donnelly,H. \ma{224}{1976}{161}.}
   \ref{ApandD}{Apps,J.S. and Dowker,J.S. \cqg{15}{1998}{1121}.}
   \ref{FandM}{Fursaev,D.V. and Miele,G. \prD{49}{1994}{987}.}
   \ref{Dowker2}{Dowker,J.S.\cqg{11}{1994}{L137}.}
   \ref{Dowker1}{Dowker,J.S.\prD{50}{1994}{6369}.}
   \ref{FNT}{Fujita,M.,Nishioka,T. and Takayanagi,T. JHEP {\bf 0809}
   (2008) 016.}
   \ref{Hund}{Hund,F. \zfp{51}{1928}{1}.}
   \ref{Elert}{Elert,W. \zfp {51}{1928}{8}.}
   \ref{Poole2}{Poole,E.G.C. \qjm{3}{1932}{183}.}
   \ref{Bellon}{Bellon,M.P. {\it On the icosahedron: from two to three
   dimensions}, arXiv:0705.3241.}
   \ref{Bellon2}{Bellon,M.P. \cqg{23}{2006}{7029}.}
   \ref{McLellan}{McLellan,A,G. \jpc{7}{1974}{3326}.}
   \ref{Boiteaux}{Boiteaux, M. \jmp{23}{1982}{1311}.}
   \ref{HHandK}{Hage Hassan,M. and Kibler,M. {\it On Hurwitz
   transformations} in {Le probl\`eme de factorisation de Hurwitz}, Eds.,
   A.Ronveaux and D.Lambert (Fac.Univ.N.D. de la Paix, Namur, 1991),
   pp.1-29.}
   \ref{Weeks2}{Weeks,Jeffrey \cqg{23}{2006}{6971}.}
   \ref{LandW}{Lachi\`eze-Rey,M. and Weeks,Jeffrey, {\it Orbifold construction of
   the modes on the Poincar\'e dodecahedral space}, arXiv:0801.4232.}
   \ref{Cayley4}{Cayley,A. \qjpam{58}{1879}{280}.}
   \ref{JMS}{Jari\'c,M.V., Michel,L. and Sharp,R.T. {\it J.Physique}
   {\bf 45} (1984) 1. }
   \ref{AandB}{Altmann,S.L. and Bradley,C.J.  {\it Phil. Trans. Roy. Soc. Lond.}
   {\bf 255} (1963) 199.}
   \ref{CandP}{Cummins,C.J. and Patera,J. \jmp{29}{1988}{1736}.}
   \ref{Sloane}{Sloane,N.J.A. \amm{84}{1977}{82}.}
   \ref{Gordan2}{Gordan,P. \ma{12}{1877}{147}.}
   \ref{DandSh}{Desmier,P.E. and Sharp,R.T. \jmp{20}{1979}{74}.}
   \ref{Kramer}{Kramer,P., \jpa{38}{2005}{3517}.}
   \ref{Klein2}{Klein, F.\ma{9}{1875}{183}.}
   \ref{Hodgkinson}{Hodgkinson,J. \jlms{10}{1935}{221}.}
   \ref{ZandD}{Zheng,Y. and Doerschuk, P.C. {\it Acta Cryst.} {\bf A52}
   (1996) 221.}
   \ref{EPM}{Elcoro,L., Perez--Mato,J.M. and Madariaga,G.
   {\it Acta Cryst.} {\bf A50} (1994) 182.}
    \ref{PSW2}{Prandl,W., Schiebel,P. and Wulf,K.
   {\it Acta Cryst.} {\bf A52} (1999) 171.}
    \ref{FCD}{Fan,P--D., Chen,J--Q. and Draayer,J.P.
   {\it Acta Cryst.} {\bf A55} (1999) 871.}
   \ref{FCD2}{Fan,P--D., Chen,J--Q. and Draayer,J.P.
   {\it Acta Cryst.} {\bf A55} (1999) 1049.}
   \ref{Honl}{H\"onl,H. \zfp{89}{1934}{244}.}
   \ref{PSW}{Patera,J., Sharp,R.T. and Winternitz,P. \jmp{19}{1978}{2362}.}
   \ref{LandH}{Lohe,M.A. and Hurst,C.A. \jmp{12}{1971}{1882}.}
   \ref{RandSA}{Ronveaux,A. and Saint-Aubin,Y. \jmp{24}{1983}{1037}.}
   \ref{JandDeV}{Jonker,J.E. and De Vries,E. \npa{105}{1967}{621}.}
   \ref{Rowe}{Rowe, E.G.Peter. \jmp{19}{1978}{1962}.}
   \ref{KNR}{Kibler,M., N\'egadi,T. and Ronveaux,A. {\it The Kustaanheimo-Stiefel
   transformation and certain special functions} \lnm{1171}{1985}{497}.}
   \ref{GLP}{Gilkey,P.B., Leahy,J.V. and Park,J-H, \jpa{29}{1996}{5645}.}
   \ref{Kohler}{K\"ohler,K.: Equivariant Reidemeister torsion on
   symmetric spaces. Math.Ann. {\bf 307}, 57-69 (1997)}
   \ref{Kohler2}{K\"ohler,K.: Equivariant analytic torsion on ${\bf P^nC}$.
   Math.Ann.{\bf 297}, 553-565 (1993) }
   \ref{Kohler3}{K\"ohler,K.: Holomorphic analytic torsion on Hermitian
   symmetric spaces. J.Reine Angew.Math. {\bf 460}, 93-116 (1995)}
   \ref{Zagierzf}{Zagier,D. {\it Zetafunktionen und Quadratische
   K\"orper}, (Springer--Verlag, Berlin, 1981).}
   \ref{Stek}{Stekholschkik,R. {\it Notes on Coxeter transformations and the McKay
   correspondence.} (Springer, Berlin, 2008).}
   \ref{Pesce}{Pesce,H. \cmh {71}{1996}{243}.}
   \ref{Pesce2}{Pesce,H. {\it Contemp. Math} {\bf 173} (1994) 231.}
   \ref{Sutton}{Sutton,C.J. {\it Equivariant isospectrality
   and isospectral deformations on spherical orbifolds}, ArXiv:math/0608567.}
   \ref{Sunada}{Sunada,T. \aom{121}{1985}{169}.}
   \ref{GoandM}{Gornet,R, and McGowan,J. {\it J.Comp. and Math.}
   {\bf 9} (2006) 270.}
   \ref{Suter}{Suter,R. {\it Manusc.Math.} {\bf 122} (2007) 1-21.}
   \ref{Lomont}{Lomont,J.S. {\it Applications of finite groups} (Academic
   Press, New York, 1959).}
   \ref{DandCh2}{Dowker,J.S. and Chang,Peter {\it Analytic torsion on
   spherical factors and tessellations}, arXiv:math.DG/0904.0744 .}
   \ref{Mackey}{Mackey,G. {\it Induced representations}
   (Benjamin, New York, 1968).}
   \ref{Koca}{Koca, {\it Turkish J.Physics}.}
   \ref{Brylinski}{Brylinski, J-L., {\it A correspondence dual to McKay's}
    ArXiv alg-geom/9612003.}
   \ref{Rossmann}{Rossman,W. {\it McKay's correspondence
   and characters of finite subgroups of\break SU(2)} {\it Progress in Math.}
      Birkhauser  (to appear) .}
   \ref{JandL}{James, G. and Liebeck, M. {\it Representations and
   characters of groups} (CUP, Cambridge, 2001).}
   \ref{IandR}{Ito,Y. and Reid,M. {\it The Mckay correspondence for finite
   subgroups of SL(3,C)} Higher dimensional varieties, (Trento 1994),
   221-240, (Berlin, de Gruyter 1996).}
   \ref{BandF}{Bauer,W. and Furutani, K. {\it J.Geom. and Phys.} {\bf
   58} (2008) 64.}
   \ref{Luck}{L\"uck,W. \jdg{37}{1993}{263}.}
   \ref{LandR}{Lott,J. and Rothenberg,M. \jdg{34}{1991}{431}.}
   \ref{DoandKi} {Dowker.J.S. and Kirsten, K. {\it Analysis and Appl.}
   {\bf 3} (2005) 45.}
   \ref{dowtess1}{Dowker,J.S. \cqg{23}{2006}{1}.}
   \ref{dowtess2}{Dowker,J.S. {\it J.Geom. and Phys.} {\bf 57} (2007) 1505.}
   \ref{MHS}{De Melo,T., Hartmann,L. and Spreafico,M. {\it Reidemeister
   Torsion and analytic torsion of discs}, ArXiv:0811.3196.}
   \ref{Vertman}{Vertman, B. {\it Analytic Torsion of a  bounded
   generalized cone}, ArXiv:0808.0449.}
   \ref{WandY} {Weng,L. and You,Y., {\it Int.J. of Math.}{\bf 7} (1996)
   109.}
   \ref{ScandT}{Schwartz, A.S. and Tyupkin,Yu.S. \np{242}{1984}{436}.}
   \ref{AAR}{Andrews, G.E., Askey,R. and Roy,R. {\it Special functions}
  (CUP, Cambridge, 1999).}
   \ref{Tsuchiya}{Tsuchiya, N.: R-torsion and analytic torsion for spherical
   Clifford-Klein manifolds.: J. Fac.Sci., Tokyo Univ. Sect.1 A, Math.
   {\bf 23}, 289-295 (1976).}
   \ref{Tsuchiya2}{Tsuchiya, N. J. Fac.Sci., Tokyo Univ. Sect.1 A, Math.
   {\bf 23}, 289-295 (1976).}
  \ref{Lerch}{Lerch,M. \am{11}{1887}{19}.}
  \ref{Lerch2}{Lerch,M. \am{29}{1905}{333}.}
  \ref{TandS}{Threlfall, W. and Seifert, H. \ma{104}{1930}{1}.}
  \ref{RandS}{Ray, D.B., and Singer, I. \aim{7}{1971}{145}.}
  \ref{RandS2}{Ray, D.B., and Singer, I. {\it Proc.Symp.Pure Math.}
  {\bf 23} (1973) 167.}
  \ref{Jensen}{Jensen,J.L.W.V. \aom{17}{1915-1916}{124}.}
  \ref{Rosenberg}{Rosenberg, S. {\it The Laplacian on a Riemannian Manifold}
  (CUP, Cambridge, 1997).}
  \ref{Nando2}{Nash, C. and O'Connor, D-J. {\it Int.J.Mod.Phys.}
  {\bf A10} (1995) 1779.}
  \ref{Fock}{Fock,V. \zfp{98}{1935}{145}.}
  \ref{Levy}{Levy,M. \prs {204}{1950}{145}.}
  \ref{Schwinger2}{Schwinger,J. \jmp{5}{1964}{1606}.}
  \ref{Muller}{M\"uller, \lnm{}{}{}.}
  \ref{VMK}{Varshalovich.}
  \ref{DandWo}{Dowker,J.S. and Wolski, A. \prA{46}{1992}{6417}.}
  \ref{Zeitlin1}{Zeitlin,V. {\it Physica D} {\bf 49} (1991).  }
  \ref{Zeitlin0}{Zeitlin,V. {\it Nonlinear World} Ed by
   V.Baryakhtar {\it et al},  Vol.I p.717,  (World Scientific, Singapore, 1989).}
  \ref{Zeitlin2}{Zeitlin,V. \prl{93}{2004}{264501}. }
  \ref{Zeitlin3}{Zeitlin,V. \pla{339}{2005}{316}. }
  \ref{Groenewold}{Groenewold, H.J. {\it Physica} {\bf 12} (1946) 405.}
  \ref{Cohen}{Cohen, L. \jmp{7}{1966}{781}.}
  \ref{AandW}{Argawal G.S. and Wolf, E. \prD{2}{1970}{2161,2187,2206}.}
  \ref{Jantzen}{Jantzen,R.T. \jmp{19}{1978}{1163}.}
  \ref{Moses2}{Moses,H.E. \aop{42}{1967}{343}.}
  \ref{Carmeli}{Carmeli,M. \jmp{9}{1968}{1987}.}
  \ref{SHS}{Siemans,M., Hancock,J. and Siminovitch,D. {\it Solid State
  Nuclear Magnetic Resonance} {\bf 31}(2007)35.}
 \ref{Dowk}{Dowker,J.S. \prD{28}{1983}{3013}.}
 \ref{Heine}{Heine, E. {\it Handbuch der Kugelfunctionen}
  (G.Reimer, Berlin. 1878, 1881).}
  \ref{Pockels}{Pockels, F. {\it \"Uber die Differentialgleichung $\De
  u+k^2u=0$} (Teubner, Leipzig. 1891).}
  \ref{Hamermesh}{Hamermesh, M., {\it Group Theory} (Addison--Wesley,
  Reading. 1962).}
  \ref{Racah}{Racah, G. {\it Group Theory and Spectroscopy}
  (Princeton Lecture Notes, 1951). }
  \ref{Gourdin}{Gourdin, M. {\it Basics of Lie Groups} (Editions
  Fronti\'eres, Gif sur Yvette. 1982.)}
  \ref{Clifford}{Clifford, W.K. \plms{2}{1866}{116}.}
  \ref{Story2}{Story, W.E. \plms{23}{1892}{265}.}
  \ref{Story}{Story, W.E. \ma{41}{1893}{469}.}
  \ref{Poole}{Poole, E.G.C. \plms{33}{1932}{435}.}
  \ref{Dickson}{Dickson, L.E. {\it Algebraic Invariants} (Wiley, N.Y.
  1915).}
  \ref{Dickson2}{Dickson, L.E. {\it Modern Algebraic Theories}
  (Sanborn and Co., Boston. 1926).}
  \ref{Hilbert2}{Hilbert, D. {\it Theory of algebraic invariants} (C.U.P.,
  Cambridge. 1993).}
  \ref{Olver}{Olver, P.J. {\it Classical Invariant Theory} (C.U.P., Cambridge.
  1999.)}
  \ref{AST}{A\v{s}erova, R.M., Smirnov, J.F. and Tolsto\v{i}, V.N. {\it
  Teoret. Mat. Fyz.} {\bf 8} (1971) 255.}
  \ref{AandS}{A\v{s}erova, R.M., Smirnov, J.F. \np{4}{1968}{399}.}
  \ref{Shapiro}{Shapiro, J. \jmp{6}{1965}{1680}.}
  \ref{Shapiro2}{Shapiro, J.Y. \jmp{14}{1973}{1262}.}
  \ref{NandS}{Noz, M.E. and Shapiro, J.Y. \np{51}{1973}{309}.}
  \ref{Cayley2}{Cayley, A. {\it Phil. Trans. Roy. Soc. Lond.}
  {\bf 144} (1854) 244.}
  \ref{Cayley3}{Cayley, A. {\it Phil. Trans. Roy. Soc. Lond.}
  {\bf 146} (1856) 101.}
  \ref{Wigner}{Wigner, E.P. {\it Gruppentheorie} (Vieweg, Braunschweig. 1931).}
  \ref{Sharp}{Sharp, R.T. \ajop{28}{1960}{116}.}
  \ref{Laporte}{Laporte, O. {\it Z. f. Naturf.} {\bf 3a} (1948) 447.}
  \ref{Lowdin}{L\"owdin, P-O. \rmp{36}{1964}{966}.}
  \ref{Ansari}{Ansari, S.M.R. {\it Fort. d. Phys.} {\bf 15} (1967) 707.}
  \ref{SSJR}{Samal, P.K., Saha, R., Jain, P. and Ralston, J.P. {\it
  Testing Isotropy of Cosmic Microwave Background Radiation},
  astro-ph/0708.2816.}
  \ref{Lachieze}{Lachi\'eze-Rey, M. {\it Harmonic projection and
  multipole Vectors}. astro- \break ph/0409081.}
  \ref{CHS}{Copi, C.J., Huterer, D. and Starkman, G.D.
  \prD{70}{2003}{043515}.}
  \ref{Jaric}{Jari\'c, J.P. {\it Int. J. Eng. Sci.} {\bf 41} (2003) 2123.}
  \ref{RandD}{Roche, J.A. and Dowker, J.S. \jpa{1}{1968}{527}.}
  \ref{KandW}{Katz, G. and Weeks, J.R. \prD{70}{2004}{063527}.}
  \ref{Waerden}{van der Waerden, B.L. {\it Die Gruppen-theoretische
  Methode in der Quantenmechanik} (Springer, Berlin. 1932).}
  \ref{EMOT}{Erdelyi, A., Magnus, W., Oberhettinger, F. and Tricomi, F.G. {
  \it Higher Transcendental Functions} Vol.1 (McGraw-Hill, N.Y. 1953).}
  \ref{Dowzilch}{Dowker, J.S. {\it Proc. Phys. Soc.} {\bf 91} (1967) 28.}
  \ref{DandD}{Dowker, J.S. and Dowker, Y.P. {\it Proc. Phys. Soc.}
  {\bf 87} (1966) 65.}
  \ref{DandD2}{Dowker, J.S. and Dowker, Y.P. \prs{}{}{}.}
  \ref{Dowk3}{Dowker,J.S. \cqg{7}{1990}{1241}.}
  \ref{Dowk5}{Dowker,J.S. \cqg{7}{1990}{2353}.}
  \ref{CoandH}{Courant, R. and Hilbert, D. {\it Methoden der
  Mathematischen Physik} vol.1 \break (Springer, Berlin. 1931).}
  \ref{Applequist}{Applequist, J. \jpa{22}{1989}{4303}.}
  \ref{Torruella}{Torruella, \jmp{16}{1975}{1637}.}
  \ref{Weinberg}{Weinberg, S.W. \pr{133}{1964}{B1318}.}
  \ref{Meyerw}{Meyer, W.F. {\it Apolarit\"at und rationale Curven}
  (Fues, T\"ubingen. 1883.) }
  \ref{Ostrowski}{Ostrowski, A. {\it Jahrsb. Deutsch. Math. Verein.} {\bf
  33} (1923) 245.}
  \ref{Kramers}{Kramers, H.A. {\it Grundlagen der Quantenmechanik}, (Akad.
  Verlag., Leipzig, 1938).}
  \ref{ZandZ}{Zou, W.-N. and Zheng, Q.-S. \prs{459}{2003}{527}.}
  \ref{Weeks1}{Weeks, J.R. {\it Maxwell's multipole vectors
  and the CMB}.  astro-ph/0412231.}
  \ref{Corson}{Corson, E.M. {\it Tensors, Spinors and Relativistic Wave
  Equations} (Blackie, London. 1950).}
  \ref{Rosanes}{Rosanes, J. \jram{76}{1873}{312}.}
  \ref{Salmon}{Salmon, G. {\it Lessons Introductory to the Modern Higher
  Algebra} 3rd. edn. \break (Hodges,  Dublin. 1876.)}
  \ref{Milnew}{Milne, W.P. {\it Homogeneous Coordinates} (Arnold. London. 1910).}
  \ref{Niven}{Niven, W.D. {\it Phil. Trans. Roy. Soc.} {\bf 170} (1879) 393.}
  \ref{Scott}{Scott, C.A. {\it An Introductory Account of
  Certain Modern Ideas and Methods in Plane Analytical Geometry,}
  (MacMillan, N.Y. 1896).}
  \ref{Bargmann}{Bargmann, V. \rmp{34}{1962}{300}.}
  \ref{Maxwell}{Maxwell, J.C. {\it A Treatise on Electricity and
  Magnetism} 2nd. edn. (Clarendon Press, Oxford. 1882).}
  \ref{BandL}{Biedenharn, L.C. and Louck, J.D.
  {\it Angular Momentum in Quantum Physics} (Addison-Wesley, Reading. 1981).}
  \ref{Weylqm}{Weyl, H. {\it The Theory of Groups and Quantum Mechanics}
  (Methuen, London. 1931).}
  \ref{Robson}{Robson, A. {\it An Introduction to Analytical Geometry} Vol I
  (C.U.P., Cambridge. 1940.)}
  \ref{Sommerville}{Sommerville, D.M.Y. {\it Analytical Conics} 3rd. edn.
   (Bell, London. 1933).}
  \ref{Coolidge}{Coolidge, J.L. {\it A Treatise on Algebraic Plane Curves}
  (Clarendon Press, Oxford. 1931).}
  \ref{SandK}{Semple, G. and Kneebone. G.T. {\it Algebraic Projective
  Geometry} (Clarendon Press, Oxford. 1952).}
  \ref{AandC}{Abdesselam A., and Chipalkatti, J. {\it The Higher
  Transvectants are redundant}, arXiv:0801.1533 [math.AG] 2008.}
  \ref{Elliott}{Elliott, E.B. {\it The Algebra of Quantics} 2nd edn.
  (Clarendon Press, Oxford. 1913).}
  \ref{Elliott2}{Elliott, E.B. \qjpam{48}{1917}{372}.}
  \ref{Howe}{Howe, R. \tams{313}{1989}{539}.}
  \ref{Clebsch}{Clebsch, A. \jram{60}{1862}{343}.}
  \ref{Prasad}{Prasad, G. \ma{72}{1912}{136}.}
  \ref{Dougall}{Dougall, J. \pems{32}{1913}{30}.}
  \ref{Penrose}{Penrose, R. \aop{10}{1960}{171}.}
  \ref{Penrose2}{Penrose, R. \prs{273}{1965}{171}.}
  \ref{Burnside}{Burnside, W.S. \qjm{10}{1870}{211}. }
  \ref{Lindemann}{Lindemann, F. \ma{23} {1884}{111}.}
  \ref{Backus}{Backus, G. {\it Rev. Geophys. Space Phys.} {\bf 8} (1970) 633.}
  \ref{Baerheim}{Baerheim, R. {\it Q.J. Mech. appl. Math.} {\bf 51} (1998) 73.}
  \ref{Lense}{Lense, J. {\it Kugelfunktionen} (Akad.Verlag, Leipzig. 1950).}
  \ref{Littlewood}{Littlewood, D.E. \plms{50}{1948}{349}.}
  \ref{Fierz}{Fierz, M. {\it Helv. Phys. Acta} {\bf 12} (1938) 3.}
  \ref{Williams}{Williams, D.N. {\it Lectures in Theoretical Physics} Vol. VII,
  (Univ.Colorado Press, Boulder. 1965).}
  \ref{Dennis}{Dennis, M. \jpa{37}{2004}{9487}.}
  \ref{Pirani}{Pirani, F. {\it Brandeis Lecture Notes on
  General Relativity,} edited by S. Deser and K. Ford. (Brandeis, Mass. 1964).}
  \ref{Sturm}{Sturm, R. \jram{86}{1878}{116}.}
  \ref{Schlesinger}{Schlesinger, O. \ma{22}{1883}{521}.}
  \ref{Askwith}{Askwith, E.H. {\it Analytical Geometry of the Conic
  Sections} (A.\&C. Black, London. 1908).}
  \ref{Todd}{Todd, J.A. {\it Projective and Analytical Geometry}.
  (Pitman, London. 1946).}
  \ref{Glenn}{Glenn. O.E. {\it Theory of Invariants} (Ginn \& Co, N.Y. 1915).}
  \ref{DowkandG}{Dowker, J.S. and Goldstone, M. \prs{303}{1968}{381}.}
  \ref{Turnbull}{Turnbull, H.A. {\it The Theory of Determinants,
  Matrices and Invariants} 3rd. edn. (Dover, N.Y. 1960).}
  \ref{MacMillan}{MacMillan, W.D. {\it The Theory of the Potential}
  (McGraw-Hill, N.Y. 1930).}
   \ref{Hobson}{Hobson, E.W. {\it The Theory of Spherical
   and Ellipsoidal Harmonics} (C.U.P., Cambridge. 1931).}
  \ref{Hobson1}{Hobson, E.W. \plms {24}{1892}{55}.}
  \ref{GandY}{Grace, J.H. and Young, A. {\it The Algebra of Invariants}
  (C.U.P., Cambridge, 1903).}
  \ref{FandR}{Fano, U. and Racah, G. {\it Irreducible Tensorial Sets}
  (Academic Press, N.Y. 1959).}
  \ref{TandT}{Thomson, W. and Tait, P.G. {\it Treatise on Natural Philosophy}
   (Clarendon Press, Oxford. 1867).}
  \ref{Brinkman}{Brinkman, H.C. {\it Applications of spinor invariants in
atomic physics}, North Holland, Amsterdam 1956.}
  \ref{Kramers1}{Kramers, H.A. {\it Proc. Roy. Soc. Amst.} {\bf 33} (1930) 953.}
  \ref{DandP2}{Dowker,J.S. and Pettengill,D.F. \jpa{7}{1974}{1527}}
  \ref{Dowk1}{Dowker,J.S. \jpa{}{}{45}.}
  \ref{Dowk2}{Dowker,J.S. \aop{71}{1972}{577}}
  \ref{DandA}{Dowker,J.S. and Apps, J.S. \cqg{15}{1998}{1121}.}
  \ref{Weil}{Weil,A., {\it Elliptic functions according to Eisenstein
  and Kronecker}, Springer, Berlin, 1976.}
  \ref{Ling}{Ling,C-H. {\it SIAM J.Math.Anal.} {\bf5} (1974) 551.}
  \ref{Ling2}{Ling,C-H. {\it J.Math.Anal.Appl.}(1988).}
 \ref{BMO}{Brevik,I., Milton,K.A. and Odintsov, S.D. \aop{302}{2002}{120}.}
 \ref{KandL}{Kutasov,D. and Larsen,F. {\it JHEP} 0101 (2001) 1.}
 \ref{KPS}{Klemm,D., Petkou,A.C. and Siopsis {\it Entropy
 bounds, monoticity properties and scaling in CFT's}. hep-th/0101076.}
 \ref{DandC}{Dowker,J.S. and Critchley,R. \prD{15}{1976}{1484}.}
 \ref{AandD}{Al'taie, M.B. and Dowker, J.S. \prD{18}{1978}{3557}.}
 \ref{Dow1}{Dowker,J.S. \prD{37}{1988}{558}.}
 \ref{Dow30}{Dowker,J.S. \prD{28}{1983}{3013}.}
 \ref{DandK}{Dowker,J.S. and Kennedy,G. \jpa{}{1978}{895}.}
 \ref{Dow2}{Dowker,J.S. \cqg{1}{1984}{359}.}
 \ref{DandKi}{Dowker,J.S. and Kirsten, K. {\it Comm. in Anal. and Geom.
 }{\bf7} (1999) 641.}
 \ref{DandKe}{Dowker,J.S. and Kennedy,G.\jpa{11}{1978}{895}.}
 \ref{Gibbons}{Gibbons,G.W. \pl{60A}{1977}{385}.}
 \ref{Cardy}{Cardy,J.L. \np{366}{1991}{403}.}
 \ref{ChandD}{Chang,P. and Dowker,J.S. \np{395}{1993}{407}.}
 \ref{DandC2}{Dowker,J.S. and Critchley,R. \prD{13}{1976}{224}.}
 \ref{Camporesi}{Camporesi,R. \prp{196}{1990}{1}.}
 \ref{BandM}{Brown,L.S. and Maclay,G.J. \pr{184}{1969}{1272}.}
 \ref{CandD}{Candelas,P. and Dowker,J.S. \prD{19}{1979}{2902}.}
 \ref{Unwin1}{Unwin,S.D. Thesis. University of Manchester. 1979.}
 \ref{Unwin2}{Unwin,S.D. \jpa{13}{1980}{313}.}
 \ref{DandB}{Dowker,J.S. and Banach,R. \jpa{11}{1978}{2255}.}
 \ref{Obhukov}{Obhukov,Yu.N. \pl{109B}{1982}{195}.}
 \ref{Kennedy}{Kennedy,G. \prD{23}{1981}{2884}.}
 \ref{CandT}{Copeland,E. and Toms,D.J. \np {255}{1985}{201}.}
  \ref{CandT2}{Copeland,E. and Toms,D.J. \cqg {3}{1986}{431}.}
 \ref{ELV}{Elizalde,E., Lygren, M. and Vassilevich,
 D.V. \jmp {37}{1996}{3105}.}
 \ref{Malurkar}{Malurkar,S.L. {\it J.Ind.Math.Soc} {\bf16} (1925/26) 130.}
 \ref{Glaisher}{Glaisher,J.W.L. {\it Messenger of Math.} {\bf18} (1889) 1.}
  \ref{Anderson}{Anderson,A. \prD{37}{1988}{536}.}
 \ref{CandA}{Cappelli,A. and D'Appollonio,\pl{487B}{2000}{87}.}
 \ref{Wot}{Wotzasek,C. \jpa{23}{1990}{1627}.}
 \ref{RandT}{Ravndal,F. and Tollesen,D. \prD{40}{1989}{4191}.}
 \ref{SandT}{Santos,F.C. and Tort,A.C. \pl{482B}{2000}{323}.}
 \ref{FandO}{Fukushima,K. and Ohta,K. {\it Physica} {\bf A299} (2001) 455.}
 \ref{GandP}{Gibbons,G.W. and Perry,M. \prs{358}{1978}{467}.}
 \ref{Dow4}{Dowker,J.S..}
  \ref{Rad}{Rademacher,H. {\it Topics in analytic number theory,}
Springer-Verlag,  Berlin,1973.}
  \ref{Halphen}{Halphen,G.-H. {\it Trait\'e des Fonctions Elliptiques},
  Vol 1, Gauthier-Villars, Paris, 1886.}
  \ref{CandW}{Cahn,R.S. and Wolf,J.A. {\it Comm.Mat.Helv.} {\bf 51}
  (1976) 1.}
  \ref{Berndt}{Berndt,B.C. \rmjm{7}{1977}{147}.}
  \ref{Hurwitz}{Hurwitz,A. \ma{18}{1881}{528}.}
  \ref{Hurwitz2}{Hurwitz,A. {\it Mathematische Werke} Vol.I. Basel,
  Birkhauser, 1932.}
  \ref{Berndt2}{Berndt,B.C. \jram{303/304}{1978}{332}.}
  \ref{RandA}{Rao,M.B. and Ayyar,M.V. \jims{15}{1923/24}{150}.}
  \ref{Hardy}{Hardy,G.H. \jlms{3}{1928}{238}.}
  \ref{TandM}{Tannery,J. and Molk,J. {\it Fonctions Elliptiques},
   Gauthier-Villars, Paris, 1893--1902.}
  \ref{schwarz}{Schwarz,H.-A. {\it Formeln und
  Lehrs\"atzen zum Gebrauche..},Springer 1893.(The first edition was 1885.)
  The French translation by Henri Pad\'e is {\it Formules et Propositions
  pour L'Emploi...},Gauthier-Villars, Paris, 1894}
  \ref{Hancock}{Hancock,H. {\it Theory of elliptic functions}, Vol I.
   Wiley, New York 1910.}
  \ref{watson}{Watson,G.N. \jlms{3}{1928}{216}.}
  \ref{MandO}{Magnus,W. and Oberhettinger,F. {\it Formeln und S\"atze},
  Springer-Verlag, Berlin 1948.}
  \ref{Klein}{Klein,F. {\it Lectures on the Icosohedron}
  (Methuen, London. 1913).}
  \ref{AandL}{Appell,P. and Lacour,E. {\it Fonctions Elliptiques},
  Gauthier-Villars,
  Paris. 1897.}
  \ref{HandC}{Hurwitz,A. and Courant,C. {\it Allgemeine Funktionentheorie},
  Springer,
  Berlin. 1922.}
  \ref{WandW}{Whittaker,E.T. and Watson,G.N. {\it Modern analysis},
  Cambridge. 1927.}
  \ref{SandC}{Selberg,A. and Chowla,S. \jram{227}{1967}{86}. }
  \ref{zucker}{Zucker,I.J. {\it Math.Proc.Camb.Phil.Soc} {\bf 82 }(1977)
  111.}
  \ref{glasser}{Glasser,M.L. {\it Maths.of Comp.} {\bf 25} (1971) 533.}
  \ref{GandW}{Glasser, M.L. and Wood,V.E. {\it Maths of Comp.} {\bf 25}
  (1971)
  535.}
  \ref{greenhill}{Greenhill,A,G. {\it The Applications of Elliptic
  Functions}, MacMillan. London, 1892.}
  \ref{Weierstrass}{Weierstrass,K. {\it J.f.Mathematik (Crelle)}
{\bf 52} (1856) 346.}
  \ref{Weierstrass2}{Weierstrass,K. {\it Mathematische Werke} Vol.I,p.1,
  Mayer u. M\"uller, Berlin, 1894.}
  \ref{Fricke}{Fricke,R. {\it Die Elliptische Funktionen und Ihre Anwendungen},
    Teubner, Leipzig. 1915, 1922.}
  \ref{Konig}{K\"onigsberger,L. {\it Vorlesungen \"uber die Theorie der
 Elliptischen Funktionen},  \break Teubner, Leipzig, 1874.}
  \ref{Milne}{Milne,S.C. {\it The Ramanujan Journal} {\bf 6} (2002) 7-149.}
  \ref{Schlomilch}{Schl\"omilch,O. {\it Ber. Verh. K. Sachs. Gesell. Wiss.
  Leipzig}  {\bf 29} (1877) 101-105; {\it Compendium der h\"oheren
  Analysis}, Bd.II, 3rd Edn, Vieweg, Brunswick, 1878.}
  \ref{BandB}{Briot,C. and Bouquet,C. {\it Th\`eorie des Fonctions
  Elliptiques}, Gauthier-Villars, Paris, 1875.}
  \ref{Dumont}{Dumont,D. \aim {41}{1981}{1}.}
  \ref{Andre}{Andr\'e,D. {\it Ann.\'Ecole Normale Superior} {\bf 6} (1877)
  265;
  {\it J.Math.Pures et Appl.} {\bf 5} (1878) 31.}
  \ref{Raman}{Ramanujan,S. {\it Trans.Camb.Phil.Soc.} {\bf 22} (1916) 159;
 {\it Collected Papers}, Cambridge, 1927}
  \ref{Weber}{Weber,H.M. {\it Lehrbuch der Algebra} Bd.III, Vieweg,
  Brunswick 190  3.}
  \ref{Weber2}{Weber,H.M. {\it Elliptische Funktionen und algebraische
  Zahlen},
  Vieweg, Brunswick 1891.}
  \ref{ZandR}{Zucker,I.J. and Robertson,M.M.
  {\it Math.Proc.Camb.Phil.Soc} {\bf 95 }(1984) 5.}
  \ref{JandZ1}{Joyce,G.S. and Zucker,I.J.
  {\it Math.Proc.Camb.Phil.Soc} {\bf 109 }(1991) 257.}
  \ref{JandZ2}{Zucker,I.J. and Joyce.G.S.
  {\it Math.Proc.Camb.Phil.Soc} {\bf 131 }(2001) 309.}
  \ref{zucker2}{Zucker,I.J. {\it SIAM J.Math.Anal.} {\bf 10} (1979) 192.}
  \ref{BandZ}{Borwein,J.M. and Zucker,I.J. {\it IMA J.Math.Anal.} {\bf 12}
  (1992) 519.}
  \ref{Cox}{Cox,D.A. {\it Primes of the form $x^2+n\,y^2$}, Wiley,
  New York, 1989.}
  \ref{BandCh}{Berndt,B.C. and Chan,H.H. {\it Mathematika} {\bf42} (1995)
  278.}
  \ref{EandT}{Elizalde,R. and Tort.hep-th/}
  \ref{KandS}{Kiyek,K. and Schmidt,H. {\it Arch.Math.} {\bf 18} (1967) 438.}
  \ref{Oshima}{Oshima,K. \prD{46}{1992}{4765}.}
  \ref{greenhill2}{Greenhill,A.G. \plms{19} {1888} {301}.}
  \ref{Russell}{Russell,R. \plms{19} {1888} {91}.}
  \ref{BandB}{Borwein,J.M. and Borwein,P.B. {\it Pi and the AGM}, Wiley,
  New York, 1998.}
  \ref{Resnikoff}{Resnikoff,H.L. \tams{124}{1966}{334}.}
  \ref{vandp}{Van der Pol, B. {\it Indag.Math.} {\bf18} (1951) 261,272.}
  \ref{Rankin}{Rankin,R.A. {\it Modular forms} C.U.P. Cambridge}
  \ref{Rankin2}{Rankin,R.A. {\it Proc. Roy.Soc. Edin.} {\bf76 A} (1976) 107.}
  \ref{Skoruppa}{Skoruppa,N-P. {\it J.of Number Th.} {\bf43} (1993) 68 .}
  \ref{Down}{Dowker.J.S. {\it Nucl.Phys.}B (Proc.Suppl) ({\bf 104})(2002)153;
  also Dowker,J.S. hep-th/ 0007129.}
  \ref{Eichler}{Eichler,M. \mz {67}{1957}{267}.}
  \ref{Zagier}{Zagier,D. \invm{104}{1991}{449}.}
  \ref{Lang}{Lang,S. {\it Modular Forms}, Springer, Berlin, 1976.}
  \ref{Kosh}{Koshliakov,N.S. {\it Mess.of Math.} {\bf 58} (1928) 1.}
  \ref{BandH}{Bodendiek, R. and Halbritter,U. \amsh{38}{1972}{147}.}
  \ref{Smart}{Smart,L.R., \pgma{14}{1973}{1}.}
  \ref{Grosswald}{Grosswald,E. {\it Acta. Arith.} {\bf 21} (1972) 25.}
  \ref{Kata}{Katayama,K. {\it Acta Arith.} {\bf 22} (1973) 149.}
  \ref{Ogg}{Ogg,A. {\it Modular forms and Dirichlet series} (Benjamin,
  New York,
   1969).}
  \ref{Bol}{Bol,G. \amsh{16}{1949}{1}.}
  \ref{Epstein}{Epstein,P. \ma{56}{1903}{615}.}
  \ref{Petersson}{Petersson.}
  \ref{Serre}{Serre,J-P. {\it A Course in Arithmetic}, Springer,
  New York, 1973.}
  \ref{Schoenberg}{Schoenberg,B., {\it Elliptic Modular Functions},
  Springer, Berlin, 1974.}
  \ref{Apostol}{Apostol,T.M. \dmj {17}{1950}{147}.}
  \ref{Ogg2}{Ogg,A. {\it Lecture Notes in Math.} {\bf 320} (1973) 1.}
  \ref{Knopp}{Knopp,M.I. \dmj {45}{1978}{47}.}
  \ref{Knopp2}{Knopp,M.I. \invm {}{1994}{361}.}
  \ref{LandZ}{Lewis,J. and Zagier,D. \aom{153}{2001}{191}.}
  \ref{DandK1}{Dowker,J.S. and Kirsten,K. {\it Elliptic functions and
  temperature inversion symmetry on spheres} hep-th/.}
  \ref{HandK}{Husseini and Knopp.}
  \ref{Kober}{Kober,H. \mz{39}{1934-5}{609}.}
  \ref{HandL}{Hardy,G.H. and Littlewood, \am{41}{1917}{119}.}
  \ref{Watson}{Watson,G.N. \qjm{2}{1931}{300}.}
  \ref{SandC2}{Chowla,S. and Selberg,A. {\it Proc.Nat.Acad.} {\bf 35}
  (1949) 371.}
  \ref{Landau}{Landau, E. {\it Lehre von der Verteilung der Primzahlen},
  (Teubner, Leipzig, 1909).}
  \ref{Berndt4}{Berndt,B.C. \tams {146}{1969}{323}.}
  \ref{Berndt3}{Berndt,B.C. \tams {}{}{}.}
  \ref{Bochner}{Bochner,S. \aom{53}{1951}{332}.}
  \ref{Weil2}{Weil,A.\ma{168}{1967}{}.}
  \ref{CandN}{Chandrasekharan,K. and Narasimhan,R. \aom{74}{1961}{1}.}
  \ref{Rankin3}{Rankin,R.A. {} {} ().}
  \ref{Berndt6}{Berndt,B.C. {\it Trans.Edin.Math.Soc}.}
  \ref{Elizalde}{Elizalde,E. {\it Ten Physical Applications of Spectral
  Zeta Function Theory}, \break (Springer, Berlin, 1995).}
  \ref{Allen}{Allen,B., Folacci,A. and Gibbons,G.W. \pl{189}{1987}{304}.}
  \ref{Krazer}{Krazer}
  \ref{Elizalde3}{Elizalde,E. {\it J.Comp.and Appl. Math.} {\bf 118}
  (2000) 125.}
  \ref{Elizalde2}{Elizalde,E., Odintsov.S.D, Romeo, A. and Bytsenko,
  A.A and
  Zerbini,S.
  {\it Zeta function regularisation}, (World Scientific, Singapore,
  1994).}
  \ref{Eisenstein}{Eisenstein}
  \ref{Hecke}{Hecke,E. \ma{112}{1936}{664}.}
  \ref{Hecke2}{Hecke,E. \ma{112}{1918}{398}.}
  \ref{Terras}{Terras,A. {\it Harmonic analysis on Symmetric Spaces} (Springer,
  New York, 1985).}
  \ref{BandG}{Bateman,P.T. and Grosswald,E. {\it Acta Arith.} {\bf 9}
  (1964) 365.}
  \ref{Deuring}{Deuring,M. \aom{38}{1937}{585}.}
  \ref{Mordell}{Mordell,J. \prs{}{}{}.}
  \ref{GandZ}{Glasser,M.L. and Zucker, {}.}
  \ref{Landau2}{Landau,E. \jram{}{1903}{64}.}
  \ref{Kirsten1}{Kirsten,K. \jmp{35}{1994}{459}.}
  \ref{Sommer}{Sommer,J. {\it Vorlesungen \"uber Zahlentheorie}
  (1907,Teubner,Leipzig).
  French edition 1913 .}
  \ref{Reid}{Reid,L.W. {\it Theory of Algebraic Numbers},
  (1910,MacMillan,New York).}
  \ref{Milnor}{Milnor, J. {\it Is the Universe simply--connected?},
  IAS, Princeton, 1978.}
  \ref{Milnor2}{Milnor, J. \ajm{79}{1957}{623}.}
  \ref{Opechowski}{Opechowski,W. {\it Physica} {\bf 7} (1940) 552.}
  \ref{Bethe}{Bethe, H.A. \zfp{3}{1929}{133}.}
  \ref{LandL}{Landau, L.D. and Lishitz, E.M. {\it Quantum
  Mechanics} (Pergamon Press, London, 1958).}
  \ref{GPR}{Gibbons, G.W., Pope, C. and R\"omer, H., \np{157}{1979}{377}.}
  \ref{Jadhav}{Jadhav,S.P. PhD Thesis, University of Manchester 1990.}
  \ref{DandJ}{Dowker,J.S. and Jadhav, S. \prD{39}{1989}{1196}.}
  \ref{CandM}{Coxeter, H.S.M. and Moser, W.O.J. {\it Generators and
  relations of finite groups} (Springer. Berlin. 1957).}
  \ref{Coxeter2}{Coxeter, H.S.M. {\it Regular Complex Polytopes},
   (Cambridge University Press, \break Cambridge, 1975).}
  \ref{Coxeter}{Coxeter, H.S.M. {\it Regular Polytopes}.}
  \ref{Stiefel}{Stiefel, E., J.Research NBS {\bf 48} (1952) 424.}
  \ref{BandS}{Brink, D.M. and Satchler, G.R. {\it Angular momentum theory}.
  (Clarendon Press, Oxford. 1962.).}
  %\ref{Racah1}
  \ref{Rose}{Rose}
  \ref{Schwinger}{Schwinger, J. {\it On Angular Momentum}
  in {\it Quantum Theory of Angular Momentum} edited by
  Biedenharn,L.C. and van Dam, H. (Academic Press, N.Y. 1965).}
  
  \ref{Ray}{Ray,D.B. \aim{4}{1970}{109}.}
  \ref{Ikeda}{Ikeda,A. {\it Kodai Math.J.} {\bf 18} (1995) 57.}
  \ref{Kennedy}{Kennedy,G. \prD{23}{1981}{2884}.}
  \ref{Ellis}{Ellis,G.F.R. {\it General Relativity} {\bf2} (1971) 7.}
  \ref{Dow8}{Dowker,J.S. \cqg{20}{2003}{L105}.}
  \ref{IandY}{Ikeda, A and Yamamoto, Y. \ojm {16}{1979}{447}.}
  \ref{BandI}{Bander,M. and Itzykson,C. \rmp{18}{1966}{2}.}
  \ref{Schulman}{Schulman, L.S. \pr{176}{1968}{1558}.}
  \ref{Bar1}{B\"ar,C. {\it Arch.d.Math.}{\bf 59} (1992) 65.}
  \ref{Bar2}{B\"ar,C. {\it Geom. and Func. Anal.} {\bf 6} (1996) 899.}
  \ref{Vilenkin}{Vilenkin, N.J. {\it Special functions},
  (Am.Math.Soc., Providence, 1968).}
  \ref{Talman}{Talman, J.D. {\it Special functions} (Benjamin,N.Y.,1968).}
  \ref{Miller}{Miller, W. {\it Symmetry groups and their applications}
  (Wiley, N.Y., 1972).}
  \ref{Dow3}{Dowker,J.S. \cmp{162}{1994}{633}.}
  \ref{Cheeger}{Cheeger, J. \jdg {18}{1983}{575}.}
  \ref{Cheeger2}{Cheeger, J. \aom {109}{1979}{259}.}
  \ref{Dow6}{Dowker,J.S. \jmp{30}{1989}{770}.}
  \ref{Dow20}{Dowker,J.S. \jmp{35}{1994}{6076}.}
  \ref{Dowjmp}{Dowker,J.S. \jmp{35}{1994}{4989}.}
  \ref{Dow21}{Dowker,J.S. {\it Heat kernels and polytopes} in {\it
   Heat Kernel Techniques and Quantum Gravity}, ed. by S.A.Fulling,
   Discourses in Mathematics and its Applications, No.4, Dept.
   Maths., Texas A\&M University, College Station, Texas, 1995.}
  \ref{Dow9}{Dowker,J.S. \jmp{42}{2001}{1501}.}
  \ref{Dow7}{Dowker,J.S. \jpa{25}{1992}{2641}.}
  \ref{Warner}{Warner.N.P. \prs{383}{1982}{379}.}
  \ref{Wolf}{Wolf, J.A. {\it Spaces of constant curvature},
  (McGraw--Hill,N.Y., 1967).}
  \ref{Meyer}{Meyer,B. \cjm{6}{1954}{135}.}
  \ref{BandB}{B\'erard,P. and Besson,G. {\it Ann. Inst. Four.} {\bf 30}
  (1980) 237.}
  \ref{PandM}{P\'{o}lya,G. and Meyer,B. \cras{228}{1948}{28}.}
  \ref{Springer}{Springer, T.A. Lecture Notes in Math. vol 585 (Springer,
  Berlin,1977).}
  \ref{SeandT}{Threlfall, H. and Seifert, W. \ma{104}{1930}{1}.}
  \ref{Hopf}{Hopf,H. \ma{95}{1925}{313}. }
  \ref{Dow}{Dowker,J.S. \jpa{5}{1972}{936}.}
  \ref{LLL}{Lehoucq,R., Lachi\'eze-Rey,M. and Luminet, J.--P. {\it
  Astron.Astrophys.} {\bf 313} (1996) 339.}
  \ref{LaandL}{Lachi\'eze-Rey,M. and Luminet, J.--P.
  \prp{254}{1995}{135}.}
  \ref{Schwarzschild}{Schwarzschild, K., {\it Vierteljahrschrift der
  Ast.Ges.} {\bf 35} (1900) 337.}
  \ref{Starkman}{Starkman,G.D. \cqg{15}{1998}{2529}.}
  \ref{LWUGL}{Lehoucq,R., Weeks,J.R., Uzan,J.P., Gausman, E. and
  Luminet, J.--P. \cqg{19}{2002}{4683}.}
  \ref{Dow10}{Dowker,J.S. \prD{28}{1983}{3013}.}
  \ref{BandD}{Banach, R. and Dowker, J.S. \jpa{12}{1979}{2527}.}
  \ref{Jadhav2}{Jadhav,S. \prD{43}{1991}{2656}.}
  \ref{Gilkey}{Gilkey,P.B. {\it Invariance theory,the heat equation and
  the Atiyah--Singer Index theorem} (CRC Press, Boca Raton, 1994).}
  \ref{BandY}{Berndt,B.C. and Yeap,B.P. {\it Adv. Appl. Math.}
  {\bf29} (2002) 358.}
  \ref{HandR}{Hanson,A.J. and R\"omer,H. \pl{80B}{1978}{58}.}
  \ref{Hill}{Hill,M.J.M. {\it Trans.Camb.Phil.Soc.} {\bf 13} (1883) 36.}
  \ref{Cayley}{Cayley,A. {\it Quart.Math.J.} {\bf 7} (1866) 304.}
  \ref{Seade}{Seade,J.A. {\it Anal.Inst.Mat.Univ.Nac.Aut\'on
  M\'exico} {\bf 21} (1981) 129.}
  \ref{CM}{Cisneros--Molina,J.L. {\it Geom.Dedicata} {\bf84} (2001)
  \ref{Goette1}{Goette,S. \jram {526} {2000} 181.}
  207.}
  \ref{NandO}{Nash,C. and O'Connor,D--J, \jmp {36}{1995}{1462}.}
  \ref{Dows}{Dowker,J.S. \aop{71}{1972}{577}; Dowker,J.S. and Pettengill,D.F.
  \jpa{7}{1974}{1527}; J.S.Dowker in {\it Quantum Gravity}, edited by
  S. C. Christensen (Hilger,Bristol,1984)}
  \ref{Jadhav2}{Jadhav,S.P. \prD{43}{1991}{2656}.}
  \ref{Dow11}{Dowker,J.S. \cqg{21}{2004}4247.}
  \ref{Dow12}{Dowker,J.S. \cqg{21}{2004}4977.}
  \ref{Dow13}{Dowker,J.S. \jpa{38}{2005}1049.}
  \ref{Zagier}{Zagier,D. \ma{202}{1973}{149}}
  \ref{RandG}{Rademacher, H. and Grosswald,E. {\it Dedekind Sums},
  (Carus, MAA, 1972).}
  \ref{Berndt7}{Berndt,B, \aim{23}{1977}{285}.}
  \ref{HKMM}{Harvey,J.A., Kutasov,D., Martinec,E.J. and Moore,G.
  {\it Localised Tachyons and RG Flows}, hep-th/0111154.}
  \ref{Beck}{Beck,M., {\it Dedekind Cotangent Sums}, {\it Acta Arithmetica}
  {\bf 109} (2003) 109-139 ; math.NT/0112077.}
  \ref{McInnes}{McInnes,B. {\it APS instability and the topology of the brane
  world}, hep-th/0401035.}
  \ref{BHS}{Brevik,I, Herikstad,R. and Skriudalen,S. {\it Entropy Bound for the
  TM Electromagnetic Field in the Half Einstein Universe}; hep-th/0508123.}
  \ref{BandO}{Brevik,I. and Owe,C.  \prD{55}{4689}{1997}.}
  \ref{Kenn}{Kennedy,G. Thesis. University of Manchester 1978.}
  \ref{KandU}{Kennedy,G. and Unwin S. \jpa{12}{L253}{1980}.}
  \ref{BandO1}{Bayin,S.S.and Ozcan,M.
  \prD{48}{2806}{1993}; \prD{49}{5313}{1994}.}
  \ref{Chang}{Chang, P., {\it Quantum Field Theory on Regular Polytopes}.
   Thesis. University of Manchester, 1993.}
  \ref{Barnesa}{Barnes,E.W. {\it Trans. Camb. Phil. Soc.} {\bf 19} (1903) 374.}
  \ref{Barnesb}{Barnes,E.W. {\it Trans. Camb. Phil. Soc.}
  {\bf 19} (1903) 426.}
  \ref{Stanley1}{Stanley,R.P. \joa {49Hilf}{1977}{134}.}
  \ref{Stanley2}{Stanley,R.P. {\it Enumerative Combinatorics} vols.1,2
  (C.U.P., Cambridge, 1997, 1999.}
  \ref{Stanley}{Stanley,R.P. \bams {1}{1979}{475}.}
  \ref{Hurley}{Hurley,A.C. \pcps {47}{1951}{51}.}
  \ref{IandK}{Iwasaki,I. and Katase,K. {\it Proc.Japan Acad. Ser} {\bf A55}
  (1979) 141.}
  \ref{IandT}{Ikeda,A. and Taniguchi,Y. {\it Osaka J. Math.} {\bf 15} (1978)
  515.}
  \ref{GandM}{Gallot,S. and Meyer,D. \jmpa{54}{1975}{259}.}
  \ref{Flatto}{Flatto,L. {\it Enseign. Math.} {\bf 24} (1978) 237.}
  \ref{OandT}{Orlik,P and Terao,H. {\it Arrangements of Hyperplanes},
  Grundlehren der Math. Wiss. {\bf 300}, (Springer--Verlag, 1992).}
  \ref{Shepler}{Shepler,A.V. \joa{220}{1999}{314}.}
  \ref{SandT}{Solomon,L. and Terao,H. \cmh {73}{1998}{237}.}
  \ref{Vass}{Vassilevich, D.V. \plb {348}{1995}39.}
  \ref{Vass2}{Vassilevich, D.V. \jmp {36}{1995}3174.}
  \ref{CandH}{Camporesi,R. and Higuchi,A. {\it J.Geom. and Physics}
  {\bf 15} (1994) 57.}
  \ref{Solomon2}{Solomon,L. \tams{113}{1964}{274}.}
  \ref{Solomon}{Solomon,L. {\it Nagoya Math. J.} {\bf 22} (1963) 57.}
  \ref{Obukhov}{Obukhov,Yu.N. \pl{109B}{1982}{195}.}
  \ref{BGH}{Bernasconi,F., Graf,G.M. and Hasler,D. {\it The heat kernel
  expansion for the electromagnetic field in a cavity}; math-ph/0302035.}
  \ref{Baltes}{Baltes,H.P. \prA {6}{1972}{2252}.}
  \ref{BaandH}{Baltes.H.P and Hilf,E.R. {\it Spectra of Finite Systems}
  (Bibliographisches Institut, Mannheim, 1976).}
  \ref{Ray}{Ray,D.B. \aim{4}{1970}{109}.}
  \ref{Hirzebruch}{Hirzebruch,F. {\it Topological methods in algebraic
  geometry} (Springer-- Verlag,\break  Berlin, 1978). }
  \ref{BBG}{Bla\v{z}i\'c,N., Bokan,N. and Gilkey, P.B. {\it Ind.J.Pure and
  Appl.Math.} {\bf 23} (1992) 103.}
  \ref{WandWi}{Weck,N. and Witsch,K.J. {\it Math.Meth.Appl.Sci.} {\bf 17}
  (1994) 1017.}
  \ref{Norlund}{N\"orlund,N.E. \am{43}{1922}{121}.}
   \ref{Norlund1}{N\"orlund,N.E. {\it Differenzenrechnung} (Springer--Verlag, 1924, Berlin.)}
  \ref{Duff}{Duff,G.F.D. \aom{56}{1952}{115}.}
  \ref{DandS}{Duff,G.F.D. and Spencer,D.C. \aom{45}{1951}{128}.}
  \ref{BGM}{Berger, M., Gauduchon, P. and Mazet, E. {\it Lect.Notes.Math.}
  {\bf 194} (1971) 1. }
  \ref{Patodi}{Patodi,V.K. \jdg{5}{1971}{233}.}
  \ref{GandS}{G\"unther,P. and Schimming,R. \jdg{12}{1977}{599}.}
  \ref{MandS}{McKean,H.P. and Singer,I.M. \jdg{1}{1967}{43}.}
  \ref{Conner}{Conner,P.E. {\it Mem.Am.Math.Soc.} {\bf 20} (1956).}
  \ref{Gilkey2}{Gilkey,P.B. \aim {15}{1975}{334}.}
  \ref{MandP}{Moss,I.G. and Poletti,S.J. \plb{333}{1994}{326}.}
  \ref{BKD}{Bordag,M., Kirsten,K. and Dowker,J.S. \cmp{182}{1996}{371}.}
  \ref{RandO}{Rubin,M.A. and Ordonez,C. \jmp{25}{1984}{2888}.}
  \ref{BaandD}{Balian,R. and Duplantier,B. \aop {112}{1978}{165}.}
  \ref{Kennedy2}{Kennedy,G. \aop{138}{1982}{353}.}
  \ref{DandKi2}{Dowker,J.S. and Kirsten, K. {\it Analysis and Appl.}
 {\bf 3} (2005) 45.}
  \ref{Dow40}{Dowker,J.S. \cqg{23}{2006}{1}.}
  \ref{BandHe}{Br\"uning,J. and Heintze,E. {\it Duke Math.J.} {\bf 51} (1984)
   959.}
  \ref{Dowl}{Dowker,J.S. {\it Functional determinants on M\"obius corners};
    Proceedings, `Quantum field theory under
    the influence of external conditions', 111-121,Leipzig 1995.}
  \ref{Dowqg}{Dowker,J.S. in {\it Quantum Gravity}, edited by
  S. C. Christensen (Hilger, Bristol, 1984).}
  \ref{Dowit}{Dowker,J.S. \jpa{11}{1978}{347}.}
  \ref{Kane}{Kane,R. {\it Reflection Groups and Invariant Theory} (Springer,
  New York, 2001).}
  \ref{Sturmfels}{Sturmfels,B. {\it Algorithms in Invariant Theory}
  (Springer, Vienna, 1993).}
  \ref{Bourbaki}{Bourbaki,N. {\it Groupes et Alg\`ebres de Lie}  Chap.III, IV
  (Hermann, Paris, 1968).}
  \ref{SandTy}{Schwarz,A.S. and Tyupkin, Yu.S. \np{242}{1984}{436}.}
  \ref{Reuter}{Reuter,M. \prD{37}{1988}{1456}.}
  \ref{EGH}{Eguchi,T. Gilkey,P.B. and Hanson,A.J. \prp{66}{1980}{213}.}
  \ref{DandCh}{Dowker,J.S. and Chang,Peter, \prD{46}{1992}{3458}.}
  \ref{APS}{Atiyah M., Patodi and Singer,I.\mpcps{77}{1975}{43}.}
  \ref{Donnelly}{Donnelly.H. {\it Indiana U. Math.J.} {\bf 27} (1978) 889.}
  \ref{Katase}{Katase,K. {\it Proc.Jap.Acad.} {\bf 57} (1981) 233.}
  \ref{Gilkey3}{Gilkey,P.B.\invm{76}{1984}{309}.}
  \ref{Degeratu}{Degeratu.A. {\it Eta--Invariants and Molien Series for
  Unimodular Groups}, Thesis MIT, 2001.}
  \ref{Seeley}{Seeley,R. \ijmp {A\bf18}{2003}{2197}.}
  \ref{Seeley2}{Seeley,R. .}
  \ref{melrose}{Melrose}
  \ref{DandW}{Douglas,R.G. and Wojciekowski,K.P. \cmp{142}{1991}{139}.}
  \ref{Dai}{Dai,X. \tams{354}{2001}{107}.}
\end{putreferences}

\bye